\documentclass{article}
\usepackage[utf8]{inputenc}
\usepackage[margin=1in]{geometry}
\usepackage{empheq}
\usepackage{float}
\usepackage{amssymb}
\usepackage{stmaryrd}
\usepackage{comment}
\usepackage{dsfont}
\usepackage{parskip}
\usepackage{authblk}
\usepackage{verbatim}
\usepackage{mdframed}
\usepackage{hyperref}

\newtheorem{prop}{Proposition}

\title{A Divergence-Conforming Hybridized Discontinuous Galerkin Method for the Incompressible Magnetohydrodynamics Equations}
\author[1]{Thad A. Gleason\footnote{thad.gleason@tamu.edu}}
\author[2]{Eric L. Peters\footnote{epeters1@ball.com}}
\author[3]{John A. Evans\footnote{john.a.evans@colorado.edu}}
\affil[1]{\scriptsize\textit{Department of Mathematics, Texas A\&M University}}
\affil[2]{\scriptsize\textit{Ball Aerospace \& Technologies}}
\affil[3]{\scriptsize\textit{Ann and H.J. Smead Aerospace Engineering Sciences, University of Colorado}}
\date{}

\begin{document}

\maketitle
\begin{abstract}
    \noindent We introduce a new hybridized discontinuous Galerkin method for the incompressible magnetohydrodynamics equations. If particular velocity, pressure, magnetic field, and magnetic pressure spaces are employed for both element and trace solution fields, we arrive at an energy stable method which returns pointwise divergence-free velocity fields and magnetic fields and properly balances linear momentum. We discretize in time using a second-order-in-time generalized-$\alpha$ method, and we present a block iterative method for solving the resulting nonlinear system of equations at each time step. We numerically examine the effectiveness of our method using a manufactured solution and observe our method yields optimal convergence rates in the $L_2$ norm for the velocity field, pressure field, magnetic field, and magnetic pressure field. We further find our method is pressure robust. We then apply our method to a selection of benchmark problems and numerically confirm our method is energy stable.
\end{abstract}


\section{Introduction}
\label{intro}

Magnetohydrodynamics (MHD) is the study of the behavior of electrically conducting fluids \cite{mhd-book}. This field is important when studying fluids such as plasmas, liquid metals, or electrolytic solutions. Some examples of when the MHD equations would apply is in the formation of stars, behavior of cosmic dust, or plasma confinement in a nuclear fusion reactor. The fundamental equations involved in MHD are a combination of the Navier-Stokes equations from fluid dynamics and Maxwell's equations from electromagnetism. Note that in this paper we focus only on incompressible MHD. For this we require the Mach number be sufficiently small enough to assume the density of the fluid remains constant. Even given this restriction, incompressible MHD is useful in studying problems such as coolant in a liquid metal fast-fission reactor \cite{fast-fission} and the plasma confinement and liquid metal coolant of a fusion reactor \cite{tokamak, metal-fusion}.

Some popular methods for solving the MHD equations include finite difference methods (FDM), finite volume methods (FVM), and finite element methods (FEM). In the category of FEM, we typically see a stabilized continuous Galerkin (CG) method \cite{Shadid1, Shadid2, Shadid3, stable-cg1, stable-cg2} or a discontinuous Galerkin (DG) method \cite{local-cdg, local-dg, exact-cdg, mhd-dg} implemented. DG methods combine the advantages of FVM with FEM. DG methods are ideally suited for solving hyperbolic partial differential equations (PDEs) defined over complex geometries, and they are robust in the presence of large solution gradients, including shocks. The major drawback to DG methods is an increased number of degrees of freedom (DOF) when compared to CG methods.

In order to reduce the number of global DOFs associated with DG methods, hybridized discontinuous Galerkin (HDG) methods were introduced. With an HDG method, there are both interior DOFs that reside on the interior of elements and trace DOFs that reside on the boundaries of elements. With the use of static condensation, the interior DOFs can be written in terms of the trace DOFs, allowing the interior DOFs to be removed completely from the system of equations. This significantly reduces the number of global DOFs. This method was first introduced for symmetric elliptic problems in \cite{OG}, and it has since been extended to other PDEs such as the advection-diffusion equation \cite{ad1, ad2, ad3}, the incompressible and compressible Navier-Stokes equations \cite{ Wells,ns-e-hdg}, the incompressible Reynolds Averaged Navier-Stokes equations with the Spalart-Allmaras model \cite{Eric}, and the PDEs governing incompressible and compressible MHD \cite{mhd-hdg2,mhd-hdg1}.

While HDG methods have been developed for the incompressible MHD equations, to the best of the authors’ knowledge, no HDG method presented in the literature to date yields both pointwise divergence-free velocity fields and pointwise divergence-free magnetic fields. However, divergence-conforming methods that exactly preserve divergence-free constraints harbor a number of advantageous properties over methods that satisfy divergence-free constraints in an approximate manner \cite{siam-review, div-free}. For instance, CG and DG methods for the incompressible Navier-Stokes equations that yield pointwise divergence-free velocity fields are energy stable and properly balance momentum \cite{evans2013isogeometric}, and they further yield velocity field approximations whose error does not depend on the pressure field \cite{linke2016pressure}. This latter property is referred to as pressure robustness in the literature. A recent paper further proved velocity field error estimates for such methods that are independent of both the pressure and the Reynolds number \cite{schroeder2018towards}. Divergence-conforming CG and DG methods for the incompressible MHD equations exhibit similar properties \cite{fu2019explicit}, though error estimates independent of the Reynolds number and magnetic Reynolds number do not exist yet for such discretizations. The above inspires us to construct an HDG method for the incompressible MHD equations that yields pointwise divergence-free velocity and magnetic fields by extending a divergence-conforming HDG method for the incompressible Navier-Stokes equations \cite{Wells}. This requires special treatment of the terms that couple the velocity and magnetic fields to arrive at a method that is energy stable. It should be noted that a previously presented HDG method for the incompressible MHD equations does yield pointwise divergence-free velocity fields \cite{qiu2020mixed}, but it yields only discretely divergence-free magnetic fields.

An outline of this paper is as follows. In the next section, we recall the strong form of the incompressible MHD equations. In Section 3, we present our semi-discrete divergence-conforming HDG method for the incompressible MHD equations and show this method is consistent. In Section 4, we prove our method pointwise conserves mass, is pointwise absent of magnetic monopoles, and globally conserves momentum. In Section 5, we prove an energy stability result for our method. In Section 6, we adopt a second-order time integration scheme for our semi-discrete method, and we present a block iterative scheme for solving the nonlinear algebraic equations at each time step. In Section 7, we demonstrate how static condensation may be employed to reduce the size of the linear systems associated with our block iterative scheme. In Section 8, we assess the spatial accuracy, temporal accuracy, pressure robustness, and energy stability of our method using a suite of numerical experiments. In Section 9, we present concluding remarks.

\section{The Strong Form of the Incompressible MHD Equations}

We begin by presenting the strong form of the incompressible MHD equations. Let $\Omega$ denote a bounded domain in $\mathbb{R}^d$, where $d$ denotes the dimensionality of the domain. We assume $d=2$ or $d=3$. Let $\partial\Omega$ denote the boundary of $\Omega$ with outward unit normal $\mathbf{n}$, and let $\partial\Omega$ be partitioned into a Dirichlet boundary $\Gamma_D$ and a Neumann boundary $\Gamma_N$ such that $\partial\Omega=\overline{\Gamma_D\cup\Gamma_N}$ and $\Gamma_D\cap\Gamma_N= \emptyset$. Let $\mathbf{u}$ denote the velocity field, $p$ denote the pressure field, $\mathbf{B}$ denote the magnetic field, and $r$ denote the magnetic pressure field. We also denote the permeability of free space as $\mu_0$. The strong form of the incompressible MHD equations is then as follows:\\

\begin{mdframed}
Find $\mathbf{u}:\bar{\Omega}\times[0,\infty)\rightarrow\mathbb{R}^d$,  $p:\Omega\times(0,\infty)\rightarrow\mathbb{R}$, $\mathbf{B}:\bar{\Omega}\times[0,\infty)\rightarrow\mathbb{R}^d$, and $r:\Omega\times(0,\infty)\rightarrow\mathbb{R}$ such that:

\textbf{Conservation of Momentum}
\begin{equation}
\frac{\partial\mathbf{u}}{\partial t}+\nabla\cdot\left(\mathbf{u}\otimes\mathbf{u}\right)+\frac{1}{\rho}\nabla p-\nabla\cdot\left(2\nu\nabla^s\mathbf{u}\right)-\nabla\cdot\left(\frac{1}{\rho \mu_0}\left(\mathbf{B}\otimes\mathbf{B}-\frac{1}{2}|\mathbf{B}|^2\mathds{1}\right)\right)=\mathbf{f}_v\qquad\text{in }\Omega\times(0,\infty)
\label{eq:strong_mom}
\end{equation}
\textbf{Conservation of Mass}
\begin{equation}
\nabla\cdot\mathbf{u}=0\qquad\text{in }\Omega\times(0,\infty)
\label{eq:strong_mass}
\end{equation}
\textbf{The Magnetic Induction Equation}
\begin{equation}
\frac{\partial \mathbf{B}}{\partial t}+\nabla\cdot(\mathbf{u}\otimes\mathbf{B})-\nabla\cdot(\mathbf{B}\otimes\mathbf{u})-\nabla\cdot\left(2\frac{\eta}{\mu_0}\nabla^a\mathbf{B}\right)+\nabla r = \mathbf{f}_m \qquad \text{in }\Omega\times(0,\infty)
\label{eq:strong_induction}
\end{equation}
\textbf{Gauss's Law for Magnetism}
\begin{equation}
\nabla\cdot\mathbf{B}=0\qquad\text{in }\Omega\times(0,\infty)
\label{eq:strong_gauss}
\end{equation}
\textbf{Boundary Conditions}
\begin{align}
    \mathbf{u}&=\mathbf{g}_v && \text{on }\Gamma_D\times(0,\infty) \label{eq:strong_BC_1}\\
    \mathbf{n}\cdot\left(-\frac{1}{\rho}p\mathds{1}+2\nu\nabla^s\mathbf{u}+\frac{1}{\rho \mu_0}\left(\mathbf{B}\otimes\mathbf{B}-\frac{1}{2}|\mathbf{B}|^2\mathds{1}\right)\right) - \min(\mathbf{u}\cdot\mathbf{n},0)\mathbf{u}&=\mathbf{h}_v && \text{on }\Gamma_N\times(0,\infty) \label{eq:strong_BC_2}\\
    \mathbf{B}&=\mathbf{g}_m && \text{on }\Gamma_D\times(0,\infty) \label{eq:strong_BC_3}\\
    \mathbf{n}\cdot\left(-r\mathds{1}+2\frac{\eta}{\mu_0} \nabla^a\mathbf{B}\right) - \min(\mathbf{u}\cdot\mathbf{n},0)\mathbf{B}&=\mathbf{h}_m && \text{on }\Gamma_N\times(0,\infty) \label{eq:strong_BC_4}
\end{align}
\textbf{Initial Conditions}
\begin{align}
        \mathbf{u}(:,0)&=\mathbf{u}_0 \qquad \hspace{1.5pt} \text{in }\Omega \label{eq:strong_IC_1}\\
        \mathbf{B}(:,0)&=\mathbf{B}_0 \qquad \text{in }\Omega \label{eq:strong_IC_2}
\end{align}
\end{mdframed}

In the above equations, $\nabla^s*=\frac{1}{2}((\nabla*)+(\nabla*)^T)$ denotes the symmetric gradient operator, $\nabla^a*=\frac{1}{2}((\nabla*)-(\nabla*)^T)$ denotes the antisymmetric gradient operator, $\otimes$ denotes the outer product operator, and $\mathds{1}$ denotes the identity matrix of degree $d$. We take the gradient of a vector $\textbf{a} = \sum_{j=1}^d a_{j} \textbf{e}_j$ to be $\nabla \textbf{a} := \sum_{i=1}^d \sum_{j=1}^d \frac{\partial a_j}{\partial x_i} \textbf{e}_i \otimes \textbf{e}_j$ where $\left\{ \textbf{e}_j \right\}_{j=1}^d$ is the standard basis for $\mathbb{R}^d$, and we take the divergence of a second-order tensor $\textbf{A} = \sum_{i=1}^d \sum_{j=1}^d A_{ij} \textbf{e}_i \otimes \textbf{e}_j$ to be $\nabla \cdot \textbf{A} := \sum_{i=1}^d \sum_{j=1}^d \frac{\partial A_{ij}}{\partial x_i} \textbf{e}_j$. We later adopt the convention $\textbf{C} : \textbf{D} = C_{ij} D_{ij}$ for two second-order tensors $\textbf{C} = \sum_{i=1}^d \sum_{j=1}^d C_{ij} \textbf{e}_i \otimes \textbf{e}_j$ and $\textbf{D} = \sum_{i=1}^d \sum_{j=1}^d D_{ij} \textbf{e}_i \otimes \textbf{e}_j$. We assume constant density $\rho\in\mathbb{R}^+$, variable kinematic viscosity $\nu : \Omega\times(0,\infty)\rightarrow\mathbb{R}^+$, variable resistivity $\eta : \Omega\times(0,\infty)\rightarrow\mathbb{R}^+$, variable body forces $\mathbf{f}_v:\Omega\times(0,\infty)\rightarrow\mathbb{R}^d$ and $\mathbf{f}_m:\Omega\times(0,\infty)\rightarrow\mathbb{R}^d$, variable velocity and magnetic field specifications on the Dirichlet boundary $\mathbf{g}_v:\Gamma_D\times(0,\infty)\rightarrow\mathbb{R}^d$ and $\mathbf{g}_m:\Gamma_D\times(0,\infty)\rightarrow\mathbb{R}^d$,  variable traction specifications on the Neumann boundary $\mathbf{h}_v:\Gamma_N\times(0,\infty)\rightarrow\mathbb{R}^d$ and $\mathbf{h}_m:\Gamma_N\times(0,\infty)\rightarrow\mathbb{R}^d$, and variable initial conditions $\mathbf{u}_0:\Omega\rightarrow\mathbb{R}^d$ and $\mathbf{B}_0:\Omega\rightarrow\mathbb{R}^d$. The magnetic pressure $r$ is not a physical quantity but is instead included in order to enforce the $\nabla\cdot\mathbf{B}=0$ condition in our HDG  method. The magnetic pressure may be interpreted as a Lagrange multiplier associated with the $\nabla\cdot\mathbf{B}=0$ condition and, for a suitable choice of boundary conditions, it is identically zero \cite{ben1999conservative}. Inclusion of a magnetic pressure to enforce the $\nabla\cdot\mathbf{B}=0$ condition is common in both continuous Galerkin \cite{Shadid2,stable-cg2,ben1999conservative,codina2006stabilized} and discontinuous Galerkin \cite{mhd-hdg2,qiu2020mixed,houston2009mixed} finite element methods.

\section{A Semi-Discrete Divergence-Conforming HDG Method for the Incompressible MHD Equations}

\begin{figure}
\centering
\includegraphics[width=0.55\textwidth]{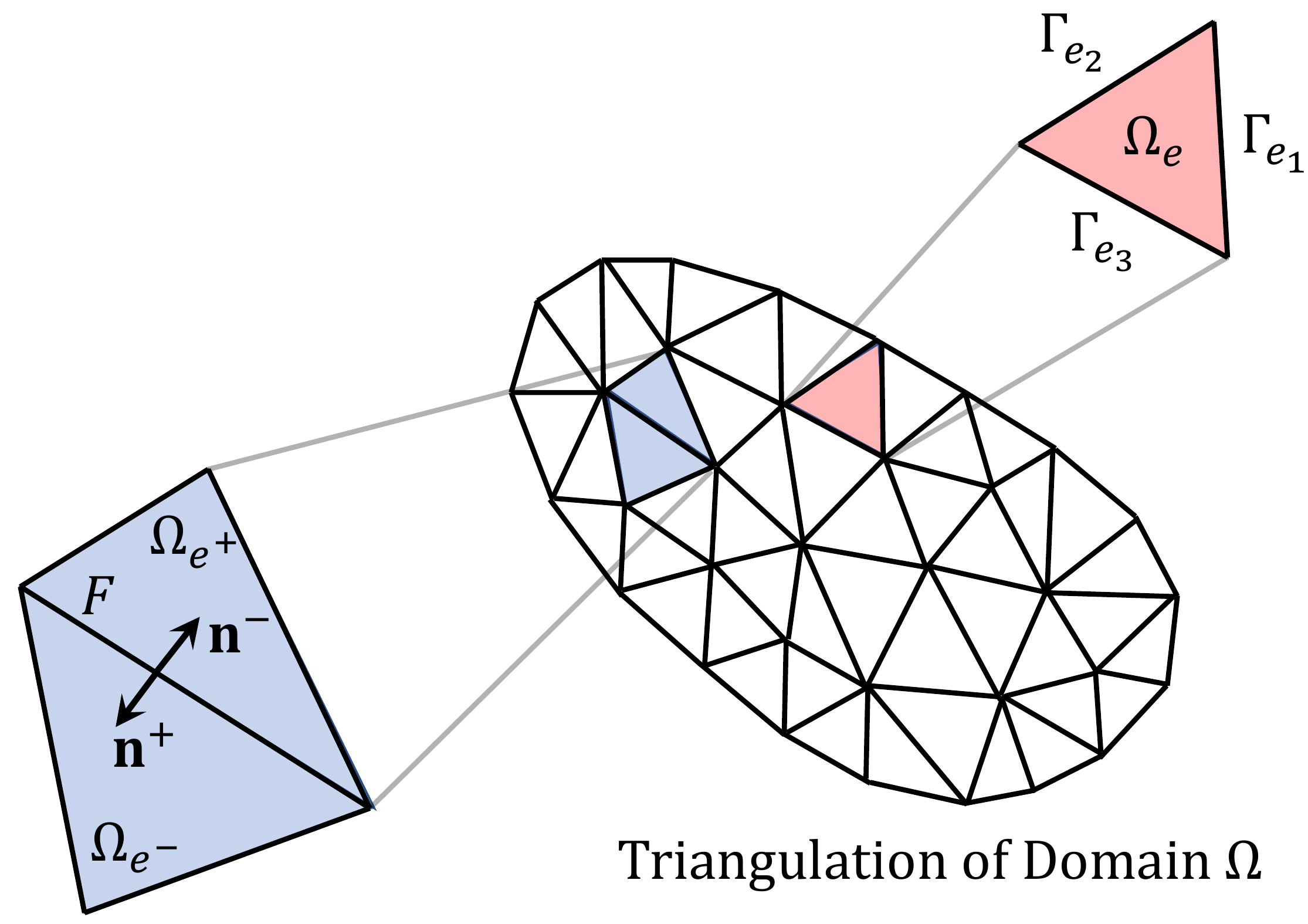}
\caption{Graphical depiction of a triangulation $\mathcal{T}$, the notation associated with an element $\Omega_e \in \mathcal{T}$, and the notation associated with a facet $F \in \mathcal{F}$.}
\label{mesh}
\end{figure}

We are now ready to construct our divergence-conforming HDG method for the incompressible MHD equations. In this section, we discretize in space, and later, we discretize in time. As is typically done with HDG methods, we first introduce a mesh over which the incompressible MHD equations will be discretized. Let $\mathcal{T}=\{\Omega_e\}_{e=1}^\text{nel}$ be a triangulation of the domain such that $\bar{\Omega}= \overline{\cup_{e=1}^\text{nel} \Omega_e}$. The $i^\text{th}$ facet of an element $\Omega_e \in \mathcal{T}$ is denoted by $\Gamma_{e_i}$, and the outward unit normal vector on $\Gamma_{e_i}$ to $\Omega_e$ is denoted by $\mathbf{n}$. Let $\mathcal{F}$ be the set of all facets. We define the mesh skeleton as $\tilde{\Gamma} = \overline{\cup_{F \in \mathcal{F}} F}$. Facets that lie on the boundary of the domain are called boundary facets, while facets that do not lie on the boundary of the domain are called interior facets. We denote the sets of interior and boundary facets as $\mathcal{F}_\textup{int}$ and $\mathcal{F}_\textup{bdy}$ respectively. Each interior facet $F$ is shared by two adjacent elements $\Omega_{e^+}$ and $\Omega_{e^-}$. We denote the outward facing normals on $F$ to $\Omega_{e^+}$ and $\Omega_{e^-}$ as $\mathbf{n}^+$ and $\mathbf{n}^-$ respectively. A visual representation of the above objects in the two-dimensional setting is displayed in Figure \ref{mesh}.

%We denote the $e$th element as $\Omega_e$, its $i$th edge as $\Gamma_{e_i}$, and the outward unit normal vector on $\Gamma_{e_i}$ to $\Omega_e$ as $\mathbf{n}$. {\color{red} John: Issue here. We use bold here and vector notation later. Be consistent. Thad: Thanks for letting me know. I had missed some vectors when I changed everything to bold. I think I have fixed them.}   Let $\mathcal{F}$ be the set of all facets, $F$ such that $\Gamma=\cup_{F\in\mathcal{F}}\bar{F}$. And facet $F$ that lies on the boundary of the domain $\partial\Omega$ is called a boundary facet. A facet $F$ that lies between two adjacent elements, denoted $\Omega_{e^+}$ and $\Omega_{e^-}$ is called an interior facet. We denote the outward facing normals on F to $\Omega_{e^+}$ and $\Omega_{e^-}$ as $\mathbf{n}^+$ and $\mathbf{n}^-$ respectively. We define the jump operator across $F$ as $\llbracket*\rrbracket=*^+\cdot\mathbf{n}^++*^-\cdot\mathbf{n}^-$ A two-dimensional representation of this breakdown can be seen in Figure \ref{mesh}. 

Let $P_k(D)$ denote the space of polynomials of degree $k\geq0$ on a domain $D$. For $k\geq1$, we consider the following finite element spaces
for the velocity field
\begin{align}
V^h &:=\{\mathbf{v}^h\in[L^2(\Omega)]^d:\mathbf{v}^h\in[P_k(\Omega_e)]^d~ \forall\Omega_e\in\mathcal{T}\} \\
\hat V^h &:=\{\mathbf{\hat{v}}^h\in[L^2(\tilde{\Gamma})]^d:\mathbf{\hat{v}}^h\in[P_{k}(F)]^d~ \forall F\in\mathcal{F}\},
\label{Vspace}
\end{align}
the pressure field
\begin{align}
Q^h_v &:=\{q^h_v\in L^2(\Omega):q^h_v\in P_{k-1}(\Omega_e)~\forall\Omega_e\in\mathcal{T}\} \\
\hat Q^h_v &:=\{\hat{q}^h_v\in L^2(\tilde{\Gamma}):\hat{q}^h_v\in P_{k}(F)~\forall F\in\mathcal{F}\},
\label{pspace}
\end{align}
the magnetic field
\begin{align}
W^h &:=\{\mathbf{w}^h\in[L^2(\Omega)]^d:\mathbf{w}^h\in[P_k(\Omega_e)]^d~ \forall\Omega_e\in\mathcal{T}\} \\
\hat W^h &:=\{\mathbf{\hat{w}}^h\in[L^2(\tilde{\Gamma})]^d:\mathbf{\hat{w}}^h\in[P_{k}(F)]^d~ \forall F\in\mathcal{F}\},
\label{Bspace}
\end{align}
and the magnetic pressure field
\begin{equation}
\begin{split}
&Q^h_m:=\{q^h_m\in L^2(\Omega):q^h_m\in P_{k-1}(\Omega_e)~\forall\Omega_e\in\mathcal{T}\} \\
&\hat Q^h_m:=\{\hat{q}^h_m\in L^2(\tilde{\Gamma}):\hat{q}^h_m\in P_{k}(F)~\forall F\in\mathcal{F}\}.
\label{rspace}
\end{split}
\end{equation}
We approximate the velocity field, pressure field, magnetic field, and magnetic pressure field over element interiors using the spaces $V^h$, $Q_v^h$, $W^h$, and $Q_m^h$, and we approximate the velocity field, pressure field, magnetic field, and magnetic pressure field over the mesh skeleton using the spaces $\hat{V}^h$, $\hat{Q}_v^h$, $\hat{W}^h$, and $\hat{Q}_m^h$. We also introduce the spaces
\begin{align}
        \hat V_g^h(t) &:= \{\mathbf{\hat{v}}^h\in\hat V^h : \mathbf{\hat{v}}^h=\mathbf{g}_v(:,t) \text{ on } \Gamma_D\} \\
        \hat V_0^h &:= \{\mathbf{\hat{v}}^h\in\hat V^h : \mathbf{\hat{v}}^h=\mathbf{0} \text{ on } \Gamma_D\} \\
        \hat W_g^h(t) &:= \{\mathbf{\hat{w}}^h\in\hat W^h : \mathbf{\hat{w}}^h=\mathbf{g}_m(:,t) \text{ on } \Gamma_D\} \\
        \hat W_0^h &:= \{\mathbf{\hat{w}}^h\in\hat W^h : \mathbf{\hat{w}}^h=\mathbf{0} \text{ on } \Gamma_D\}.
\end{align}
We employ the spaces $\hat V_g^h(t)$ and $\hat V_0^h$ as velocity field trace trial and test spaces, and we employ the spaces $\hat W_g^h(t)$ and $\hat W_0^h$ as magnetic field trace trial and test spaces.

As vector-valued functions in the spaces $V^h$ and $W^h$ are discontinuous across element boundaries, we introduce an operator to measure the jump of the normal component of these functions across element boundaries. Let $F$ be an interior facet shared by two adjacent elements $\Omega_{e^+}$ and $\Omega_{e^-}$ with outward facing normals $\textbf{n}^+$ and $\textbf{n}^-$. For vector-valued functions $\textbf{y}$ lying in either $V^h$ and $W^h$, we denote the jump of the normal component of $\textbf{y}$ across $F$ as $\llbracket \textbf{y} \rrbracket = \textbf{y}^+ \cdot \textbf{n}^+ + \textbf{y}^- \cdot \textbf{n}^-$ where $\textbf{y}^+ = \textbf{y}|_{\Omega_{e^+}}$ and $\textbf{y}^- = \textbf{y}|_{\Omega_{e^-}}$.

With all of the above notation in hand, we are ready to present our semi-discrete HDG method for the incompressible MHD equations. Our method may be interpreted as an extension of Rhebergen and Wells's divergence-conforming HDG method for the incompressible Navier-Stokes equations \cite{Wells}, wherein advective fluxes are treated using upwinding and diffusive fluxes using the symmetric interior penalty method, to the incompressible MHD equations. We further take special care in discretizing the Lorentz force $-\nabla\cdot\left(\frac{1}{\rho \mu_0}\left(\mathbf{B}\otimes\mathbf{B}-\frac{1}{2}|\mathbf{B}|^2\mathds{1}\right)\right)$ appearing in the conservation of momentum equation and the coupling term $-\nabla\cdot\left(\mathbf{B}\otimes\mathbf{u}\right)$ appearing in the magnetic induction equation to arrive at a method that is energy stable. These two terms are responsible for the transfer of energy between the velocity field and the magnetic field.\\

\begin{mdframed}
Find $(\mathbf{u}^h(t),\mathbf{\hat{u}}^h(t),p^h(t),\hat{p}^h(t),\mathbf{B}^h(t),\mathbf{\hat{B}}^h(t),r^h(t),\hat{r}^h(t)) \in V^h\times\hat{V}^h_g(t)\times Q^h_v\times\hat{Q}^h_v \times W^h\times\hat{W}^h_g(t)\times Q^h_m\times\hat{Q}^h_m$ for each $t \in [0,\infty)$ such that:\\

\textbf{Conservation of Momentum}
\begin{equation}
\begin{split}
\Sigma_e&\int_{\Omega_e}\frac{\partial\mathbf{u}^h}{\partial t}\cdot\mathbf{v}^h-\Sigma_e\int_{\Omega_e}(\mathbf{u}^h\otimes\mathbf{u}^h):\nabla\mathbf{v}^h
-\Sigma_e\int_{\Omega_e}\frac{1}{\rho}p^h\mathds{1}:\nabla\mathbf{v}^h+\Sigma_e\int_{\Omega_e}2\nu\nabla^s\mathbf{u}^h:\nabla^s\mathbf{v}^h\\
+\Sigma_e&\int_{\Omega_e}\frac{1}{\rho \mu_0}\left(\mathbf{B}^h\otimes\mathbf{B}^h-\frac{1}{2}|\mathbf{B}^h|^2\mathds{1}\right):\nabla\mathbf{v}^h
+\Sigma_e\Sigma_i\int_{\Gamma_{e_i}} (1-\lambda) (\mathbf{u}^h\otimes\mathbf{u}^h):(\mathbf{n}\otimes\mathbf{v}^h)\\
+\Sigma_e\Sigma_i&\int_{\Gamma_{e_i}}\lambda (\mathbf{u}^h\otimes\mathbf{\hat{u}}^h):(\mathbf{n}\otimes\mathbf{v}^h)
+\Sigma_e\Sigma_i\int_{\Gamma_{e_i}}\frac{1}{\rho}\hat{p}^h\mathds{1}:(\mathbf{n}\otimes\mathbf{v}^h)\\
-\Sigma_e\Sigma_i&\int_{\Gamma_{e_i}}2\nu\nabla^s\mathbf{u}^h:(\mathbf{n}\otimes\mathbf{v}^h)
+\Sigma_e\Sigma_i\int_{\Gamma_{e_i}}2\frac{C_\textup{pen}}{h_e}\nu\left( \mathbf{n} \otimes (\mathbf{u}^h-\mathbf{\hat{u}}^h) \right) :(\mathbf{n} \otimes \mathbf{v}^h)\\
-\Sigma_e\Sigma_i&\int_{\Gamma_{e_i}}2\nu\left(\mathbf{n} \otimes (\mathbf{u}^h-\mathbf{\hat{u}}^h)\right):\nabla^s\mathbf{v}^h
-\Sigma_e\Sigma_i\int_{\Gamma_{e_i}}\frac{1}{\rho \mu_0}\left(\mathbf{\hat{B}}^h\otimes\mathbf{\hat{B}}^h-\frac{1}{2}|\mathbf{\hat{B}}^h|^2\mathds{1}\right):(\mathbf{n} \otimes \mathbf{v}^h)\\
-\Sigma_e&\int_{\Omega_e}\mathbf{f}_v\cdot\mathbf{v}^h
=0 \hspace{86.5mm} \forall\mathbf{v}^h\in V^h \text{ and } t \in (0,\infty)
\end{split}
\label{p}
\end{equation}

\textbf{Conservation of Momentum Flux}
\begin{equation}
\begin{split}
\Sigma_e\Sigma_i&\int_{\Gamma_{e_i}} (1-\lambda) (\mathbf{u}^h\otimes\mathbf{u}^h):(\mathbf{n}\otimes\mathbf{\hat{v}}^h) +\Sigma_e\Sigma_i\int_{\Gamma_{e_i}}\lambda (\mathbf{u}^h\otimes\mathbf{\hat{u}}^h):(\mathbf{n}\otimes\mathbf{\hat{v}}^h) \\
+\Sigma_e\Sigma_i&\int_{\Gamma_{e_i}}\frac{1}{\rho}\hat{p}^h\mathds{1}:(\mathbf{n}\otimes\mathbf{\hat{v}}^h)
-\Sigma_e\Sigma_i\int_{\Gamma_{e_i}}2\nu\nabla^s\mathbf{u}^h:(\mathbf{n}\otimes\mathbf{\hat{v}}^h) \\
+\Sigma_e\Sigma_i&\int_{\Gamma_{e_i}}2\frac{C_\textup{pen}}{h_e}\nu\left( \mathbf{n} \otimes (\mathbf{u}^h-\mathbf{\hat{u}}^h) \right) :(\mathbf{n} \otimes \mathbf{\hat{v}}^h) -\Sigma_e\Sigma_i\int_{\Gamma_{e_i}}\frac{1}{\rho \mu_0}\left(\mathbf{\hat{B}}^h\otimes\mathbf{\hat{B}}^h-\frac{1}{2}|\mathbf{\hat{B}}^h|^2\mathds{1}\right):(\mathbf{n} \otimes \mathbf{\hat{v}}^h)\\
-\Sigma_e\Sigma_i&\int_{\Gamma_N\cap\Gamma_{e_i}}(1-\lambda)(\mathbf{u}^h\otimes\mathbf{\hat{u}}^h):(\mathbf{n}\otimes\mathbf{\hat{v}}^h)
+\Sigma_e\Sigma_i\int_{\Gamma_N\cap\Gamma_{e_i}}\mathbf{h}_v\cdot\mathbf{\hat{v}}^h
=0\hspace{12mm} \forall\mathbf{\hat{v}}^h\in\hat{V}^h_0 \text{ and } t \in (0,\infty)
    \end{split}
    \label{pc}
\end{equation}

\textbf{Conservation of Mass}
\begin{equation}
    \Sigma_e\int_{\Omega_e}\frac{1}{\rho}\nabla\cdot\mathbf{u}^hq^h_v=0 \hspace{10mm} \forall q^h_v\in Q_v \text{ and } t \in (0,\infty)
\label{m}
\end{equation}

\textbf{Conservation of Mass Flux}
\begin{equation}
    \Sigma_e\Sigma_i\int_{\Gamma_{e_i}}\frac{1}{\rho} \left( \mathbf{u}^h -\mathbf{\hat{u}}^h \right) \cdot\mathbf{n}\hat{q}^h_v
    =0\hspace{10mm} \forall\hat{q}^h_v\in\hat{Q}^h_v \text{ and } t \in (0,\infty)
    \label{mc}
\end{equation}

\textbf{The Magnetic Induction Equation}
\begin{equation}
\begin{split}
\Sigma_e&\int_{\Omega_e}\frac{\partial\mathbf{B}^h}{\partial t}\cdot\mathbf{w}^h - \Sigma_e\int_{\Omega_e} (\mathbf{u}^h\otimes\mathbf{B}^h):\nabla\mathbf{w}^h
-\Sigma_e\int_{\Omega_e} r^h\mathds{1}:\nabla\mathbf{w}^h
+\Sigma_e\int_{\Omega_e} 2\frac{\eta}{\mu_0}\nabla^a\mathbf{B}^h:\nabla^a\mathbf{w}^h\\
-\Sigma_e&\int_{\Omega_e}(\mathbf{B}^h\otimes\mathbf{w}^h):\nabla\mathbf{u}^h
+\Sigma_e\Sigma_i\int_{\Gamma_{e_i}} (1-\lambda) (\mathbf{u}^h\otimes\mathbf{B}^h):(\mathbf{n}\otimes\mathbf{w}^h)\\
+\Sigma_e\Sigma_i&\int_{\Gamma_{e_i}}\lambda (\mathbf{u}^h\otimes \mathbf{\hat{B}}^h):(\mathbf{n}\otimes\mathbf{w}^h)
+\Sigma_e\Sigma_i\int_{\Gamma_{e_i}}\hat r^h\mathds{1}:(\mathbf{n}\otimes\mathbf{w}^h)\\
-\Sigma_e\Sigma_i&\int_{\Gamma_{e_i}} 2\frac{\eta}{\mu_0}\nabla^a\mathbf{B}^h:(\mathbf{n}\otimes\mathbf{w}^h)
+\Sigma_e\Sigma_i\int_{\Gamma_{e_i}} 2\frac{\eta}{\mu_0}\frac{C_\textup{pen}}{h_e}\left(\mathbf{n}\otimes(\mathbf{B}^h-\mathbf{\hat{B}}^h)\right):(\mathbf{n}\otimes\mathbf{w}^h)\\
-\Sigma_e\Sigma_i&\int_{\Gamma_{e_i}} 2\frac{\eta}{\mu_0}\left(\mathbf{n}\otimes(\mathbf{B}^h-\mathbf{\hat{B}}^h)\right):\nabla^a\mathbf{w}^h
-\Sigma_e\int_{\Omega_e}\mathbf{f}_m\cdot\mathbf{w}^h 
= 0 \hspace{21mm} \forall\mathbf{w}^h\in W^h \text{ and } t \in (0,\infty)
\end{split}
\label{B}
\end{equation}

\textbf{Conservation of Magnetic Induction Flux}
\begin{equation}
\begin{split}
\Sigma_e\Sigma_i&\int_{\Gamma_{e_i}} (1-\lambda) (\mathbf{u}^h\otimes\mathbf{B}^h):(\mathbf{n}\otimes\mathbf{\hat{w}}^h)
+\Sigma_e\Sigma_i\int_{\Gamma_{e_i}}\lambda (\mathbf{u}^h\otimes \mathbf{\hat{B}}^h):(\mathbf{n}\otimes\mathbf{\hat{w}}^h)\\
+\Sigma_e\Sigma_i&\int_{\Gamma_{e_i}}\hat r^h\mathds{1}:(\mathbf{n}\otimes\mathbf{\hat{w}}^h)
-\Sigma_e\Sigma_i\int_{\Gamma_{e_i}} 2\frac{\eta}{\mu_0}\nabla^a\mathbf{B}^h:(\mathbf{n}\otimes\mathbf{\hat{w}}^h)\\
+\Sigma_e\Sigma_i&\int_{\Gamma_{e_i}} 2\frac{\eta}{\mu_0}\frac{C_\textup{pen}}{h_e}\left(\mathbf{n}\otimes(\mathbf{B}^h-\mathbf{\hat{B}}^h)\right):(\mathbf{n}\otimes\mathbf{\hat{w}}^h)\\ 
-\Sigma_e\Sigma_i&\int_{\Gamma_{e_i}} \left(\mathbf{\hat{B}}^h\otimes\mathbf{\hat{w}}^h - \frac{1}{2} \left( \mathbf{\hat{B}}^h\cdot\mathbf{\hat{w}}^h \right) \mathds{1} \right):(\mathbf{n}\otimes\left(\mathbf{u}^h-\mathbf{\hat{u}}^h\right))\\
-\Sigma_e\Sigma_i&\int_{\Gamma_N \cap \Gamma_{e_i}} (1-\lambda) (\mathbf{u}^h\otimes\mathbf{\hat{B}}^h):(\mathbf{n}\otimes\mathbf{\hat{w}}^h)+\Sigma_e\Sigma_i\int_{\Gamma_N\cap\Gamma_{e_i}}\mathbf{h}_m\cdot\mathbf{\hat{w}}^h 
= 0 \hspace{7mm} \forall \mathbf{\hat{w}}^h\in \hat{W}^h_0  \text{ and } t \in (0,\infty)
\end{split}
\label{Bc}
\end{equation}

\textbf{Gauss's Law of Magnetism}
\begin{equation}
\Sigma_e\int_{\Omega_e} q^h_m\nabla\cdot\mathbf{B}^h = 0 \hspace{10mm} \forall q^h_m\in Q^h_m  \text{ and } t \in (0,\infty)
\label{mono}
\end{equation}

\textbf{Conservation of Magnetic Flux}
\begin{equation}
\Sigma_e\Sigma_i\int_{\Gamma_{e_i}} \left( \mathbf{B}^h - \mathbf{\hat{B}}^h \right) \cdot\mathbf{n}\hat{q}^h_m = 0 \hspace{10mm} \forall \hat{q}^h_m\in \hat{Q}^h_m  \text{ and } t \in (0,\infty)
\label{monoc}
\end{equation}

\textbf{Velocity Initial Condition}
\begin{equation}
    \Sigma_e\int_{\Omega_e} \left( \mathbf{u}^h(:,0) - \mathbf{u}_0 \right) \cdot \mathbf{v}^h = 0 \hspace{10mm} \forall \textbf{v}^h \in V^h
\label{v_ic}
\end{equation}

\textbf{Magnetic Field Initial Condition}
\begin{equation}
    \Sigma_e\int_{\Omega_e} \left( \mathbf{B}^h(:,0) - \mathbf{B}_0 \right) \cdot \mathbf{w}^h = 0 \hspace{10mm} \forall \textbf{w}^h \in W^h
\label{b_ic}
\end{equation}

\end{mdframed}

In the above equations, on the $i$th facet $\Gamma_{e_i}$ of the $e$th element $\Omega_e$, $\lambda$ is an indicator function that takes the value of one if the facet is an inflow facet and a value of zero if it is an outflow facet, that is,
\begin{equation}
\lambda = 
\begin{cases}
1 & \mathbf{u}^h\cdot\mathbf{n}<0 \\
0 & \mathbf{u}^h\cdot\mathbf{n}\geq 0.
\end{cases}
\end{equation}
Moreover, $h_e$ denotes the diameter of element $\Omega_e$, and $C_\textup{pen}$ is a penalty constant associated with enforcement of the continuity of the velocity and magnetic fields across element boundary boundaries. As is typical with interior penalty methods, the penalty constant $C_\textup{pen}$ must be chosen sufficiently large to ensure the semi-discrete formulation is energy stable \cite{burman2007continuous}. In our later computations, we have selected $C_\textup{pen} = (k+1)(k+2)$ as we found this choice to be sufficient for energy stability. The above semi-discrete formulation holds for both $d = 2$ and $d = 3$, though only two-dimensional numerical examples are presented in this paper.  We also did not include the differential quantities $\partial \Omega$ and $\partial \Gamma$ in the above integrands for the sake of conciseness and will continue to do so throughout this paper.

We now show a consistency result for our semi-discrete formulation.\\

\begin{prop}[Consistency]
The semi-discrete HDG method presented in Equations \eqref{p}-\eqref{b_ic} is consistent provided the exact solution $(\mathbf{u},p,\mathbf{B},r)$ of the incompressible MHD equations is sufficiently smooth. That is, Equations \eqref{p}-\eqref{b_ic} hold if we replace $(\mathbf{u}^h,\hat{\mathbf{u}}^h,p^h,\hat{p}^h,\mathbf{B}^h,\hat{\mathbf{B}}^h,r^h,\hat{r}^h)$ with $(\mathbf{u},\mathbf{u}|_{\tilde{\Gamma}},p,p|_{\tilde{\Gamma}},\mathbf{B},\mathbf{B}|_{\tilde{\Gamma}},r,r|_{\tilde{\Gamma}})$.\\
\end{prop}

\textit{Proof:} Note that, for a sufficiently smooth exact solution $(\mathbf{u},p,\mathbf{B},r)$, Equations \eqref{eq:strong_mom}-\eqref{eq:strong_BC_2} hold in a pointwise manner. Using this information, we show that \eqref{p}-\eqref{b_ic} hold if we replace $(\mathbf{u}^h,\hat{\mathbf{u}}^h,p^h,\hat{p}^h,\mathbf{B}^h,\hat{\mathbf{B}}^h,r^h,\hat{r}^h)$ with $(\mathbf{u},\mathbf{u}|_{\tilde{\Gamma}},p,p|_{\tilde{\Gamma}},\mathbf{B},\mathbf{B}|_{\tilde{\Gamma}},r,r|_{\tilde{\Gamma}})$. Equations \eqref{m} and \eqref{mono} trivially hold since $\nabla \cdot \mathbf{u} \equiv \nabla \cdot \mathbf{B} \equiv 0 $ over $\Omega \times (0,\infty)$, Equations \eqref{mc} and \eqref{monoc} trivially hold since $\mathbf{u} \equiv \mathbf{u}|_{\tilde{\Gamma}}$ and $\mathbf{B} \equiv \mathbf{B}|_{\tilde{\Gamma}}$ over $\tilde{\Gamma} \times (0,\infty)$, and Equations \eqref{v_ic} and \eqref{b_ic} trivially hold since $\mathbf{u}(:,0) \equiv \mathbf{u}_0$ and $\mathbf{B}(:,0) \equiv \mathbf{B}_0$ over $\Omega$.
 Equation \eqref{p} holds since
\begin{equation}
    \begin{aligned}
    \Sigma_e&\int_{\Omega_e}\frac{\partial\mathbf{u}}{\partial t}\cdot\mathbf{v}^h-\Sigma_e\int_{\Omega_e}(\mathbf{u}\otimes\mathbf{u}):\nabla\mathbf{v}^h
    -\Sigma_e\int_{\Omega_e}\frac{1}{\rho}p\mathds{1}:\nabla\mathbf{v}^h+\Sigma_e\int_{\Omega_e}2\nu\nabla^s\mathbf{u}:\nabla^s\mathbf{v}^h\\
    +\Sigma_e&\int_{\Omega_e}\frac{1}{\rho \mu_0}\left(\mathbf{B}\otimes\mathbf{B}-\frac{1}{2}|\mathbf{B}|^2\mathds{1}\right):\nabla\mathbf{v}^h
    +\Sigma_e\Sigma_i\int_{\Gamma_{e_i}} (1-\lambda) (\mathbf{u}\otimes\mathbf{u}):(\mathbf{n}\otimes\mathbf{v}^h)\\
    +\Sigma_e\Sigma_i&\int_{\Gamma_{e_i}}\lambda (\mathbf{u}\otimes\mathbf{u}|_{\tilde{\Gamma}}):(\mathbf{n}\otimes\mathbf{v}^h)
    +\Sigma_e\Sigma_i\int_{\Gamma_{e_i}}\frac{1}{\rho}p|_{\tilde{\Gamma}}\mathds{1}:(\mathbf{n}\otimes\mathbf{v}^h)\\
    -\Sigma_e\Sigma_i&\int_{\Gamma_{e_i}}2\nu\nabla^s\mathbf{u}:(\mathbf{n}\otimes\mathbf{v}^h)
    +\Sigma_e\Sigma_i\int_{\Gamma_{e_i}}2\frac{C_\textup{pen}}{h_e}\nu\left( \mathbf{n} \otimes (\mathbf{u}-\mathbf{u}|_{\tilde{\Gamma}}) \right) :(\mathbf{n} \otimes \mathbf{v}^h)\\
    -\Sigma_e\Sigma_i&\int_{\Gamma_{e_i}}2\nu\left(\mathbf{n} \otimes (\mathbf{u}-\mathbf{u}|_{\tilde{\Gamma}})\right):\nabla^s\mathbf{v}^h
    -\Sigma_e\Sigma_i\int_{\Gamma_{e_i}}\frac{1}{\rho \mu_0}\left(\mathbf{B}|_{\tilde{\Gamma}}\otimes\mathbf{B}|_{\tilde{\Gamma}}-\frac{1}{2}|\mathbf{B}|_{\tilde{\Gamma}}|^2\mathds{1}\right):(\mathbf{n} \otimes \mathbf{v}^h)\\
    -\Sigma_e&\int_{\Omega_e}\mathbf{f}_v\cdot\mathbf{v}^h \\
    =\Sigma_e&\int_{\Omega_e}\left[\frac{\partial\mathbf{u}}{\partial t}+\nabla\cdot\left(\mathbf{u}\otimes\mathbf{u}\right)+\frac{1}{\rho}\nabla p-\nabla\cdot\left(2\nu\nabla^s\mathbf{u}\right)-\nabla\cdot\left(\frac{1}{\rho \mu_0}\left(\mathbf{B}\otimes\mathbf{B}-\frac{1}{2}|\mathbf{B}|^2\mathds{1}\right)\right)-\mathbf{f}_v\right]\cdot\mathbf{v}^h
    \end{aligned}
\end{equation}
for all $\mathbf{v}^h \in V^h$ by reverse integration by parts and
\begin{equation}
\begin{split}
    \frac{\partial\mathbf{u}}{\partial t}+\nabla\cdot\left(\mathbf{u}\otimes\mathbf{u}\right)+\frac{1}{\rho}\nabla p-\nabla\cdot\left(2\nu\nabla^s\mathbf{u}\right)-\nabla\cdot\left(\frac{1}{\rho \mu_0}\left(\mathbf{B}\otimes\mathbf{B}-\frac{1}{2}|\mathbf{B}|^2\mathds{1}\right)\right)&-\mathbf{f}_v=0
\end{split}
\end{equation}
in $\Omega\times(0,\infty)$, and Equation \eqref{B} holds since \begin{equation}
\begin{aligned}
\Sigma_e&\int_{\Omega_e}\frac{\partial\mathbf{B}}{\partial t}\cdot\mathbf{w}^h - \Sigma_e\int_{\Omega_e} (\mathbf{u}\otimes\mathbf{B}):\nabla\mathbf{w}^h
-\Sigma_e\int_{\Omega_e} r\mathds{1}:\nabla\mathbf{w}^h
+\Sigma_e\int_{\Omega_e} 2\frac{\eta}{\mu_0}\nabla^a\mathbf{B}:\nabla^a\mathbf{w}^h\\
-\Sigma_e&\int_{\Omega_e}(\mathbf{B}\otimes\mathbf{w}^h):\nabla\mathbf{u}
+\Sigma_e\Sigma_i\int_{\Gamma_{e_i}} (1-\lambda) (\mathbf{u}\otimes\mathbf{B}):(\mathbf{n}\otimes\mathbf{w}^h)\\
+\Sigma_e\Sigma_i&\int_{\Gamma_{e_i}}\lambda (\mathbf{u}\otimes \mathbf{B}|_{\tilde{\Gamma}}):(\mathbf{n}\otimes\mathbf{w}^h)
+\Sigma_e\Sigma_i\int_{\Gamma_{e_i}}r|_{\tilde{\Gamma}}\mathds{1}:(\mathbf{n}\otimes\mathbf{w}^h)\\
-\Sigma_e\Sigma_i&\int_{\Gamma_{e_i}} 2\frac{\eta}{\mu_0}\nabla^a\mathbf{B}:(\mathbf{n}\otimes\mathbf{w}^h)
+\Sigma_e\Sigma_i\int_{\Gamma_{e_i}} 2\frac{\eta}{\mu_0}\frac{C_\textup{pen}}{h_e}\left(\mathbf{n}\otimes(\mathbf{B}-\mathbf{B}|_{\tilde{\Gamma}})\right):(\mathbf{n}\otimes\mathbf{w}^h)\\
-\Sigma_e\Sigma_i&\int_{\Gamma_{e_i}} 2\frac{\eta}{\mu_0}\left(\mathbf{n}\otimes(\mathbf{B}-\mathbf{B}|_{\tilde{\Gamma}})\right):\nabla^a\mathbf{w}^h
-\Sigma_e\int_{\Omega_e}\mathbf{f}_m\cdot\mathbf{w}^h \\
=\Sigma_e&\int_{\Omega_e}\left[\frac{\partial \mathbf{B}}{\partial t}+\nabla\cdot(\mathbf{u}\otimes\mathbf{B})-\nabla\cdot(\mathbf{B}\otimes\mathbf{u})-\nabla\cdot\left(\frac{\eta}{\mu_0}\nabla^a\mathbf{B}\right)+\nabla r - \mathbf{f}_m\right]\cdot\mathbf{w}^h
\end{aligned}
\end{equation}
for all $\mathbf{w}^h \in W^h$  by reverse integration by parts and
\begin{equation}
    \frac{\partial \mathbf{B}}{\partial t}+\nabla\cdot(\mathbf{u}\otimes\mathbf{B})-\nabla\cdot(\mathbf{B}\otimes\mathbf{u})-\nabla\cdot\left(\frac{\eta}{\mu_0}\nabla^a\mathbf{B}\right)+\nabla r - \mathbf{f}_m=0
\end{equation}
in $\Omega\times(0,\infty)$. Finally, Equation \eqref{pc} holds since
\begin{equation}
    \begin{aligned}
        \Sigma_e\Sigma_i&\int_{\Gamma_{e_i}} (1-\lambda) (\mathbf{u}\otimes\mathbf{u}):(\mathbf{n}\otimes\mathbf{\hat{v}}^h) +\Sigma_e\Sigma_i\int_{\Gamma_{e_i}}\lambda (\mathbf{u}\otimes\mathbf{u}|_{\tilde{\Gamma}}):(\mathbf{n}\otimes\mathbf{\hat{v}}^h) \\
+\Sigma_e\Sigma_i&\int_{\Gamma_{e_i}}\frac{1}{\rho}p|_{\tilde{\Gamma}}\mathds{1}:(\mathbf{n}\otimes\mathbf{\hat{v}}^h)
-\Sigma_e\Sigma_i\int_{\Gamma_{e_i}}2\nu\nabla^s\mathbf{u}:(\mathbf{n}\otimes\mathbf{\hat{v}}^h) \\
+\Sigma_e\Sigma_i&\int_{\Gamma_{e_i}}2\frac{C_\textup{pen}}{h_e}\nu\left( \mathbf{n} \otimes (\mathbf{u}-\mathbf{u}|_{\tilde{\Gamma}}) \right) :(\mathbf{n} \otimes \mathbf{\hat{v}}^h) -\Sigma_e\Sigma_i\int_{\Gamma_{e_i}}\frac{1}{\rho \mu_0}\left(\mathbf{B}|_{\tilde{\Gamma}}\otimes\mathbf{B}|_{\tilde{\Gamma}}-\frac{1}{2}|\mathbf{B}|_{\tilde{\Gamma}}|^2\mathds{1}\right):(\mathbf{n} \otimes \mathbf{\hat{v}}^h)\\
-\Sigma_e\Sigma_i&\int_{\Gamma_N\cap\Gamma_{e_i}}(1-\lambda)(\mathbf{u}\otimes\mathbf{u}|_{\tilde{\Gamma}}):(\mathbf{n}\otimes\mathbf{\hat{v}}^h)
+\Sigma_e\Sigma_i\int_{\Gamma_N\cap\Gamma_{e_i}}\mathbf{h}_v\cdot\mathbf{\hat{v}}^h\\
        =-\Sigma_e\Sigma_i&\int_{\Gamma_N\cap\Gamma_{e_i}}\left(\mathbf{n}\cdot\left(-\frac{1}{\rho}p\mathds{1}+2\nu\nabla^s\mathbf{u}+\frac{1}{\rho \mu_0}\left(\mathbf{B}\otimes\mathbf{B}-\frac{1}{2}|\mathbf{B}|^2\mathds{1}\right)\right)-\text{min}(\mathbf{u}\cdot\mathbf{n},0)\mathbf{u}-\mathbf{h}_v\right)\cdot\mathbf{\hat{v}}^h
    \end{aligned}
\end{equation}
for all $\mathbf{\hat{v}}^h \in \hat{V}_0^h$ and
\begin{equation}
    \mathbf{n}\cdot\left(-\frac{1}{\rho}p\mathds{1}+2\nu\nabla^s\mathbf{u}+\frac{1}{\rho \mu_0}\left(\mathbf{B}\otimes\mathbf{B}-\frac{1}{2}|\mathbf{B}|^2\mathds{1}\right)\right)-\text{min}(\mathbf{u}\cdot\mathbf{n},0)\mathbf{u}-\textbf{h}_v = \mathbf{0}
\end{equation}
on $\Gamma_N\times(0,\infty)$, and Equation \eqref{Bc} holds since
\begin{equation}
\begin{aligned}
    \Sigma_e\Sigma_i&\int_{\Gamma_{e_i}} (1-\lambda) (\mathbf{u}\otimes\mathbf{B}):(\mathbf{n}\otimes\mathbf{\hat{w}}^h)
+\Sigma_e\Sigma_i\int_{\Gamma_{e_i}}\lambda (\mathbf{u}\otimes \mathbf{B}|_{\tilde{\Gamma}}):(\mathbf{n}\otimes\mathbf{\hat{w}}^h)\\
+\Sigma_e\Sigma_i&\int_{\Gamma_{e_i}}r|_{\tilde{\Gamma}}\mathds{1}:(\mathbf{n}\otimes\mathbf{\hat{w}}^h)
-\Sigma_e\Sigma_i\int_{\Gamma_{e_i}} 2\frac{\eta}{\mu_0}\nabla^a\mathbf{B}:(\mathbf{n}\otimes\mathbf{\hat{w}}^h)\\
+\Sigma_e\Sigma_i&\int_{\Gamma_{e_i}} 2\frac{\eta}{\mu_0}\frac{C_\textup{pen}}{h_e}\left(\mathbf{n}\otimes(\mathbf{B}-\mathbf{B}|_{\tilde{\Gamma}})\right):(\mathbf{n}\otimes\mathbf{\hat{w}}^h)\\
-\Sigma_e\Sigma_i&\int_{\Gamma_{e_i}} \left(\mathbf{B}|_{\tilde{\Gamma}}\otimes\mathbf{\hat{w}}^h - \frac{1}{2} \left( \mathbf{B}|_{\tilde{\Gamma}}\cdot\mathbf{\hat{w}}^h \right) \mathds{1} \right):(\mathbf{n}\otimes\left(\mathbf{u}-\mathbf{u}|_{\tilde{\Gamma}}\right))
\\
-\Sigma_e\Sigma_i&\int_{\Gamma_N \cap \Gamma_{e_i}} (1-\lambda) (\mathbf{u}\otimes\mathbf{B}|_{\tilde{\Gamma}}):(\mathbf{n}\otimes\mathbf{\hat{w}}^h)+\Sigma_e\Sigma_i\int_{\Gamma_N\cap\Gamma_{e_i}}\mathbf{h}_m\cdot\mathbf{\hat{w}}^h\\
=-\Sigma_e\Sigma_i&\int_{\Gamma_N\cap\Gamma_{e_i}}\left(\mathbf{n}\cdot\left(-r\mathds{1}+2\frac{\eta}{\mu_0}\nabla^a\mathbf{B}\right)-\min (\mathbf{u}\cdot\mathbf{n},0)\mathbf{B}-\mathbf{h}_m\right)\cdot\mathbf{\hat{w}}^h
\end{aligned}
\end{equation}
for all $\mathbf{\hat{w}}^h \in \hat{W}_0^h$ and
\begin{equation}
    \mathbf{n}\cdot\left(-r\mathds{1}+2\frac{\eta}{\mu_0}\nabla^a\mathbf{B}\right)-\min (\mathbf{u}\cdot\mathbf{n},0)\mathbf{B}-\mathbf{h}_m = \mathbf{0}
\end{equation}
on $\Gamma_N\times(0,\infty)$. $\square$

\section{Conservation Properties of the Semi-Discrete HDG Method}

In this section, we prove three conservation results for our semi-discrete HDG method. Our first result shows that our semi-discrete HDG method conserves mass in a pointwise manner.\\

\begin{prop}[Pointwise Mass Conservation]
If $\mathbf{u}^h\in V^h$ and $\mathbf{\hat{u}}^h\in\hat{V}^h$ satisfy Equations \eqref{m} and \eqref{mc}, then
\begin{align}
    \nabla\cdot\mathbf{u}^h&=0 & \forall\mathbf{x}\in\Omega_e & \textup{ and } \forall\Omega_e\in\mathcal{T}\\
    \llbracket\mathbf{u}^h\rrbracket&=0 & \forall\mathbf{x}\in F & \textup{ and } \forall F\in\mathcal{F}_{\mathrm{int}}\\
    \mathbf{u}^h\cdot\mathbf{n} &= \mathbf{\hat{u}}^h\cdot\mathbf{n} & \forall\mathbf{x}\in F & \textup{ and } \forall F\in\mathcal{F}_{\mathrm{bdy}}.
\end{align}
\label{prop:divu}
\end{prop}

\textit{Proof:} From Equation \ref{m}, it follows that
\begin{equation}
    \int_{\Omega_e} q^h_v\nabla\cdot\mathbf{u}^h=0\quad\forall q^h_v\in P_{k-1}(\Omega_e) \quad\forall\Omega_e\in\mathcal{T}.
    \label{divu}
\end{equation}
Since $\nabla\cdot\mathbf{u}^h|_{\Omega_e} \in P_{k-1}(\Omega_e)$ for every $\Omega_e \in \mathcal{T}$, we can take $q^h_v=\nabla\cdot\mathbf{u}^h|_{\Omega_e}$ in \eqref{divu}, yielding $\int_{\Omega_e}(\nabla\cdot\mathbf{u}^h)^2=0$ for every $\Omega\in\mathcal{T}$. Thus $\nabla\cdot\mathbf{u}^h\equiv 0$ in $\Omega_e$ for every $\Omega_e\in\mathcal{T}$.

From Equation \ref{mc}, it follows that
\begin{equation}
    \int_F\llbracket\mathbf{u}^h\rrbracket q^h_v=0 \quad\forall q^h_v\in P_{k}(F) \quad\forall F\in\mathit{F}_{\mathrm{int}}.
    \label{jumpu}
\end{equation}
Since $\llbracket\mathbf{u}^h\rrbracket|_F\in P_k(F)$ for all $F\in\mathcal{F}_{\mathrm{int}}$, we can take $q^h_v=\llbracket\mathbf{u}^h\rrbracket|_F$ in \eqref{jumpu}, yielding $\int_F\llbracket\mathbf{u}^h\rrbracket^2=0$ for all $F\in\mathcal{F}_{\mathrm{int}}$. Thus $\llbracket\mathbf{u}^h\rrbracket\equiv0$ on $\tilde{\Gamma}/\partial\Omega$.

From Equation \eqref{mc}, it also follows that
\begin{equation}
    \int_F\left(\mathbf{u}^h-\mathbf{\hat{u}}^h\right)\cdot\mathbf{n}q^h_v=0 \quad\forall q\in P_k(F) \quad\forall F\in\mathit{F}_{\mathrm{bdy}}.
    \label{uminusuhat}
\end{equation}
Since $\left(\mathbf{u}^h-\mathbf{\hat{u}}^h\right)\cdot\mathbf{n}|_F\in P_k(F)$ for all $F\in\mathcal{F}_{bdy}$, we can take $q^h_v=\left(\mathbf{u}^h-\mathbf{\hat{u}}^h\right)\cdot\mathbf{n}|_F$ in \eqref{uminusuhat}, yielding $\int_F\left(\left(\mathbf{u}^h-\mathbf{\hat{u}}^h\right)\cdot\mathbf{n}\right)^2=0$ for all $F\in\mathcal{F}_{\mathrm{bdy}}$. Thus $\left(\mathbf{u}^h-\mathbf{\hat{u}}^h\right)\cdot\mathbf{n}\equiv0$ on $\partial\Omega\cap\tilde{\Gamma}$. $\square$\\

Our next result shows that our semi-discrete HDG method is pointwise absent of magnetic monopoles.\\

\begin{prop}[Pointwise Absence of Magnetic Monopoles]
If $\mathbf{B}^h\in W^h$ and $\mathbf{\hat{B}}^h\in\hat{W}^h$ satisfy Equations \eqref{mono} and \eqref{monoc}, then
\begin{align}
    \nabla\cdot\mathbf{B}^h&=0 & \forall\mathbf{x}\in\Omega_e & \textup{ and } \forall\Omega_e\in\mathcal{T}\\
    \llbracket\mathbf{B}^h\rrbracket&=0 & \forall\mathbf{x}\in F & \textup{ and } \forall F\in\mathcal{F}_{\mathrm{int}}\\
    \mathbf{B}^h\cdot\mathbf{n} &= \mathbf{\hat{B}}^h\cdot\mathbf{n} & \forall\mathbf{x}\in F & \textup{ and } \forall F\in\mathcal{F}_{\mathrm{bdy}}.
\end{align}
\label{prop:divB}
\end{prop}

\textit{Proof:} This result can be proven using the same argument as the proof of Proposition \ref{prop:divu}. $\square$\\

Our final result gives a global momentum balance law for our semi-discrete HDG method.\\

\begin{prop}[Global Momentum Balance]
The solution $(\mathbf{u}^h,\hat{\mathbf{u}}^h,p^h,\hat{p}^h,\mathbf{B}^h,\hat{\mathbf{B}}^h,r^h,\hat{r}^h)$ of the semi-discrete HDG method presented in Equations \eqref{p}-\eqref{b_ic} satisfies the following if $\Gamma_D=\emptyset$:

\begin{equation}
    \frac{d}{dt} \sum_e \int_{\Omega_e} \rho \mathbf{u}^h = \sum_e \int_{\Omega_e}\rho\mathbf{f}_v - \sum_e \sum_i \int_{\Gamma_N\cap\Gamma_{e_i}}(1-\lambda)\rho(\mathbf{\hat{u}}^h\cdot\mathbf{n})\mathbf{\hat{u}}^h - \sum_e \sum_i \int_{\Gamma_N\cap\Gamma_{e_i}} \rho \mathbf{h}_v.
\end{equation}
\end{prop}

\textit{Proof:} This result follows by setting $\mathbf{v}^h= \rho\mathbf{e}_j$ in Equation \eqref{p}, $\mathbf{\hat{v}}^h=-\rho\mathbf{e}_j$ in Equation \eqref{pc}, and summing the two expressions together. $\square$

\section{Energy Stability of the Semi-Discrete HDG Formulation}

The following is a global energy stability result for our semi-discrete HDG method.\\

\begin{prop}[Global Energy Stability]
The solution $(\mathbf{u}^h,\hat{\mathbf{u}}^h,p^h,\hat{p}^h,\mathbf{B}^h,\hat{\mathbf{B}}^h,r^h,\hat{r}^h)$ of the semi-discrete HDG method presented in Equations \eqref{p}-\eqref{b_ic} satisfies the following if $\mathbf{f}_v \equiv \mathbf{0}$, $\mathbf{f}_m \equiv \mathbf{0}$, $\mathbf{g}_v \equiv \mathbf{0}$, $\mathbf{g}_m \equiv \mathbf{0}$, $\mathbf{h}_v \equiv \mathbf{0}$, $\mathbf{h}_m \equiv \mathbf{0}$, and $C_\textup{pen}$ is sufficiently large:
\begin{equation}
    \frac{d}{dt}\Sigma_e\int_{\Omega_e} \frac{1}{2}\left( \rho|\mathbf{u}^h|^2+\frac{1}{\mu_0}|\mathbf{B}^h|^2\right)\leq0.
\end{equation}
\label{prop:energy}
\end{prop}

\textit{Proof:} By setting $\mathbf{v}^h=\rho \mathbf{u}^h$, $\mathbf{\hat{v}}^h=-\rho \mathbf{\hat{u}}^h$, $\mathbf{w}^h=\frac{1}{\mu_0}\mathbf{B}^h$, and $\mathbf{\hat{w}}^h=-\frac{1}{\mu_0}\mathbf{\hat{B}}^h$ in Equations \eqref{p}, \eqref{pc}, \eqref{B}, and \eqref{Bc}, summing the resulting equations, and exploiting Propositions \ref{prop:divu} and \ref{prop:divB}, we find that
\begin{equation}
\Sigma_e \int_{\Omega_e}\rho \frac{\partial\mathbf{u}^h}{\partial t} \cdot \mathbf{u}^h + \Sigma_e \int_{\Omega_e}\frac{1}{\mu_0}\frac{\partial\mathbf{B}^h}{\partial t} \cdot \mathbf{B}^h + I + II + III + IV = 0
\end{equation}
where
\begin{equation}
\begin{aligned}
I =& -\Sigma_e\int_{\Omega_e}\rho (\mathbf{u}^h\otimes\mathbf{u}^h):\nabla\mathbf{u}^h + \Sigma_e \Sigma_i \int_{\Gamma_{e_i}} \rho \left( \mathbf{u}^h \cdot \mathbf{n} \right) \left( (1-\lambda) \mathbf{u}^h + \lambda \mathbf{\hat{u}}^h \right) \cdot \left( \mathbf{u}^h - \mathbf{\hat{u}}^h \right) \\ 
& + \Sigma_e \Sigma_i \int_{\Gamma_N \cap \Gamma_{e_i}} (1-\lambda) \rho \left( \mathbf{u}^h \cdot \mathbf{n} \right) | \mathbf{\hat{u}}^h |^2
\end{aligned}
\end{equation}
\begin{equation}
\begin{aligned}
II =& -\Sigma_e\int_{\Omega_e}\frac{1}{\mu_0}(\mathbf{u}^h\otimes\mathbf{B}^h):\nabla\mathbf{u}^h + \Sigma_e \Sigma_i \int_{\Gamma_{e_i}} \frac{1}{\mu_0} \left( \mathbf{u}^h \cdot \mathbf{n} \right) \left( (1-\lambda) \mathbf{B}^h + \lambda \mathbf{\hat{B}}^h \right) \cdot \left( \mathbf{B}^h - \mathbf{\hat{B}}^h \right) \\ 
& + \Sigma_e \Sigma_i \int_{\Gamma_N \cap \Gamma_{e_i}} \frac{1}{\mu_0} (1-\lambda) \left( \mathbf{u}^h \cdot \mathbf{n} \right) | \mathbf{\hat{B}}^h |^2
\end{aligned}
\end{equation}
\begin{equation}
III = \Sigma_e\int_{\Omega_e}2\rho\nu|\nabla^s\mathbf{u}^h|^2 + \Sigma_e\Sigma_i \int_{\Gamma_{e_i}}2\frac{C_\textup{pen}}{h_e}\rho\nu|\mathbf{u}^h-\mathbf{\hat{u}}^h|^2-\Sigma_e\Sigma_i\int_{\Gamma_{e_i}}4\rho\nu\left(\nabla^s\mathbf{u}^h\cdot\mathbf{n}\right)\cdot\left(\mathbf{u}^h-\mathbf{\hat{u}}^h\right)
\end{equation}
\begin{equation}
IV = \Sigma_e\int_{\Omega_e}2\frac{\eta}{\mu_0^2}|\nabla^a\mathbf{B}^h|^2 + \Sigma_e\Sigma_i \int_{\Gamma_{e_i}}2\frac{C_\textup{pen}}{h_e}\frac{\eta}{\mu_0^2}|\mathbf{B}^h-\mathbf{\hat{B}}^h|^2-\Sigma_e\Sigma_i\int_{\Gamma_{e_i}}4\frac{\eta}{\mu_0^2}\left(\nabla^a\mathbf{B}^h\cdot\mathbf{n}\right)\cdot\left(\mathbf{B}^h-\mathbf{\hat{B}}^h\right).
\end{equation}
By the product rule, we have
\begin{align}\Sigma_e \int_{\Omega_e}\rho\frac{\partial\mathbf{u}^h}{\partial t} \cdot \mathbf{u}^h + \Sigma_e \int_{\Omega_e} \frac{1}{\mu_0} \frac{\partial\mathbf{B}^h}{\partial t} \cdot \mathbf{B}^h
&= \Sigma_e\int_{\Omega_e}\rho\frac{1}{2}\frac{\partial|\mathbf{u}^h|^2}{\partial t} + \Sigma_e\int_{\Omega_e} \frac{1}{\mu_0} \frac{1}{2}\frac{\partial|\mathbf{B}^h|^2}{\partial t} \nonumber \\
&= \frac{d}{dt}\Sigma_e\int_{\Omega_e} \frac{1}{2}\left( \rho|\mathbf{u}^h|^2+\frac{1}{\mu_0} |\mathbf{B}^h|^2\right).
\end{align}
Thus, it suffices to show that each of $I$, $II$, $III$, and $IV$ are non-negative. By Proposition \ref{prop:divu}, the product rule, and the divergence theorem, we have
\begin{align}
-\Sigma_e\int_{\Omega_e}\rho(\mathbf{u}^h\otimes\mathbf{u}^h):\nabla\mathbf{u}^h &= -\Sigma_e\int_{\Omega_e} \nabla \cdot \left( \frac{1}{2} \rho |\mathbf{u}^h|^2 \mathbf{u}^h \right) \nonumber \\
&= -\Sigma_e \Sigma_i \int_{\Gamma_{e_i}} \frac{1}{2} \rho \left( \mathbf{u}^h \cdot \mathbf{n} \right) |\mathbf{u}^h|^2.
\label{eq:I_1}
\end{align}
Using the identity $\lambda = \frac{1}{2} \left(1 - \frac{|\mathbf{u}\cdot\mathbf{n}|}{\mathbf{u}\cdot\mathbf{n}} \right)$, the fact that $\textbf{n}^+ + \textbf{n}^- = \mathbf{0}$ for all $F \in \mathcal{F}_{\text{int}}$, and the fact that $\mathbf{\hat{u}}^h = \mathbf{0}$ on $\Gamma_D \cap \tilde{\Gamma}$, we find that
\begin{equation}
\begin{aligned}
\Sigma_e \Sigma_i \int_{\Gamma_{e_i}} \rho \left( \mathbf{u}^h \cdot \mathbf{n} \right) \left( (1-\lambda) \mathbf{u}^h + \lambda \mathbf{\hat{u}}^h \right) \cdot \left( \mathbf{u}^h - \mathbf{\hat{u}}^h \right) = & \Sigma_e \Sigma_i \int_{\Gamma_{e_i}} \frac{1}{2} \rho \left( \mathbf{u}^h \cdot \mathbf{n} \right) | \mathbf{u}^h |^2 \\
&- \Sigma_e \Sigma_i \int_{\Gamma_N \cap \Gamma_{e_i}} \frac{1}{2} \rho \left( \mathbf{u}^h \cdot \mathbf{n} \right) | \mathbf{\hat{u}}^h |^2 \\
& + \Sigma_e \Sigma_i \int_{\Gamma_{e_i}} \frac{1}{2} \rho \left| \mathbf{u}^h \cdot \mathbf{n} \right| | \mathbf{u}^h - \mathbf{\hat{u}}^h |^2
\end{aligned}
\label{eq:I_2}
\end{equation}
and
\begin{equation}
\begin{aligned}
\Sigma_e \Sigma_i \int_{\Gamma_N \cap \Gamma_{e_i}} (1-\lambda) \rho \left( \mathbf{u}^h \cdot \mathbf{n} \right) | \mathbf{\hat{u}}^h |^2 = & \Sigma_e \Sigma_i \int_{\Gamma_N \cap \Gamma_{e_i}} \frac{1}{2} \rho \left( \mathbf{u}^h \cdot \mathbf{n} \right) | \mathbf{\hat{u}}^h |^2 \\
&+ \Sigma_e \Sigma_i \int_{\Gamma_N \cap \Gamma_{e_i}} \frac{1}{2} \rho \left| \mathbf{u}^h \cdot \mathbf{n} \right| | \mathbf{\hat{u}}^h |^2.
\end{aligned}
\label{eq:I_3}
\end{equation}
Combining Equations \eqref{eq:I_1}-\eqref{eq:I_3} with the definition of $I$, we attain
\begin{equation}
I = \Sigma_e \Sigma_i \int_{\Gamma_{e_i}} \frac{1}{2} \rho \left| \mathbf{u}^h \cdot \mathbf{n} \right| | \mathbf{u}^h - \mathbf{\hat{u}}^h |^2 + \Sigma_e \Sigma_i \int_{\Gamma_N \cap \Gamma_{e_i}} \frac{1}{2} \rho \left| \mathbf{u}^h \cdot \mathbf{n} \right| | \mathbf{\hat{u}}^h |^2 \geq 0.
\end{equation}
Using a similar argument, we find that
\begin{equation}
II = \Sigma_e \Sigma_i \int_{\Gamma_{e_i}} \frac{1}{2\mu_0} \left| \mathbf{u}^h \cdot \mathbf{n} \right| | \mathbf{B}^h - \mathbf{\hat{B}}^h |^2 + \Sigma_e \Sigma_i \int_{\Gamma_N \cap \Gamma_{e_i}} \frac{1}{2\mu_0} \left| \mathbf{u}^h \cdot \mathbf{n} \right| | \mathbf{\hat{B}}^h |^2 \geq 0.
\end{equation}
Using the Cauchy-Schwarz inequality and Young's inequality, we can show
\begin{equation}
    -\Sigma_e\Sigma_i\int_{\Gamma_{e_i}} 4\rho\nu\left(\nabla^s\mathbf{u}^h\cdot\mathbf{n}\right)\cdot\left(\mathbf{u}^h-\mathbf{\hat{u}}^h\right)\leq \Sigma_e\Sigma_i\int_{\Gamma_{e_i}}2\frac{h_e}{C_\textup{pen}}\rho\nu|\nabla^s\mathbf{u}^h \cdot \mathbf{n}|^2 + \Sigma_e\Sigma_i\int_{\Gamma_{e_i}} 2\frac{C_\textup{pen}}{h_e} \rho\nu|\mathbf{u}^h-\mathbf{\hat{u}}^h|^2.
\end{equation}
If $C_\textup{pen}$ is large enough so the trace inequality
\begin{equation}
\Sigma_e\Sigma_i\int_{\Gamma_{e_i}}2\frac{h_e}{C_\textup{pen}}\rho\nu|\nabla^s\mathbf{u}^h \cdot \mathbf{n}|^2 \leq \Sigma_e \int_{\Omega_e}2\rho\nu|\nabla^s\mathbf{u}^h|^2
\end{equation}
holds (see, e.g., \cite{warburton2003constants,evans2013explicit}), it follows that $III \geq 0$. If $C_\textup{pen}$ is also large enough so the trace inequality
\begin{equation}
\Sigma_e\Sigma_i\int_{\Gamma_{e_i}}2\frac{h_e}{C_\textup{pen}}\frac{\eta}{\mu^2_0}|\nabla^a\mathbf{B}^h \cdot \mathbf{n}|^2 \leq \Sigma_e \int_{\Omega_e}2\frac{\eta}{\mu^2_0}|\nabla^a\mathbf{B}^h|^2
\end{equation}
holds, a similar argument can be used to show $IV \geq 0$. As each of $I$, $II$, $III$, and $IV$ are non-negative, the desired result holds. $\square$

\section{Fully-Discrete HDG Formulation for the MHD Equations}

To discretize in time, we employ the generalized-$\alpha$ method. The generalized-$\alpha$ method was originally developed for second-order systems of differential-algebraic equations
by Chung and Hulbert \cite{genalpha1} and then later extended to first-order systems of differential-algebraic equations by Jansen et al. \cite{jansen2000generalized}. We adopt the method of Jansen et al. here. Provided the parameters of the generalized-$\alpha$ method are chosen appropriately, it is second-order-in-time and unconditionally stable. To proceed, let $\mathrm{U}$, $\mathrm{P}$, $\mathrm{B}$, and $\mathrm{R}$ denote the degree-of-freedom vectors associated with the interior velocity, pressure, magnetic, and magnetic pressure fields, and let $\mathrm{\hat{U}}$, $\mathrm{\hat{P}}$, $\mathrm{\hat{B}}$, and $\mathrm{\hat{R}}$ denote the degree-of-freedom vectors associated with the trace velocity, pressure, magnetic, and magnetic pressure fields. We can then write Equations \eqref{p}-\eqref{monoc} as a system of differential-algebraic equations of the form
\begin{equation}
    \mathrm{R}_\text{total}(\mathrm{\dot{U}},\mathrm{U},\mathrm{P},\mathrm{\hat{U}},\mathrm{\hat{P}},\mathrm{\dot{B}},\mathrm{B},\mathrm{R},\mathrm{\hat{B}},\mathrm{\hat{R}}) =\mathrm{0}
    \label{resid}
\end{equation}
where $\dot{(\ast)}$ denotes differentiation in time. Given the values of $\mathrm{\dot{U}}$, $\mathrm{U}$, $\mathrm{P}$, $\mathrm{\hat{U}}$, $\mathrm{\hat{P}}$, $\mathrm{\dot{B}}$, $\mathrm{B}$, $\mathrm{R}$, $\mathrm{\hat{B}}$, and $\mathrm{\hat{R}}$ at the $n^\text{th}$ time step, we find the values of these quantities at the $(n+1)^\text{st}$ time step by solving
\begin{equation}
    \mathrm{R}_\text{total}(\mathrm{\dot{U}}_{n+\alpha_m},\mathrm{U}_{n+\alpha_f},\mathrm{P}_{n+\alpha_f},\mathrm{\hat{U}}_{n+\alpha_f},\mathrm{\hat{P}}_{n+\alpha_f},\mathrm{\dot{B}}_{n+\alpha_m},\mathrm{B}_{n+\alpha_f},\mathrm{R}_{n+\alpha_f},\mathrm{\hat{B}}_{n+\alpha_f},\mathrm{\hat{R}}_{n+\alpha_f}) =\mathrm{0}
    \label{residual}
\end{equation}
together with the Newmark equations
\begin{align}
\mathrm{U}_{n+1} &= \mathrm{U}_{n} + \Delta t_n \left( (1-\gamma) \mathrm{\dot{U}}_{n} + \gamma \mathrm{\dot{U}}_{n+1} \right) \\
\mathrm{B}_{n+1} &= \mathrm{B}_{n} + \Delta t_n \left( (1-\gamma) \mathrm{\dot{B}}_{n} + \gamma \mathrm{\dot{B}}_{n+1} \right)
\label{newmark}
\end{align}
where
\begin{align}
\mathrm{\dot{U}}_{n+\alpha_m} &= \mathrm{\dot{U}}_{n} + \alpha_m \left( \mathrm{\dot{U}}_{n+1} - \mathrm{\dot{U}}_{n} \right) & \mathrm{\dot{B}}_{n+\alpha_m} &= \mathrm{\dot{B}}_{n} + \alpha_m \left( \mathrm{\dot{B}}_{n+1} - \mathrm{\dot{B}}_{n} \right) \\
\mathrm{U}_{n+\alpha_f} &= \mathrm{U}_{n} + \alpha_f \left( \mathrm{U}_{n+1} - \mathrm{U}_{n} \right) & \mathrm{B}_{n+\alpha_f} &= \mathrm{B}_{n} + \alpha_f \left( \mathrm{B}_{n+1} - \mathrm{B}_{n} \right) \\
\mathrm{P}_{n+\alpha_f} &= \mathrm{P}_{n} + \alpha_f \left( \mathrm{P}_{n+1} - \mathrm{P}_{n} \right) & \mathrm{R}_{n+\alpha_f} &= \mathrm{R}_{n} + \alpha_f \left( \mathrm{R}_{n+1} - \mathrm{R}_{n} \right) \\
\mathrm{\hat{U}}_{n+\alpha_f} &= \mathrm{\hat{U}}_{n} + \alpha_f \left( \mathrm{\hat{U}}_{n+1} - \mathrm{\hat{U}}_{n} \right) & \mathrm{\hat{B}}_{n+\alpha_f} &= \mathrm{\hat{B}}_{n} + \alpha_f \left( \mathrm{\hat{B}}_{n+1} - \mathrm{\hat{B}}_{n} \right) \\
\mathrm{\hat{P}}_{n+\alpha_f} &= \mathrm{\hat{P}}_{n} + \alpha_f \left( \mathrm{\hat{P}}_{n+1} - \mathrm{\hat{P}}_{n} \right) & \mathrm{\hat{R}}_{n+\alpha_f} &= \mathrm{\hat{R}}_{n} + \alpha_f \left( \mathrm{\hat{R}}_{n+1} - \mathrm{\hat{R}}_{n} \right), \label{alpha_equations}
\end{align}
$\Delta t_n = t_{n+1} - t_n$ is the time step size between the $n^\text{th}$ and $(n+1)^\text{st}$ time steps, and $\gamma$, $\alpha_m$, and $\alpha_f$ are algorithmic parameters. Second-order accuracy is attained if $\gamma = \frac{1}{2} - \alpha_f + \alpha_m$ while unconditional stability requires $\alpha_m \geq \alpha_f \geq \frac{1}{2}$. We in particular set $\alpha_m$ and $\alpha_f$ using
\begin{align}
\alpha_m &= \frac{1}{2} \left( \frac{3-\rho_\infty}{1+\rho_\infty} \right) & \alpha_f &= \frac{1}{1+\rho_\infty}
\end{align}
where $\rho_{\infty} \in [0,1]$ is a free parameter,
and we set $\gamma = \frac{1}{2} - \alpha_f + \alpha_m$. This gives rise to a second-order-in-time and unconditionally stable method with tunable numerical dissipation. In particular, for a linear model problem, the generalized-$\alpha$ method annihilates the highest frequency in one time-step if $\rho_\infty = 0$ and preserves the highest frequency if $\rho_\infty = 1$. It should be noted that conventionally the pressure field is evaluated at the $(n+1)^\text{st}$ time step rather than the $(n+\alpha_f)^\text{th}$ time step. However, as recently shown in \cite{liu2021note}, this limits the accuracy of the pressure field to first-order-in-time. Moreover, we have found that evaluating the trace fields at the $(n+1)^\text{st}$ time step rather than the $(n+\alpha_f)^\text{th}$ time step limits the accuracy of all fields to first-order-in-time. It should be noted that application of the generalized-$\alpha$ method at the first step of a simulation requires knowledge of $\mathrm{\dot{U}}$, $\mathrm{P}$, $\mathrm{\hat{U}}$, $\mathrm{\hat{P}}$, $\mathrm{\dot{B}}$, $\mathrm{R}$, $\mathrm{\hat{B}}$, and $\mathrm{\hat{R}}$ at the initial time unless $\gamma = \alpha_m = \alpha_f = 1$, which corresponds to the backward Euler method. Generally speaking, such knowledge is not known, so it is recommended that one instead apply the backward Euler method for the first time step of a simulation before applying the generalized-$\alpha$ method for subsequent time steps. This does not upset the accuracy or stability properties of the generalized-$\alpha$ method.

To solve Equations \eqref{residual}-\eqref{alpha_equations}, one can Newton's method. This gives rise to the following predictor-multicorrector algorithm:\\

\textbf{Predictor Stage:}
Set:
\begin{align}
    \mathrm{\dot{U}}_{n+1}^{(0)} &= \frac{\gamma-1}{\gamma}\mathrm{\dot{U}}_n & \mathrm{U}_{n+1}^{(0)} &= \mathrm{U}_n & \mathbf{\mathrm{P}}_{n+1}^{(0)} &= \mathbf{\mathrm{P}}_n & \mathrm{\hat{U}}_{n+1}^{(0)} &=\mathrm{\hat{U}}_n & \mathrm{\hat{P}}_{n+1}^{(0)} &= \mathrm{\hat{P}}_n \\
    \mathrm{\dot{B}}_{n+1}^{(0)} &= \frac{\gamma-1}{\gamma}\mathrm{\dot{B}}_n & \mathrm{B}_{n+1}^{(0)} &= \mathrm{B}_n & \mathbf{\mathrm{R}}_{n+1}^{(0)} &= \mathbf{\mathrm{R}}_n & \mathrm{\hat{B}}_{n+1}^{(0)} &=\mathrm{\hat{B}}_n & \mathrm{\hat{R}}_{n+1}^{(0)} &= \mathrm{\hat{R}}_n
\end{align}

\textbf{Multicorrector Stage:}
Repeat the following steps for $i=1,2,...,i_{\text{max}}$.

\textit{Step 1:} Evaluate iterates at the $\alpha$-levels:
\begin{align}
\mathrm{\dot{U}}^{(i)}_{n+\alpha_m} &= \mathrm{\dot{U}}_{n} + \alpha_m \left( \mathrm{\dot{U}}^{(i-1)}_{n+1} - \mathrm{\dot{U}}_{n} \right) & \mathrm{\dot{B}}^{(i)}_{n+\alpha_m} &= \mathrm{\dot{B}}_{n} + \alpha_m \left( \mathrm{\dot{B}}^{(i-1)}_{n+1} - \mathrm{\dot{B}}_{n} \right) \\
\mathrm{U}^{(i)}_{n+\alpha_f} &= \mathrm{U}_{n} + \alpha_f \left( \mathrm{U}^{(i-1)}_{n+1} - \mathrm{U}_{n} \right) & \mathrm{B}^{(i)}_{n+\alpha_f} &= \mathrm{B}_{n} + \alpha_f \left( \mathrm{B}^{(i-1)}_{n+1} - \mathrm{B}_{n} \right) \\
\mathrm{P}^{(i)}_{n+\alpha_f} &= \mathrm{P}_{n} + \alpha_f \left( \mathrm{P}^{(i-1)}_{n+1} - \mathrm{P}_{n} \right) & \mathrm{R}^{(i)}_{n+\alpha_f} &= \mathrm{R}_{n} + \alpha_f \left( \mathrm{R}^{(i-1)}_{n+1} - \mathrm{R}_{n} \right) \\
\mathrm{\hat{U}}^{(i)}_{n+\alpha_f} &= \mathrm{\hat{U}}_{n} + \alpha_f \left( \mathrm{\hat{U}}^{(i-1)}_{n+1} - \mathrm{\hat{U}}_{n} \right) & \mathrm{\hat{B}}^{(i)}_{n+\alpha_f} &= \mathrm{\hat{B}}_{n} + \alpha_f \left( \mathrm{\hat{B}}^{(i-1)}_{n+1} - \mathrm{\hat{B}}_{n} \right) \\
\mathrm{\hat{P}}^{(i)}_{n+\alpha_f} &= \mathrm{\hat{P}}_{n} + \alpha_f \left( \mathrm{\hat{P}}^{(i-1)}_{n+1} - \mathrm{\hat{P}}_{n} \right) & \mathrm{\hat{R}}^{(i)}_{n+\alpha_f} &= \mathrm{\hat{R}}_{n} + \alpha_f \left( \mathrm{\hat{R}}^{(i-1)}_{n+1} - \mathrm{\hat{R}}_{n} \right) \label{alpha_equations2}
\end{align}

\textit{Step 2:} Use the solutions at the $\alpha$-levels to assemble the residual and tangent matrix of the linear system:
\begin{equation}
\left[
\begin{array}{cccccccc}
\mathrm{K}^{(i)}_{\mathbf{u}\mathbf{u}} & \mathrm{K}^{(i)}_{\mathbf{u}p} & \mathrm{K}^{(i)}_{\mathbf{u}\mathbf{\hat{u}}} & \mathrm{K}^{(i)}_{\mathbf{u}\hat{p}} & \mathrm{K}^{(i)}_{\mathbf{u}\mathbf{B}} & \mathrm{K}^{(i)}_{\mathbf{u}r} & \mathrm{K}^{(i)}_{\mathbf{u}\mathbf{\hat{B}}} & \mathrm{K}^{(i)}_{\mathbf{u}\hat{r}} \\

\mathrm{K}^{(i)}_{p\mathbf{u}} & \mathrm{K}^{(i)}_{pp} & \mathrm{K}^{(i)}_{p\mathbf{\hat{u}}} & \mathrm{K}^{(i)}_{p\hat{p}} & \mathrm{K}^{(i)}_{p\mathbf{B}} & \mathrm{K}^{(i)}_{pr} & \mathrm{K}^{(i)}_{p\mathbf{\hat{B}}} & \mathrm{K}^{(i)}_{p\hat{r}} \\

\mathrm{K}^{(i)}_{\mathbf{\hat{u}}\mathbf{u}} & \mathrm{K}^{(i)}_{\mathbf{\hat{u}}p} & \mathrm{K}^{(i)}_{\mathbf{\hat{u}}\mathbf{\hat{u}}} & \mathrm{K}^{(i)}_{\mathbf{\hat{u}}\hat{p}} & \mathrm{K}^{(i)}_{\mathbf{\hat{u}}\mathbf{B}} & \mathrm{K}^{(i)}_{\mathbf{\hat{u}}r} & \mathrm{K}^{(i)}_{\mathbf{\hat{u}}\mathbf{\hat{B}}} & \mathrm{K}^{(i)}_{\mathbf{\hat{u}}\hat{r}} \\

\mathrm{K}^{(i)}_{\hat{p}\mathbf{u}} & \mathrm{K}^{(i)}_{\hat{p}p} & \mathrm{K}^{(i)}_{\hat{p}\mathbf{\hat{u}}} & \mathrm{K}^{(i)}_{\hat{p}\hat{p}} & \mathrm{K}^{(i)}_{\hat{p}\mathbf{B}} & \mathrm{K}^{(i)}_{\hat{p}r} & \mathrm{K}^{(i)}_{\hat{p}\mathbf{\hat{B}}} & \mathrm{K}^{(i)}_{\hat{p}\hat{r}} \\

\mathrm{K}^{(i)}_{\mathbf{B}\mathbf{u}} & \mathrm{K}^{(i)}_{\mathbf{B}p} & \mathrm{K}^{(i)}_{\mathbf{B}\mathbf{\hat{u}}} & \mathrm{K}^{(i)}_{\mathbf{B}\hat{p}} & \mathrm{K}^{(i)}_{\mathbf{B}\mathbf{B}} & \mathrm{K}^{(i)}_{\mathbf{B}r} & \mathrm{K}^{(i)}_{\mathbf{B}\mathbf{\hat{B}}} & \mathrm{K}^{(i)}_{\mathbf{B}\hat{r}} \\

\mathrm{K}^{(i)}_{r\mathbf{u}} & \mathrm{K}^{(i)}_{rp} & \mathrm{K}^{(i)}_{r\mathbf{\hat{u}}} & \mathrm{K}^{(i)}_{r\hat{p}} & \mathrm{K}^{(i)}_{r\mathbf{B}} & \mathrm{K}^{(i)}_{rr} & \mathrm{K}^{(i)}_{r\mathbf{\hat{B}}} & \mathrm{K}^{(i)}_{r\hat{r}} \\

\mathrm{K}^{(i)}_{\mathbf{\hat{B}}\mathbf{u}} & \mathrm{K}^{(i)}_{\mathbf{\hat{B}}p} & \mathrm{K}^{(i)}_{\mathbf{\hat{B}}\mathbf{\hat{u}}} & \mathrm{K}^{(i)}_{\mathbf{\hat{B}}\hat{p}} & \mathrm{K}^{(i)}_{\mathbf{\hat{B}}\mathbf{B}} & \mathrm{K}^{(i)}_{\mathbf{\hat{B}}r} & \mathrm{K}^{(i)}_{\mathbf{\hat{B}}\mathbf{\hat{B}}} & \mathrm{K}^{(i)}_{\mathbf{\hat{B}}\hat{r}} \\

\mathrm{K}^{(i)}_{\hat{r}\mathbf{u}} & \mathrm{K}^{(i)}_{\hat{r}p} & \mathrm{K}^{(i)}_{\hat{r}\mathbf{\hat{u}}} & \mathrm{K}^{(i)}_{\hat{r}\hat{p}} & \mathrm{K}^{(i)}_{\hat{r}\mathbf{B}} & \mathrm{K}^{(i)}_{\hat{r}r} & \mathrm{K}^{(i)}_{\hat{r}\mathbf{\hat{B}}} & \mathrm{K}^{(i)}_{\hat{r}\hat{r}}

\end{array}
\right] \left[
\begin{array}{c}
\Delta \mathrm{\dot{U}}^{(i)}_{n+1} \\
\Delta \mathrm{P}^{(i)}_{n+1} \\
\Delta \mathrm{\hat{U}}^{(i)}_{n+1} \\
\Delta \mathrm{\hat{P}}^{(i)}_{n+1} \\
\Delta \mathrm{\dot{B}}^{(i)}_{n+1} \\
\Delta \mathrm{R}^{(i)}_{n+1} \\
\Delta \mathrm{\hat{B}}^{(i)}_{n+1} \\
\Delta \mathrm{\hat{R}}^{(i)}_{n+1}
\end{array}
\right] = -\left[
\begin{array}{c}
\mathrm{R}_\mathbf{u}^{(i)} \\
\mathrm{R}_p^{(i)} \\
\mathrm{R}_\mathbf{\hat{u}}^{(i)} \\
\mathrm{R}_{\hat{p}}^{(i)} \\
\mathrm{R}_\mathbf{B}^{(i)} \\
\mathrm{R}_r^{(i)} \\
\mathrm{R}_\mathbf{\hat{B}}^{(i)} \\
\mathrm{R}_{\hat{r}}^{(i)}
\end{array}
\right]
\label{eq:gen_alpha_system}
\end{equation}
where $\mathrm{R}_\mathbf{u}^{(i)}$, $\mathrm{R}_p^{(i)}$, $\mathrm{R}_\mathbf{\hat{u}}^{(i)}$, $\mathrm{R}_{\hat{p}}^{(i)}$, $\mathrm{R}_\mathbf{B}^{(i)}$, $\mathrm{R}_r^{(i)}$, $\mathrm{R}_\mathbf{\hat{B}}^{(i)}$, and $\mathrm{R}_{\hat{r}}^{(i)}$ are the components of
\begin{equation}
\mathrm{R}^{(i)}_\text{total} = \mathrm{R}_\text{total}\left(\mathrm{\dot{U}}^{(i)}_{n+\alpha_m},\mathrm{U}^{(i)}_{n+\alpha_f},\mathrm{P}^{(i)}_{n+\alpha_f},\mathrm{\hat{U}}^{(i)}_{n+\alpha_f},\mathrm{\hat{P}}^{(i)}_{n+\alpha_f},\mathrm{\dot{B}}^{(i)}_{n+\alpha_m},\mathrm{B}^{(i)}_{n+\alpha_f},\mathrm{R}^{(i)}_{n+\alpha_f},\mathrm{\hat{B}}^{(i)}_{n+\alpha_f},\mathrm{\hat{R}}^{(i)}_{n+\alpha_f}\right)
\end{equation}
associated with Equations \eqref{p}, \eqref{m}, \eqref{pc}, \eqref{mc}, \eqref{B}, \eqref{mono}, \eqref{Bc}, and \eqref{monoc} respectively and
\begin{align}
\mathrm{K}^{(i)}_{\mathbf{u}\mathbf{u}} &= \alpha_m \frac{\partial \mathrm{R}_\mathbf{u}^{(i)}}{\partial \mathrm{\dot{U}}_{n+\alpha_m}} + \alpha_f \gamma \Delta t_n \frac{\partial \mathrm{R}_\mathbf{u}^{(i)}}{\partial \mathrm{U}_{n+\alpha_f}} & \mathrm{K}^{(i)}_{p\mathbf{u}} &= \alpha_m \frac{\partial \mathrm{R}_p^{(i)}}{\partial \mathrm{\dot{U}}_{n+\alpha_m}} + \alpha_f \gamma \Delta t_n \frac{\partial \mathrm{R}_p^{(i)}}{\partial \mathrm{U}_{n+\alpha_f}} \\
\mathrm{K}^{(i)}_{\mathbf{\hat{u}}\mathbf{u}} &= \alpha_m \frac{\partial \mathrm{R}_\mathbf{\hat{u}}^{(i)}}{\partial \mathrm{\dot{U}}_{n+\alpha_m}} + \alpha_f \gamma \Delta t_n \frac{\partial \mathrm{R}_\mathbf{\hat{u}}^{(i)}}{\partial \mathrm{U}_{n+\alpha_f}} & \mathrm{K}^{(i)}_{\hat{p}\mathbf{u}} &= \alpha_m \frac{\partial \mathrm{R}_{\hat{p}}^{(i)}}{\partial \mathrm{\dot{U}}_{n+\alpha_m}} + \alpha_f \gamma \Delta t_n \frac{\partial \mathrm{R}_{\hat{p}}^{(i)}}{\partial \mathrm{U}_{n+\alpha_f}} \\
\mathrm{K}^{(i)}_{\mathbf{B}\mathbf{u}} &= \alpha_m \frac{\partial \mathrm{R}_\mathbf{B}^{(i)}}{\partial \mathrm{\dot{U}}_{n+\alpha_m}} + \alpha_f \gamma \Delta t_n \frac{\partial \mathrm{R}_\mathbf{B}^{(i)}}{\partial \mathrm{U}_{n+\alpha_f}} & \mathrm{K}^{(i)}_{r\mathbf{u}} &= \alpha_m \frac{\partial \mathrm{R}_r^{(i)}}{\partial \mathrm{\dot{U}}_{n+\alpha_m}} + \alpha_f \gamma \Delta t_n \frac{\partial \mathrm{R}_r^{(i)}}{\partial \mathrm{U}_{n+\alpha_f}} \\
\mathrm{K}^{(i)}_{\mathbf{\hat{B}}\mathbf{u}} &= \alpha_m \frac{\partial \mathrm{R}_\mathbf{\hat{B}}^{(i)}}{\partial \mathrm{\dot{U}}_{n+\alpha_m}} + \alpha_f \gamma \Delta t_n \frac{\partial \mathrm{R}_\mathbf{\hat{B}}^{(i)}}{\partial \mathrm{U}_{n+\alpha_f}} & \mathrm{K}^{(i)}_{\hat{r}\mathbf{u}} &= \alpha_m \frac{\partial \mathrm{R}_{\hat{r}}^{(i)}}{\partial \mathrm{\dot{U}}_{n+\alpha_m}} + \alpha_f \gamma \Delta t_n \frac{\partial \mathrm{R}_{\hat{r}}^{(i)}}{\partial \mathrm{U}_{n+\alpha_f}}
\end{align}
\begin{align}
\mathrm{K}^{(i)}_{\mathbf{u}p} &= \alpha_f \frac{\partial \mathrm{R}_\mathbf{u}^{(i)}}{\partial \mathrm{p}_{n+\alpha_f}} & \mathrm{K}^{(i)}_{pp} &= \alpha_f \frac{\partial \mathrm{R}_p^{(i)}}{\partial \mathrm{p}_{n+\alpha_f}} &
\mathrm{K}^{(i)}_{\mathbf{\hat{u}}p} &= \alpha_f \frac{\partial \mathrm{R}_\mathbf{\hat{u}}^{(i)}}{\partial \mathrm{p}_{n+\alpha_f}} & \mathrm{K}^{(i)}_{\hat{p}p} &= \alpha_f \frac{\partial \mathrm{R}_{\hat{p}}^{(i)}}{\partial \mathrm{p}_{n+\alpha_f}} \\
\mathrm{K}^{(i)}_{\mathbf{B}p} &= \alpha_f \frac{\partial \mathrm{R}_\mathbf{B}^{(i)}}{\partial \mathrm{p}_{n+\alpha_f}} & \mathrm{K}^{(i)}_{rp} &= \alpha_f \frac{\partial \mathrm{R}_r^{(i)}}{\partial \mathrm{p}_{n+\alpha_f}} &
\mathrm{K}^{(i)}_{\mathbf{\hat{B}}p} &= \alpha_f \frac{\partial \mathrm{R}_\mathbf{\hat{B}}^{(i)}}{\partial \mathrm{p}_{n+\alpha_f}} & \mathrm{K}^{(i)}_{\hat{r}p} &= \alpha_f \frac{\partial \mathrm{R}_{\hat{r}}^{(i)}}{\partial \mathrm{p}_{n+\alpha_f}}
\end{align}
\begin{align}
\mathrm{K}^{(i)}_{\mathbf{u}\mathbf{\hat{u}}} &= \alpha_f \frac{\partial \mathrm{R}_\mathbf{u}^{(i)}}{\partial \mathrm{\mathbf{\hat{u}}}_{n+\alpha_f}} & \mathrm{K}^{(i)}_{p\mathbf{\hat{u}}} &= \alpha_f \frac{\partial \mathrm{R}_p^{(i)}}{\partial \mathrm{\mathbf{\hat{u}}}_{n+\alpha_f}} &
\mathrm{K}^{(i)}_{\mathbf{\hat{u}}\mathbf{\hat{u}}} &= \alpha_f \frac{\partial \mathrm{R}_\mathbf{\hat{u}}^{(i)}}{\partial \mathrm{\mathbf{\hat{u}}}_{n+\alpha_f}} & \mathrm{K}^{(i)}_{\hat{p}\mathbf{\hat{u}}} &= \alpha_f \frac{\partial \mathrm{R}_{\hat{p}}^{(i)}}{\partial \mathrm{\mathbf{\hat{u}}}_{n+\alpha_f}} \\
\mathrm{K}^{(i)}_{\mathbf{B}\mathbf{\hat{u}}} &= \alpha_f \frac{\partial \mathrm{R}_\mathbf{B}^{(i)}}{\partial \mathrm{\mathbf{\hat{u}}}_{n+\alpha_f}} & \mathrm{K}^{(i)}_{r\mathbf{\hat{u}}} &= \alpha_f \frac{\partial \mathrm{R}_r^{(i)}}{\partial \mathrm{\mathbf{\hat{u}}}_{n+\alpha_f}} &
\mathrm{K}^{(i)}_{\mathbf{\hat{B}}\mathbf{\hat{u}}} &= \alpha_f \frac{\partial \mathrm{R}_\mathbf{\hat{B}}^{(i)}}{\partial \mathrm{\mathbf{\hat{u}}}_{n+\alpha_f}} & \mathrm{K}^{(i)}_{\hat{r}\mathbf{\hat{u}}} &= \alpha_f \frac{\partial \mathrm{R}_{\hat{r}}^{(i)}}{\partial \mathrm{\mathbf{\hat{u}}}_{n+\alpha_f}}
\end{align}
\begin{align}
\mathrm{K}^{(i)}_{\mathbf{u}\hat{p}} &= \alpha_f \frac{\partial \mathrm{R}_\mathbf{u}^{(i)}}{\partial \mathrm{\hat{p}}_{n+\alpha_f}} & \mathrm{K}^{(i)}_{p\hat{p}} &= \alpha_f \frac{\partial \mathrm{R}_p^{(i)}}{\partial \mathrm{\hat{p}}_{n+\alpha_f}} &
\mathrm{K}^{(i)}_{\mathbf{\hat{u}}\hat{p}} &= \alpha_f \frac{\partial \mathrm{R}_\mathbf{\hat{u}}^{(i)}}{\partial \mathrm{\hat{p}}_{n+\alpha_f}} & \mathrm{K}^{(i)}_{\hat{p}\hat{p}} &= \alpha_f \frac{\partial \mathrm{R}_{\hat{p}}^{(i)}}{\partial \mathrm{\hat{p}}_{n+\alpha_f}} \\
\mathrm{K}^{(i)}_{\mathbf{B}\hat{p}} &= \alpha_f \frac{\partial \mathrm{R}_\mathbf{B}^{(i)}}{\partial \mathrm{\hat{p}}_{n+\alpha_f}} & \mathrm{K}^{(i)}_{r\hat{p}} &= \alpha_f \frac{\partial \mathrm{R}_r^{(i)}}{\partial \mathrm{\hat{p}}_{n+\alpha_f}} &
\mathrm{K}^{(i)}_{\mathbf{\hat{B}}\hat{p}} &= \alpha_f \frac{\partial \mathrm{R}_\mathbf{\hat{B}}^{(i)}}{\partial \mathrm{\hat{p}}_{n+\alpha_f}} & \mathrm{K}^{(i)}_{\hat{r}\hat{p}} &= \alpha_f \frac{\partial \mathrm{R}_{\hat{r}}^{(i)}}{\partial \mathrm{\hat{p}}_{n+\alpha_f}}
\end{align}
\begin{align}
\mathrm{K}^{(i)}_{\mathbf{u}\mathbf{B}} &= \alpha_m \frac{\partial \mathrm{R}_\mathbf{u}^{(i)}}{\partial \mathrm{\dot{B}}_{n+\alpha_m}} + \alpha_f \gamma \Delta t_n \frac{\partial \mathrm{R}_\mathbf{u}^{(i)}}{\partial \mathrm{B}_{n+\alpha_f}} & \mathrm{K}^{(i)}_{p\mathbf{B}} &= \alpha_m \frac{\partial \mathrm{R}_p^{(i)}}{\partial \mathrm{\dot{B}}_{n+\alpha_m}} + \alpha_f \gamma \Delta t_n \frac{\partial \mathrm{R}_p^{(i)}}{\partial \mathrm{B}_{n+\alpha_f}} \\
\mathrm{K}^{(i)}_{\mathbf{\hat{u}}\mathbf{B}} &= \alpha_m \frac{\partial \mathrm{R}_\mathbf{\hat{u}}^{(i)}}{\partial \mathrm{\dot{B}}_{n+\alpha_m}} + \alpha_f \gamma \Delta t_n \frac{\partial \mathrm{R}_\mathbf{\hat{u}}^{(i)}}{\partial \mathrm{B}_{n+\alpha_f}} & \mathrm{K}^{(i)}_{\hat{p}\mathbf{B}} &= \alpha_m \frac{\partial \mathrm{R}_{\hat{p}}^{(i)}}{\partial \mathrm{\dot{B}}_{n+\alpha_m}} + \alpha_f \gamma \Delta t_n \frac{\partial \mathrm{R}_{\hat{p}}^{(i)}}{\partial \mathrm{B}_{n+\alpha_f}} \\
\mathrm{K}^{(i)}_{\mathbf{B}\mathbf{B}} &= \alpha_m \frac{\partial \mathrm{R}_\mathbf{B}^{(i)}}{\partial \mathrm{\dot{B}}_{n+\alpha_m}} + \alpha_f \gamma \Delta t_n \frac{\partial \mathrm{R}_\mathbf{B}^{(i)}}{\partial \mathrm{B}_{n+\alpha_f}} & \mathrm{K}^{(i)}_{r\mathbf{B}} &= \alpha_m \frac{\partial \mathrm{R}_r^{(i)}}{\partial \mathrm{\dot{B}}_{n+\alpha_m}} + \alpha_f \gamma \Delta t_n \frac{\partial \mathrm{R}_r^{(i)}}{\partial \mathrm{B}_{n+\alpha_f}} \\
\mathrm{K}^{(i)}_{\mathbf{\hat{B}}\mathbf{B}} &= \alpha_m \frac{\partial \mathrm{R}_\mathbf{\hat{B}}^{(i)}}{\partial \mathrm{\dot{B}}_{n+\alpha_m}} + \alpha_f \gamma \Delta t_n \frac{\partial \mathrm{R}_\mathbf{\hat{B}}^{(i)}}{\partial \mathrm{B}_{n+\alpha_f}} & \mathrm{K}^{(i)}_{\hat{r}\mathbf{B}} &= \alpha_m \frac{\partial \mathrm{R}_{\hat{r}}^{(i)}}{\partial \mathrm{\dot{B}}_{n+\alpha_m}} + \alpha_f \gamma \Delta t_n \frac{\partial \mathrm{R}_{\hat{r}}^{(i)}}{\partial \mathrm{B}_{n+\alpha_f}}
\end{align}
\begin{align}
\mathrm{K}^{(i)}_{\mathbf{u}r} &= \alpha_f \frac{\partial \mathrm{R}_\mathbf{u}^{(i)}}{\partial \mathrm{r}_{n+\alpha_f}} & \mathrm{K}^{(i)}_{pr} &= \alpha_f \frac{\partial \mathrm{R}_p^{(i)}}{\partial \mathrm{r}_{n+\alpha_f}} &
\mathrm{K}^{(i)}_{\mathbf{\hat{u}}r} &= \alpha_f \frac{\partial \mathrm{R}_\mathbf{\hat{u}}^{(i)}}{\partial \mathrm{r}_{n+\alpha_f}} & \mathrm{K}^{(i)}_{\hat{p}r} &= \alpha_f \frac{\partial \mathrm{R}_{\hat{p}}^{(i)}}{\partial \mathrm{r}_{n+\alpha_f}} \\
\mathrm{K}^{(i)}_{\mathbf{B}r} &= \alpha_f \frac{\partial \mathrm{R}_\mathbf{B}^{(i)}}{\partial \mathrm{r}_{n+\alpha_f}} & \mathrm{K}^{(i)}_{rr} &= \alpha_f \frac{\partial \mathrm{R}_r^{(i)}}{\partial \mathrm{r}_{n+\alpha_f}} &
\mathrm{K}^{(i)}_{\mathbf{\hat{B}}r} &= \alpha_f \frac{\partial \mathrm{R}_\mathbf{\hat{B}}^{(i)}}{\partial \mathrm{r}_{n+\alpha_f}} & \mathrm{K}^{(i)}_{\hat{r}r} &= \alpha_f \frac{\partial \mathrm{R}_{\hat{r}}^{(i)}}{\partial \mathrm{r}_{n+\alpha_f}}
\end{align}
\begin{align}
\mathrm{K}^{(i)}_{\mathbf{u}\mathbf{\hat{B}}} &= \alpha_f \frac{\partial \mathrm{R}_\mathbf{u}^{(i)}}{\partial \mathrm{\mathbf{\hat{B}}}_{n+\alpha_f}} & \mathrm{K}^{(i)}_{p\mathbf{\hat{B}}} &= \alpha_f \frac{\partial \mathrm{R}_p^{(i)}}{\partial \mathrm{\mathbf{\hat{B}}}_{n+\alpha_f}} &
\mathrm{K}^{(i)}_{\mathbf{\hat{u}}\mathbf{\hat{B}}} &= \alpha_f \frac{\partial \mathrm{R}_\mathbf{\hat{u}}^{(i)}}{\partial \mathrm{\mathbf{\hat{B}}}_{n+\alpha_f}} & \mathrm{K}^{(i)}_{\hat{p}\mathbf{\hat{B}}} &= \alpha_f \frac{\partial \mathrm{R}_{\hat{p}}^{(i)}}{\partial \mathrm{\mathbf{\hat{B}}}_{n+\alpha_f}} \\
\mathrm{K}^{(i)}_{\mathbf{B}\mathbf{\hat{B}}} &= \alpha_f \frac{\partial \mathrm{R}_\mathbf{B}^{(i)}}{\partial \mathrm{\mathbf{\hat{B}}}_{n+\alpha_f}} & \mathrm{K}^{(i)}_{r\mathbf{\hat{B}}} &= \alpha_f \frac{\partial \mathrm{R}_r^{(i)}}{\partial \mathrm{\mathbf{\hat{B}}}_{n+\alpha_f}} &
\mathrm{K}^{(i)}_{\mathbf{\hat{B}}\mathbf{\hat{B}}} &= \alpha_f \frac{\partial \mathrm{R}_\mathbf{\hat{B}}^{(i)}}{\partial \mathrm{\mathbf{\hat{B}}}_{n+\alpha_f}} & \mathrm{K}^{(i)}_{\hat{r}\mathbf{\hat{B}}} &= \alpha_f \frac{\partial \mathrm{R}_{\hat{r}}^{(i)}}{\partial \mathrm{\mathbf{\hat{B}}}_{n+\alpha_f}}
\end{align}
\begin{align}
\mathrm{K}^{(i)}_{\mathbf{u}\hat{r}} &= \alpha_f \frac{\partial \mathrm{R}_\mathbf{u}^{(i)}}{\partial \mathrm{\hat{r}}_{n+\alpha_f}} & \mathrm{K}^{(i)}_{p\hat{r}} &= \alpha_f \frac{\partial \mathrm{R}_p^{(i)}}{\partial \mathrm{\hat{r}}_{n+\alpha_f}} &
\mathrm{K}^{(i)}_{\mathbf{\hat{u}}\hat{r}} &= \alpha_f \frac{\partial \mathrm{R}_\mathbf{\hat{u}}^{(i)}}{\partial \mathrm{\hat{r}}_{n+\alpha_f}} & \mathrm{K}^{(i)}_{\hat{p}\hat{r}} &= \alpha_f \frac{\partial \mathrm{R}_{\hat{p}}^{(i)}}{\partial \mathrm{\hat{r}}_{n+\alpha_f}} \\
\mathrm{K}^{(i)}_{\mathbf{B}\hat{r}} &= \alpha_f \frac{\partial \mathrm{R}_\mathbf{B}^{(i)}}{\partial \mathrm{\hat{r}}_{n+\alpha_f}} & \mathrm{K}^{(i)}_{r\hat{r}} &= \alpha_f \frac{\partial \mathrm{R}_r^{(i)}}{\partial \mathrm{\hat{r}}_{n+\alpha_f}} &
\mathrm{K}^{(i)}_{\mathbf{\hat{B}}\hat{r}} &= \alpha_f \frac{\partial \mathrm{R}_\mathbf{\hat{B}}^{(i)}}{\partial \mathrm{\hat{r}}_{n+\alpha_f}} & \mathrm{K}^{(i)}_{\hat{r}\hat{r}} &= \alpha_f \frac{\partial \mathrm{R}_{\hat{r}}^{(i)}}{\partial \mathrm{\hat{r}}_{n+\alpha_f}}
\end{align}
are the corresponding tangent matrix entries.\\

\textit{Step 3:} Solve the linear system assembled in Step 2 for the update vectors $\Delta \mathrm{\dot{U}}^{(i)}_{n+1}$, $\Delta \mathrm{P}^{(i)}_{n+1}$, $\Delta \mathrm{\hat{U}}^{(i)}_{n+1}$, $\Delta \mathrm{\hat{P}}^{(i)}_{n+1}$, $\Delta \mathrm{\dot{B}}^{(i)}_{n+1}$, $\Delta \mathrm{R}^{(i)}_{n+1}$, $\Delta \mathrm{\hat{B}}^{(i)}_{n+1}$, and $\Delta \mathrm{\hat{R}}^{(i)}_{n+1}$, and update the iterates using the relations:
\begin{align}
\mathrm{\dot{U}}^{(i)}_{n+1} &= \mathrm{\dot{U}}^{(i-1)}_{n+1} + \Delta \mathrm{\dot{U}}^{(i)}_{n+1} & \mathrm{\dot{B}}^{(i)}_{n+1} &= \mathrm{\dot{B}}^{(i-1)}_{n+1} + \Delta \mathrm{\dot{B}}^{(i)}_{n+1} \\
\mathrm{U}^{(i)}_{n+1} &= \mathrm{U}^{(i-1)}_{n+1} + \gamma \Delta t_n \Delta \mathrm{\dot{U}}^{(i)}_{n+1} & \mathrm{B}^{(i)}_{n+1} &= \mathrm{B}^{(i-1)}_{n+1} + \Delta t_n \Delta \mathrm{\dot{B}}^{(i)}_{n+1} \\
\mathrm{P}^{(i)}_{n+1} &= \mathrm{P}^{(i-1)}_{n+1} + \Delta \mathrm{P}^{(i)}_{n+1} & \mathrm{R}^{(i)}_{n+1} &= \mathrm{R}^{(i-1)}_{n+1} + \Delta \mathrm{R}^{(i)}_{n+1} \\
\mathrm{\hat{U}}^{(i)}_{n+1} &= \mathrm{\hat{U}}^{(i-1)}_{n+1} + \Delta \mathrm{\hat{U}}^{(i)}_{n+1} & \mathrm{\hat{B}}^{(i)}_{n+1} &= \mathrm{\hat{B}}^{(i-1)}_{n+1} + \Delta \mathrm{\hat{B}}^{(i)}_{n+1} \\
\mathrm{\hat{P}}^{(i)}_{n+1} &= \mathrm{\hat{P}}^{(i-1)}_{n+1} + \Delta \mathrm{\hat{P}}^{(i)}_{n+1} & \mathrm{\hat{R}}^{(i)}_{n+1} &= \mathrm{\hat{R}}^{(i-1)}_{n+1} + \Delta \mathrm{\hat{R}}^{(i)}_{n+1}.
\label{update_equations}
\end{align}

\textit{Step 4:} If the norm of $\mathrm{R}^{(i)}_\text{total}$ is less than some prescribed value, terminate the iterative loop and set $\mathrm{\dot{U}}_{n+1} = \mathrm{\dot{U}}^{(i)}_{n+1}$, $\mathrm{U}_{n+1} = \mathrm{U}^{(i)}_{n+1}$, $\mathrm{P}_{n+1} = \mathrm{P}^{(i)}_{n+1}$, $\mathrm{\hat{U}}_{n+1} = \mathrm{\hat{U}}^{(i)}_{n+1}$, $\mathrm{\hat{P}}_{n+1} = \mathrm{\hat{P}}^{(i)}_{n+1}$, $\mathrm{\dot{B}}_{n+1} = \mathrm{\dot{B}}^{(i)}_{n+1}$, $\mathrm{B}_{n+1} = \mathrm{B}^{(i)}_{n+1}$,  $\mathrm{R}_{n+1} = \mathrm{R}^{(i)}_{n+1}$, $\mathrm{\hat{B}}_{n+1} = \mathrm{\hat{B}}^{(i)}_{n+1}$, and $\mathrm{\hat{R}}_{n+1} = \mathrm{\hat{R}}^{(i)}_{n+1}$.\\

The linear system given by \eqref{eq:gen_alpha_system} is large and expensive to both assemble and solve, so in our later numerical experiments, we remove the coupling terms between the flow field variables and the magnetic field variables in the tangent matrix. That is, we replace the linear system given by \eqref{eq:gen_alpha_system} with the following block diagonal linear system:
\begin{equation}
\left[
\begin{array}{cccccccc}
\mathrm{K}^{(i)}_{\mathbf{u}\mathbf{u}} & \mathrm{K}^{(i)}_{\mathbf{u}p} & \mathrm{K}^{(i)}_{\mathbf{u}\mathbf{\hat{u}}} & \mathrm{K}^{(i)}_{\mathbf{u}\hat{p}} & \mathrm{0} & \mathrm{0} & \mathrm{0} & \mathrm{0} \\

\mathrm{K}^{(i)}_{p\mathbf{u}} & \mathrm{K}^{(i)}_{pp} & \mathrm{K}^{(i)}_{p\mathbf{\hat{u}}} & \mathrm{K}^{(i)}_{p\hat{p}} & \mathrm{0} & \mathrm{0} & \mathrm{0} & \mathrm{0} \\

\mathrm{K}^{(i)}_{\mathbf{\hat{u}}\mathbf{u}} & \mathrm{K}^{(i)}_{\mathbf{\hat{u}}p} & \mathrm{K}^{(i)}_{\mathbf{\hat{u}}\mathbf{\hat{u}}} & \mathrm{K}^{(i)}_{\mathbf{\hat{u}}\hat{p}} & \mathrm{0} & \mathrm{0} & \mathrm{0} & \mathrm{0} \\

\mathrm{K}^{(i)}_{\hat{p}\mathbf{u}} & \mathrm{K}^{(i)}_{\hat{p}p} & \mathrm{K}^{(i)}_{\hat{p}\mathbf{\hat{u}}} & \mathrm{K}^{(i)}_{\hat{p}\hat{p}} & \mathrm{0} & \mathrm{0} & \mathrm{0} & \mathrm{0} \\

\mathrm{0} & \mathrm{0} & \mathrm{0} & \mathrm{0} & \mathrm{K}^{(i)}_{\mathbf{B}\mathbf{B}} & \mathrm{K}^{(i)}_{\mathbf{B}r} & \mathrm{K}^{(i)}_{\mathbf{B}\mathbf{\hat{B}}} & \mathrm{K}^{(i)}_{\mathbf{B}\hat{r}} \\

\mathrm{0} & \mathrm{0} & \mathrm{0} & \mathrm{0} & \mathrm{K}^{(i)}_{r\mathbf{B}} & \mathrm{K}^{(i)}_{rr} & \mathrm{K}^{(i)}_{r\mathbf{\hat{B}}} & \mathrm{K}^{(i)}_{r\hat{r}} \\

\mathrm{0} & \mathrm{0} & \mathrm{0} & \mathrm{0} & \mathrm{K}^{(i)}_{\mathbf{\hat{B}}\mathbf{B}} & \mathrm{K}^{(i)}_{\mathbf{\hat{B}}r} & \mathrm{K}^{(i)}_{\mathbf{\hat{B}}\mathbf{\hat{B}}} & \mathrm{K}^{(i)}_{\mathbf{\hat{B}}\hat{r}} \\

\mathrm{0} & \mathrm{0} & \mathrm{0} & \mathrm{0} & \mathrm{K}^{(i)}_{\hat{r}\mathbf{B}} & \mathrm{K}^{(i)}_{\hat{r}r} & \mathrm{K}^{(i)}_{\hat{r}\mathbf{\hat{B}}} & \mathrm{K}^{(i)}_{\hat{r}\hat{r}}

\end{array}
\right] \left[
\begin{array}{c}
\Delta \mathrm{\dot{U}}^{(i)}_{n+1} \\
\Delta \mathrm{P}^{(i)}_{n+1} \\
\Delta \mathrm{\hat{U}}^{(i)}_{n+1} \\
\Delta \mathrm{\hat{P}}^{(i)}_{n+1} \\
\Delta \mathrm{\dot{B}}^{(i)}_{n+1} \\
\Delta \mathrm{R}^{(i)}_{n+1} \\
\Delta \mathrm{\hat{B}}^{(i)}_{n+1} \\
\Delta \mathrm{\hat{R}}^{(i)}_{n+1}
\end{array}
\right] = -\left[
\begin{array}{c}
\mathrm{R}_\mathbf{u}^{(i)} \\
\mathrm{R}_p^{(i)} \\
\mathrm{R}_\mathbf{\hat{u}}^{(i)} \\
\mathrm{R}_{\hat{p}}^{(i)} \\
\mathrm{R}_\mathbf{B}^{(i)} \\
\mathrm{R}_r^{(i)} \\
\mathrm{R}_\mathbf{\hat{B}}^{(i)} \\
\mathrm{R}_{\hat{r}}^{(i)}
\end{array}
\right].
\label{eq:alt_gen_alpha_system}
\end{equation}
Note that \eqref{eq:alt_gen_alpha_system} is comprised of two smaller decoupled linear systems, a system for the flow field variable updates,
\begin{equation}
\left[
\begin{array}{cccc}
\mathrm{K}^{(i)}_{\mathbf{u}\mathbf{u}} & \mathrm{K}^{(i)}_{\mathbf{u}p} & \mathrm{K}^{(i)}_{\mathbf{u}\mathbf{\hat{u}}} & \mathrm{K}^{(i)}_{\mathbf{u}\hat{p}} \\

\mathrm{K}^{(i)}_{p\mathbf{u}} & \mathrm{K}^{(i)}_{pp} & \mathrm{K}^{(i)}_{p\mathbf{\hat{u}}} & \mathrm{K}^{(i)}_{p\hat{p}} \\

\mathrm{K}^{(i)}_{\mathbf{\hat{u}}\mathbf{u}} & \mathrm{K}^{(i)}_{\mathbf{\hat{u}}p} & \mathrm{K}^{(i)}_{\mathbf{\hat{u}}\mathbf{\hat{u}}} & \mathrm{K}^{(i)}_{\mathbf{\hat{u}}\hat{p}} \\

\mathrm{K}^{(i)}_{\hat{p}\mathbf{u}} & \mathrm{K}^{(i)}_{\hat{p}p} & \mathrm{K}^{(i)}_{\hat{p}\mathbf{\hat{u}}} & \mathrm{K}^{(i)}_{\hat{p}\hat{p}}

\end{array}
\right] \left[
\begin{array}{c}
\Delta \mathrm{\dot{U}}^{(i)}_{n+1} \\
\Delta \mathrm{P}^{(i)}_{n+1} \\
\Delta \mathrm{\hat{U}}^{(i)}_{n+1} \\
\Delta \mathrm{\hat{P}}^{(i)}_{n+1}
\end{array}
\right] = -\left[
\begin{array}{c}
\mathrm{R}_\mathbf{u}^{(i)} \\
\mathrm{R}_p^{(i)} \\
\mathrm{R}_\mathbf{\hat{u}}^{(i)} \\
\mathrm{R}_{\hat{p}}^{(i)}
\end{array}
\right],
\label{eq:alt_gen_alpha_system_flow}
\end{equation}
and a system for the magnetic field variable updates,
\begin{equation}
\left[
\begin{array}{cccc}
\mathrm{K}^{(i)}_{\mathbf{B}\mathbf{B}} & \mathrm{K}^{(i)}_{\mathbf{B}r} & \mathrm{K}^{(i)}_{\mathbf{B}\mathbf{\hat{B}}} & \mathrm{K}^{(i)}_{\mathbf{B}\hat{r}} \\

\mathrm{K}^{(i)}_{r\mathbf{B}} & \mathrm{K}^{(i)}_{rr} & \mathrm{K}^{(i)}_{r\mathbf{\hat{B}}} & \mathrm{K}^{(i)}_{r\hat{r}} \\

\mathrm{K}^{(i)}_{\mathbf{\hat{B}}\mathbf{B}} & \mathrm{K}^{(i)}_{\mathbf{\hat{B}}r} & \mathrm{K}^{(i)}_{\mathbf{\hat{B}}\mathbf{\hat{B}}} & \mathrm{K}^{(i)}_{\mathbf{\hat{B}}\hat{r}} \\

\mathrm{K}^{(i)}_{\hat{r}\mathbf{B}} & \mathrm{K}^{(i)}_{\hat{r}r} & \mathrm{K}^{(i)}_{\hat{r}\mathbf{\hat{B}}} & \mathrm{K}^{(i)}_{\hat{r}\hat{r}}

\end{array}
\right] \left[
\begin{array}{c}
\Delta \mathrm{\dot{B}}^{(i)}_{n+1} \\
\Delta \mathrm{R}^{(i)}_{n+1} \\
\Delta \mathrm{\hat{B}}^{(i)}_{n+1} \\
\Delta \mathrm{\hat{R}}^{(i)}_{n+1}
\end{array}
\right] = -\left[
\begin{array}{c}
\mathrm{R}_\mathbf{B}^{(i)} \\
\mathrm{R}_r^{(i)} \\
\mathrm{R}_\mathbf{\hat{B}}^{(i)} \\
\mathrm{R}_{\hat{r}}^{(i)}
\end{array}
\right],
\label{eq:alt_gen_alpha_system_magnetic}
\end{equation}
which can be assembled and solved in parallel. As such, we refer to the predictor-multicorrector method that employs \eqref{eq:alt_gen_alpha_system} rather than \eqref{eq:gen_alpha_system} as a block iterative predictor-multicorrector method. Alternatively, we refer to the predictor-multicorrector method that employs \eqref{eq:gen_alpha_system} as a fully coupled predictor-multicorrector method as it maintains the coupling between the flow field variables and the magnetic flow field variables in the tangent matrix. According to the terminology used in \cite{cervera1996computational}, our block iterative predictor-multicorrector method is of Jacobi type.

\section{Static Condensation}

To solve \eqref{eq:alt_gen_alpha_system}, \eqref{eq:alt_gen_alpha_system_flow}, or \eqref{eq:alt_gen_alpha_system_magnetic}, we can statically condense the interior degree-of-freedom updates to arrive at a smaller linear system for only the trace degree-of-freedom updates. To see this, consider the linear system given by \eqref{eq:alt_gen_alpha_system_magnetic}. We can write this system as
\begin{equation}
\left[
\begin{array}{cc}
\mathrm{K}_{\text{int},\text{int}} & \mathrm{K}_{\text{int},\text{trace}} \\
\mathrm{K}_{\text{trace},\text{int}} & \mathrm{K}_{\text{trace},\text{trace}}
\end{array}
\right] \left[
\begin{array}{c}
\mathrm{d}_\text{int} \\
\mathrm{d}_\text{trace}
\end{array}
\right] = \left[
\begin{array}{c}
\mathrm{F}_\text{int} \\
\mathrm{F}_\text{trace}
\end{array}
\right]
\end{equation}
where
\begin{align}
\mathrm{K}_{\text{int},\text{int}} &= \left[
\begin{array}{cc}
\mathrm{K}^{(i)}_{\mathbf{B}\mathbf{B}} & \mathrm{K}^{(i)}_{\mathbf{B}r} \\
\mathrm{K}^{(i)}_{r\mathbf{B}} & \mathrm{K}^{(i)}_{rr}
\end{array}
\right] & \mathrm{K}_{\text{int},\text{trace}}
 &= \left[
\begin{array}{cc}
\mathrm{K}^{(i)}_{\mathbf{B}\mathbf{\hat{B}}} & \mathrm{K}^{(i)}_{\mathbf{B}\hat{r}} \\
\mathrm{K}^{(i)}_{r\mathbf{\hat{B}}} & \mathrm{K}^{(i)}_{r\hat{r}}
\end{array}
\right] \\
\mathrm{K}_{\text{trace},\text{int}} &= \left[
\begin{array}{cc}
\mathrm{K}^{(i)}_{\mathbf{\hat{B}}\mathbf{B}} & \mathrm{K}^{(i)}_{\mathbf{\hat{B}}r} \\
\mathrm{K}^{(i)}_{\hat{r}\mathbf{B}} & \mathrm{K}^{(i)}_{\hat{r}r}
\end{array}
\right] & \mathrm{K}_{\text{trace},\text{trace}}
 &= \left[
\begin{array}{cc}
\mathrm{K}^{(i)}_{\mathbf{\hat{B}}\mathbf{\hat{B}}} & \mathrm{K}^{(i)}_{\mathbf{\hat{B}}\hat{r}} \\
\mathrm{K}^{(i)}_{\hat{r}\mathbf{\hat{B}}} & \mathrm{K}^{(i)}_{\hat{r}\hat{r}}
\end{array}
\right] \\
\mathrm{d}_{\text{int}} &= \left[
\begin{array}{c}
\Delta \mathrm{\dot{B}}^{(i)}_{n+1} \\
\Delta \mathrm{R}^{(i)}_{n+1}
\end{array}
\right] & \mathrm{d}_{\text{trace}}
 &= \left[
\begin{array}{c}
\Delta \mathrm{\hat{B}}^{(i)}_{n+1} \\
\Delta \mathrm{\hat{R}}^{(i)}_{n+1}
\end{array}
\right] \\
\mathrm{F}_{\text{int}} &= -\left[
\begin{array}{c}
\mathrm{R}_\mathbf{B}^{(i)} \\
\mathrm{R}_r^{(i)}
\end{array}
\right] & \mathrm{F}_{\text{trace}}
 &= -\left[
\begin{array}{c}
\mathrm{R}_\mathbf{\hat{B}}^{(i)} \\
\mathrm{R}_{\hat{r}}^{(i)}.
\end{array}
\right]
\end{align}
Static condensation of $\mathrm{d}_{\text{int}}$ yields the smaller linear system
\begin{equation}
    \mathrm{S}_{\text{trace},\text{trace}} \mathrm{d}_\text{trace} =  \mathrm{F}_\text{trace} - \mathrm{K}_{\text{trace},\text{int}} \left( \mathrm{K}_{\text{int},\text{int}} \right)^{-1} \mathrm{F}_\text{int}
\end{equation}
where $\mathrm{S}_{\text{trace},\text{trace}}$ is the Schur complement
\begin{equation}
    \mathrm{S}_{\text{trace},\text{trace}} := \mathrm{K}_{\text{trace},\text{trace}} - \mathrm{K}_{\text{trace},\text{int}} \left( \mathrm{K}_{\text{int},\text{int}} \right)^{-1} \mathrm{K}_{\text{int},\text{trace}}.
\end{equation}
Since the interior degrees-of-freedom are local to each element, the matrix $\mathrm{K}_{\text{int},\text{int}}$ can be inverted in an element-by-element manner and the Schur complement $\mathrm{S}_{\text{trace},\text{trace}}$ has the same sparsity structure as $\mathrm{K}_{\text{trace},\text{trace}}$. In fact, the Schur complement can be formed and assembled element-wise in a standard element assembly routine as discussed in \cite{kirby2012cg}. Once the trace degree-of-freedom updates are attained, the interior degree-of-freedom updates can be recovered using the linear system
\begin{equation}
    \mathrm{K}_{\text{int},\text{int}} \mathrm{d}_\text{int} =  \mathrm{F}_\text{int} - \mathrm{K}_{\text{int},\text{trace}} \mathrm{d}_\text{trace}
\end{equation}
which can be solved in an element-by-element manner.

\section{Numerical Results}

Now that we have presented our HDG method for solving the incompressible MHD equations, we conduct a sequence of numerical experiments to verify the method. We first analyze the spatial accuracy and pressure robustness of our method with a manufactured solution. We then analyze the spatial accuracy, temporal accuracy, and energy stability of our method by considering the following benchmark problems: Hartmann channel flow, Alfv\'en wave propagation, and the Kelvin-Helmholtz instability problem.

We use time-marching in each numerical experiment, and for our steady test cases, we solve the unsteady incompressible MHD equations to steady state. We employ the backward Euler method for the first time step and the generalized-$\alpha$ method with $\rho_\infty = 0.5$ for all subsequent time steps. In each numerical experiment, we set $\rho = \mu_0 = 1$ and $C_\text{pen} = (k+1)(k+2)$, and for problems with only Dirichlet boundary conditions applied, we enforce $\int_{\Omega} p=\int_{\Omega} r=0$ as the pressure and magnetic pressure fields are unique only up to a constant. We frequently reference the Reynolds number and magnetic Reynolds number below, which are defined as $Re = \mathcal{U}\mathcal{L}/\nu$ and $Re_m = \mu_0\mathcal{U}\mathcal{L}/\eta$ where $\mathcal{U}$ and $\mathcal{L}$ are characteristic velocity and length scales.

\subsection{Two-Dimensional Manufactured Solution}
Our first numerical experiment is a steady manufactured solution. In particular, we adopt the two-dimensional manufactured vortex solution considered in Subsection 9.1 of \cite{Eric} and extend it to the setting of incompressible MHD. 
The velocity field, pressure field, magnetic field, and magnetic pressure field are set to
\begin{equation}
    \mathbf{u}=
    \begin{pmatrix}
    -2x^2e^x(-y^2+y)(2y-1)(x-1)^2 \\
    -xy^2e^x(x(x+3)-2)(x-1)(y-1)^2
    \end{pmatrix}
\end{equation}
\begin{equation}
    p=p_0\sin(\pi x)\sin(\pi y)
\end{equation}
\begin{equation}
    \mathbf{B}=
    \begin{pmatrix}
    -2x^2e^x(-y^2+y)(2y-1)(x-1)^2 \\
    -xy^2e^x(x(x+3)-2)(x-1)(y-1)^2
    \end{pmatrix}
\end{equation}
\begin{equation}
    r=\sin(\pi x)\sin(\pi y)
\end{equation}
where $p_0 \in \mathbb{R}$. The corresponding forcings in the conservation of momentum and magnetic induction equations are then
\begin{equation}
    \mathbf{f}_v=\nabla\cdot\left(\mathbf{u}\otimes\mathbf{u}\right)+\frac{1}{\rho}\nabla p-\nabla\cdot\left(2\nu\nabla^s\mathbf{u}\right)-\nabla\cdot\left(\frac{1}{\rho\mu_0}\left(\mathbf{B}\otimes\mathbf{B}-\frac{1}{2}|\mathbf{B}|^2\mathds{1}\right)\right)
\end{equation}
\begin{equation}
    \mathbf{f}_m=\nabla\cdot(\mathbf{u}\otimes\mathbf{B})-\nabla\cdot(\mathbf{B}\otimes\mathbf{u})-\nabla\cdot\left(\frac{\eta}{\mu_0}\nabla^a\mathbf{B}\right)+\nabla r.
\end{equation}
The flow domain for this experiment is set to $\Omega = [0,1]^2$,  and homogeneous Dirichlet boundary conditions are applied along all of $\partial \Omega$ for both the velocity and magnetic fields. For real problems of scientific and engineering interest, $r = 0$, but we consider $r \neq 0$ to stress our method. We further set $\nu = \eta = 10^{-2}$ which corresponds to $Re = Re_m = 10^2$ if we consider $\mathcal{U} = \mathcal{L} = 1$.

\begin{figure}[t!]
    \centering
    \includegraphics[width=0.45\textwidth]{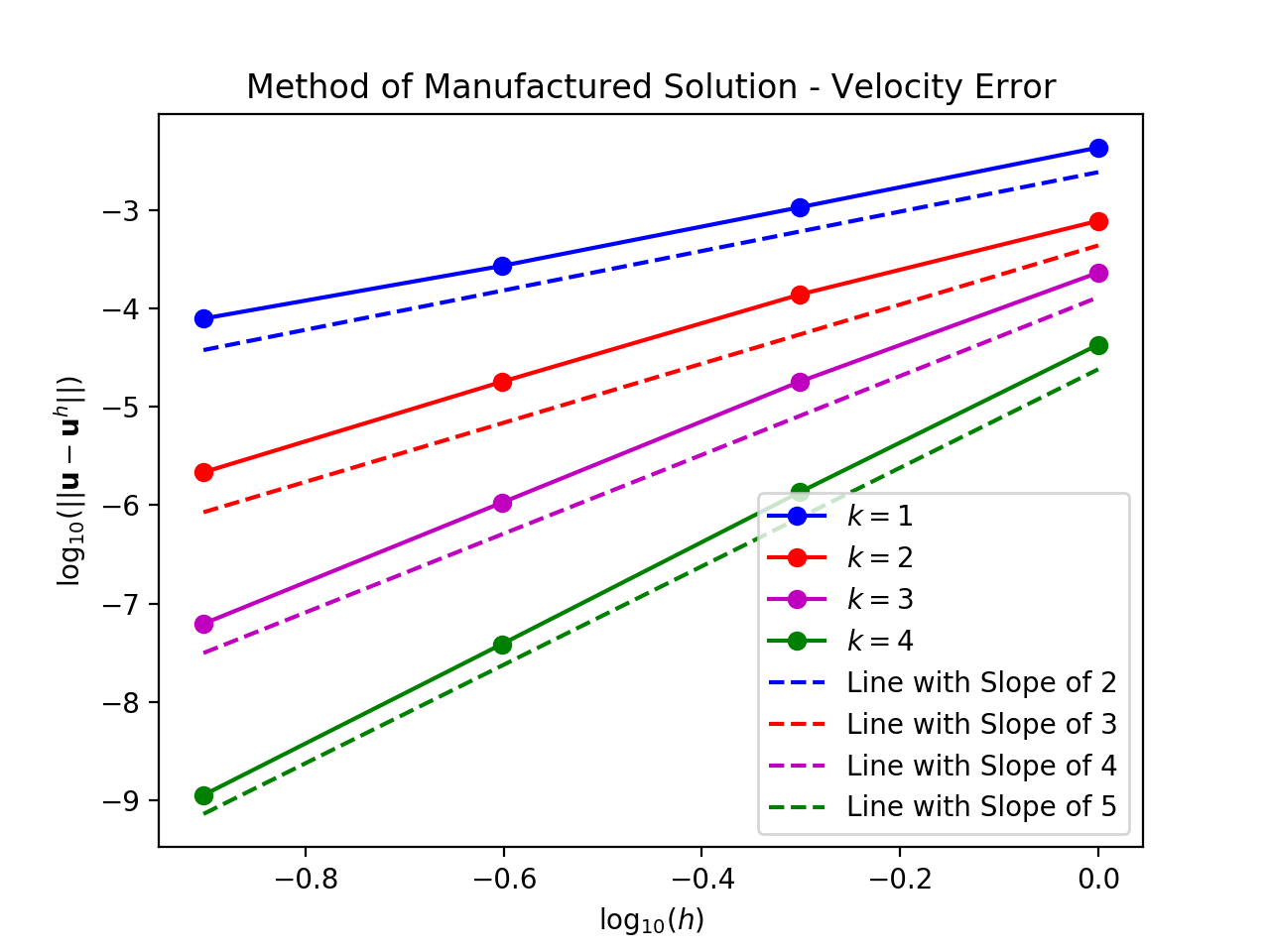}
    \includegraphics[width=0.45\textwidth]{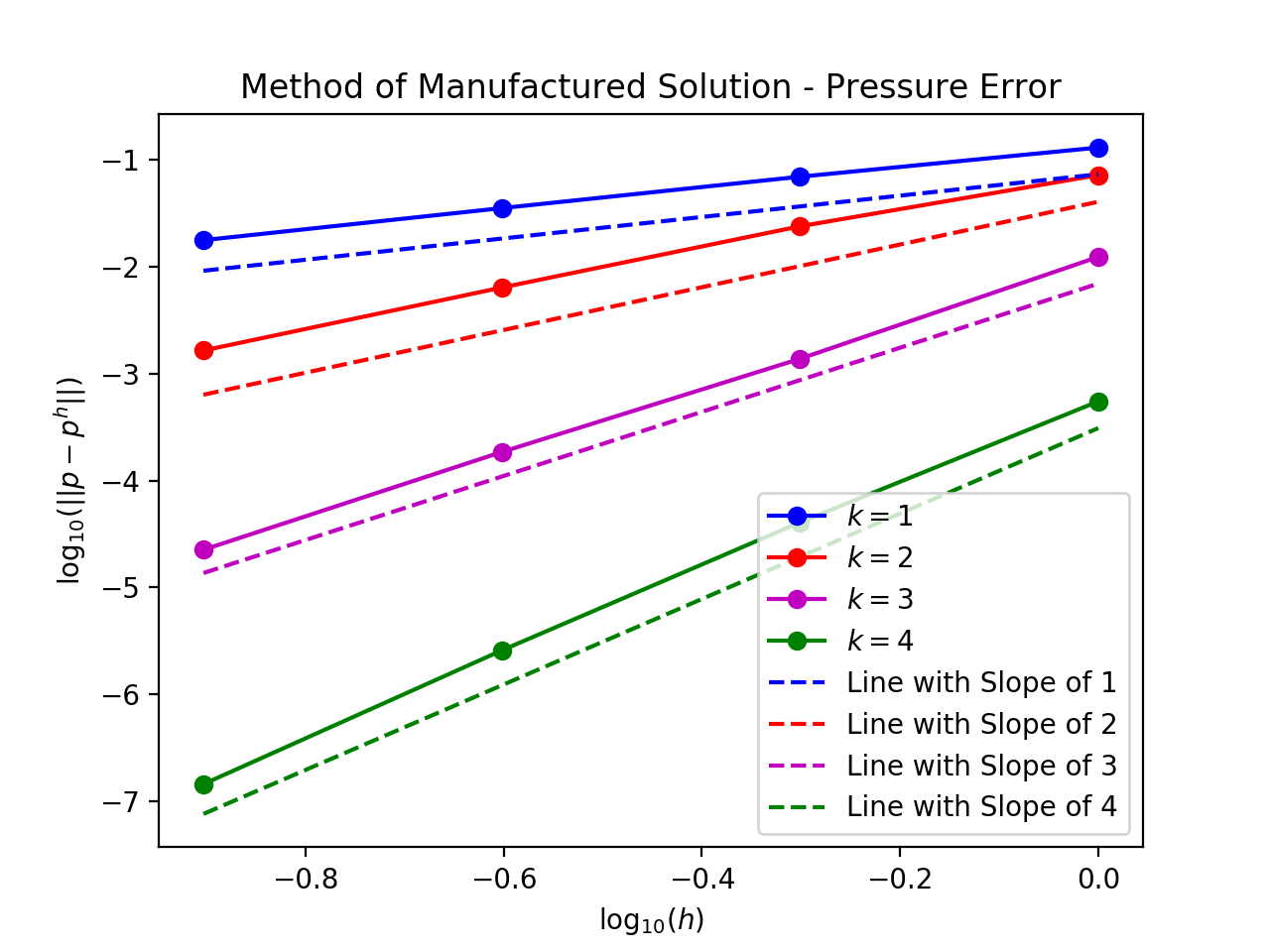}
    \includegraphics[width=0.45\textwidth]{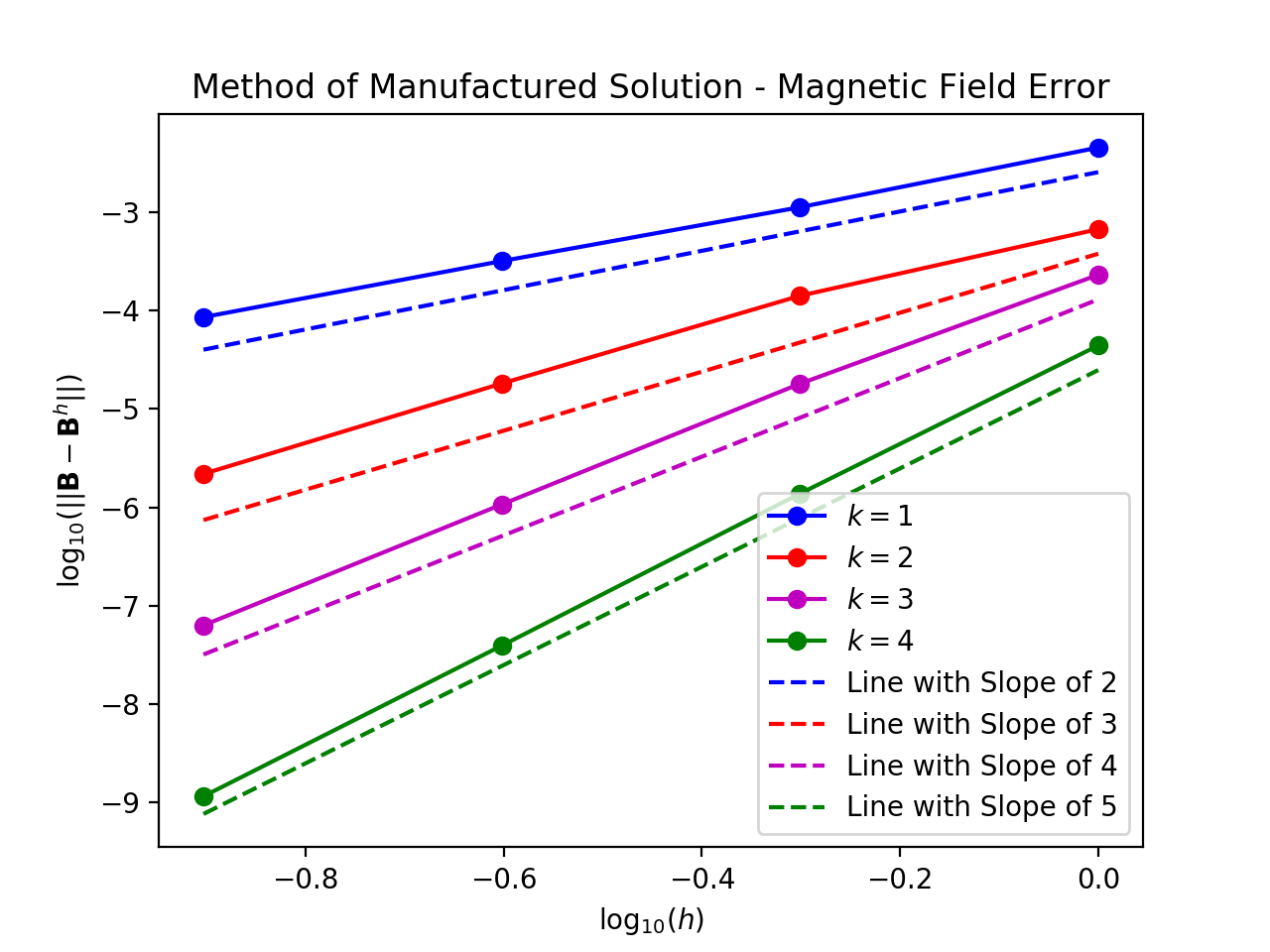}
    \includegraphics[width=0.45\textwidth]{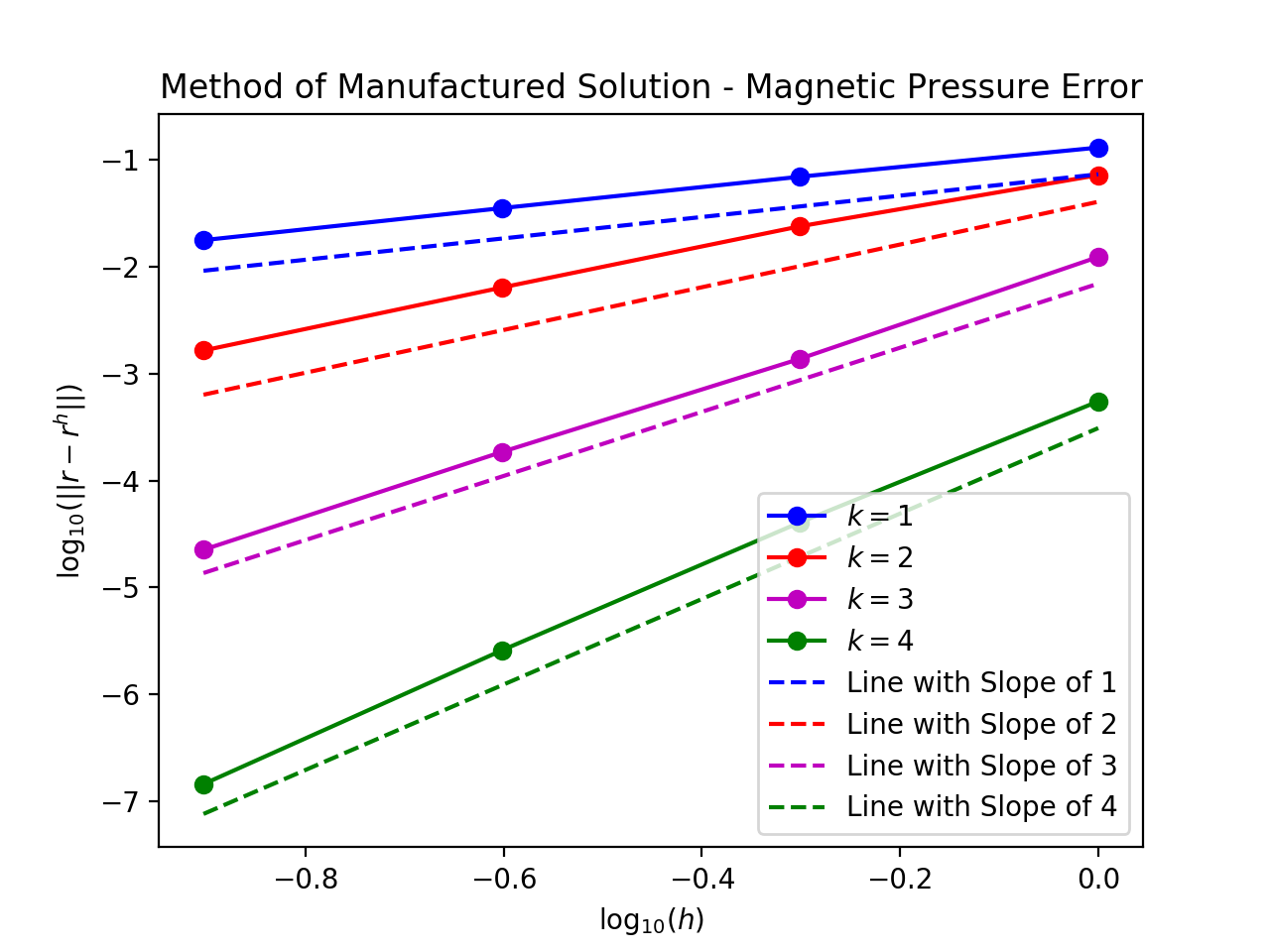}
    \caption{$L^2$-errors in the computed fields for the method of manufactured solutions problem for $p_0 = 1$.}
    \label{mms_conv}
\end{figure}

\begin{figure}[t!]
    \centering
    \includegraphics[width=0.45\textwidth]{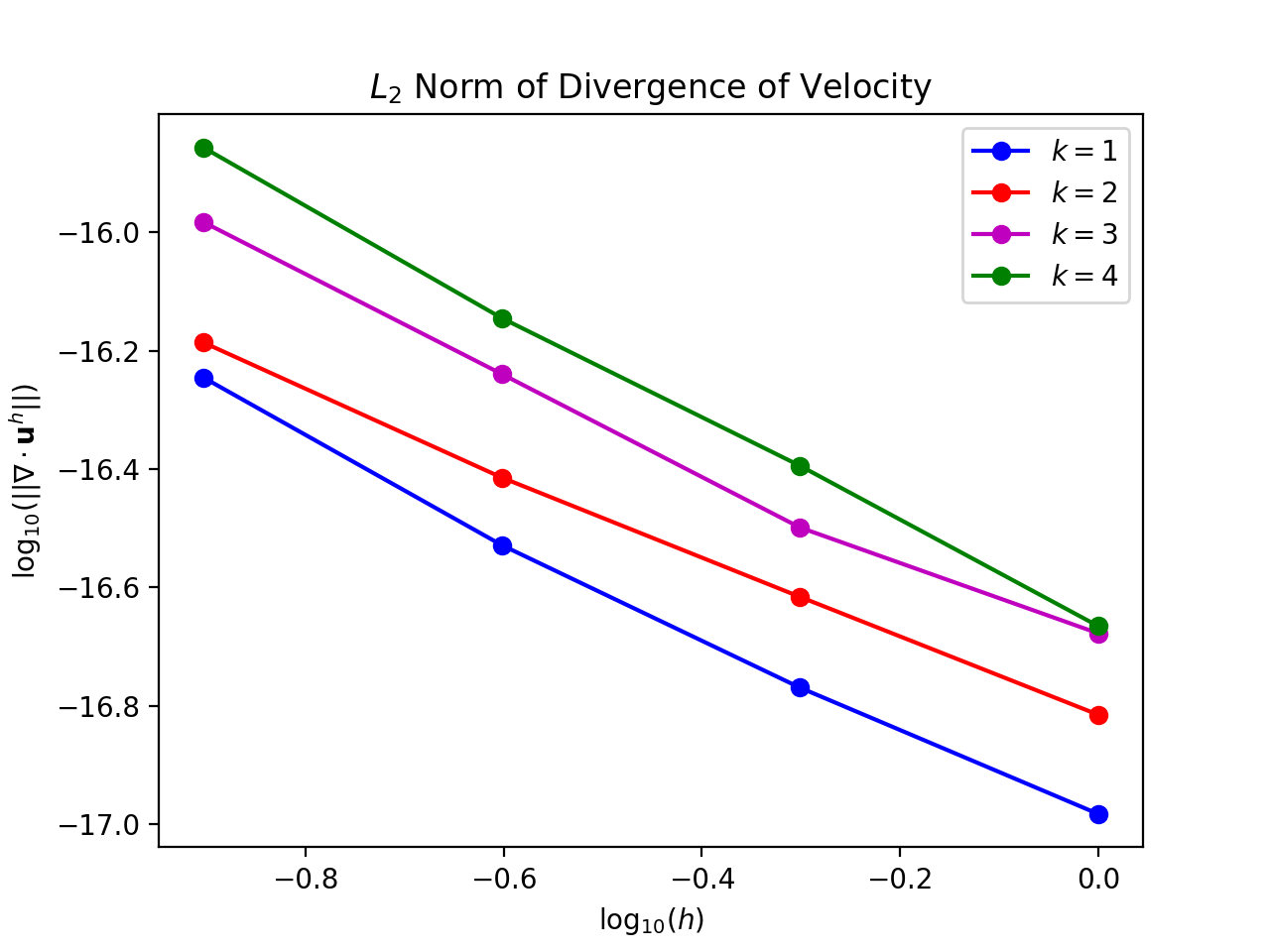}
    \includegraphics[width=0.45\textwidth]{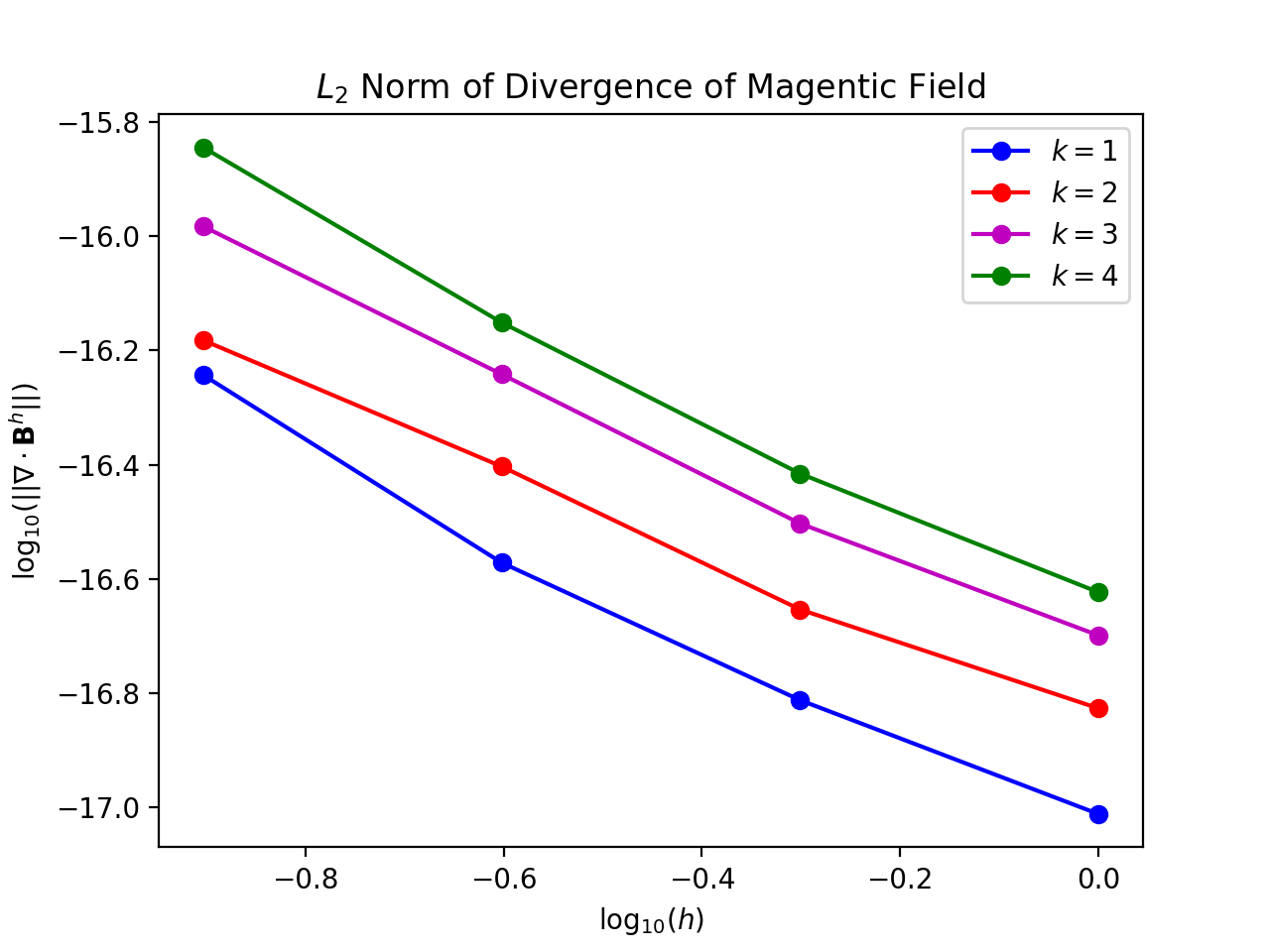}
    \caption{$L^2$-errors in the divergence of the computed velocity and magnetic fields for the method of manufactured solutions problem for $p_0 = 1$.}
    \label{diverr}
\end{figure}

To assess the spatial accuracy of our method, we have solved this manufactured solution problem for a series of meshes, polynomial degrees $k = 1, 2, 3, 4$, and $p_0 = 1$. Plots of the $L^2$-error in the computed velocity, pressure, magnetic, and magnetic pressure fields are shown in Fig. \ref{mms_conv} where $h = \max_e h_e$ denotes the global mesh size. It is apparent from the plots that optimal rates are attained for each polynomial degree. Plots of the velocity and magnetic field divergence $L^2$-errors are displayed in Fig. \ref{diverr}, and we observe that these errors are equal to zero up to machine precision for each mesh and polynomial degree.

To assess the pressure robustness of our method, we have solved the manufactured solution problem using a mesh of $1,024$ elements, a polynomial degree of $k = 2$, and a range of $p_0$ values. The results of this study are displayed in Table \ref{table:pressure_robustness}. It is clear from the table that the $L^2$-errors of the velocity and magnetic field are independent of $p_0$, suggesting that these errors indeed do not depend on the pressure field and thus our method is pressure robust. Remarkably, the $L^2$-error of the magnetic pressure field is also independent of $p_0$. We believe this is due to absence of the pressure field in the magnetic induction equation.

%%%%%%%%%%%%%%%%%%%%%%%%%%%%%%

\begin{table}[H]
\centering
\begin{tabular}{|c|c|c|c|c|c|c|}
\hline
$p_0$ & $||\mathbf{u}-\mathbf{u}^h||$ & $\|p-p^h\|$ & $||\mathbf{B}-\mathbf{B}^h||$ & $\|r-r^h\|$ & $\|\nabla\cdot\mathbf{u}^h\|$ & $\|\nabla\cdot\mathbf{B}^h\|$ \\
\hline
1 & $2.16\cdot10^{-6}$ & $1.65\cdot10^{-3}$ & $2.19\cdot10^{-6}$ & $1.65\cdot10^{-3}$ & $6.51\cdot10^{-17}$ & $6.57\cdot10^{-17}$ \\
10 & $2.16\cdot10^{-6}$ & $1.43\cdot10^{-2}$ & $2.19\cdot10^{-9}$ & $1.46\cdot10^{-3}$ & $6.60\cdot10^{-17}$ & $6.58\cdot10^{-17}$ \\
25 & $2.16\cdot10^{-6}$ & $3.58\cdot10^{-2}$ & $2.19\cdot10^{-6}$ & $1.46\cdot10^{-3}$ & $6.64\cdot10^{-17}$ & $6.57\cdot10^{-17}$ \\
100 & $2.16\cdot10^{-6}$ & $1.43\cdot10^{-1}$ & $2.19\cdot10^{-6}$ & $1.46\cdot10^{-3}$ & $6.53\cdot10^{-17}$ & $6.53\cdot10^{-17}$ \\
\hline
\end{tabular}
\caption{$L^2$-errors in the computed fields for the method of manufactured solutions problem for a mesh of $1,024$ elements, a polynomial degree of $k = 2$, and a range of $p_0$ values.}
\label{table:pressure_robustness}
\end{table}

\subsection{Hartmann Channel Flow}
We next consider Hartmann channel flow, a generalization of the classical plane Poiseuille problem to the setting of incompressible MHD. In this problem, a steady flow between horizontal plates located at $y = -L$ and $y = L$ is driven by an applied pressure gradient $\frac{\partial p}{\partial x}=-G_0$ and subject to an applied magnetic field $B_y = B_0$ in the vertical direction. The resulting velocity field $\textbf{u} = (u_x,u_y)$ and magnetic field $\textbf{B} = (B_x,B_y)$ then take the form (see, e.g., \cite{Shadid1, Shadid2})
\begin{align}
    u_x &= \frac{G_0L^2}{\rho\nu} \frac{1}{Ha}\frac{1}{\tanh(Ha)}\left(1-\frac{\cosh(Ha~y/L)}{\cosh(Ha)}\right) &
    u_y &= 0
    \label{Hartmann_ux} \\
    B_x &= -\frac{G_0L^2\mu_0}{\sqrt{\rho\nu\eta}}\frac{1}{Ha}\left(\frac{y}{L}-\frac{1}{\tanh(Ha)}\frac{\sinh(Ha~y/L)}{\cosh(Ha)}\right) &
    B_y &= B_0.
    \label{Hartmann_Bx}
\end{align}
The pressure field can be found by inserting the above into the conservation of momentum equation, solving for $\nabla p$, and integrating, yielding
\begin{equation}
p = -G_0 x - \frac{G_0^2L^4\mu_0}{2Ha^2\rho\nu\eta}\left(\frac{\sinh^2(Ha~y/L)}{\cosh^2(Ha)\tanh^2(Ha)}-\frac{2\sinh(Ha~y/L)}{\cosh(Ha)\tanh(Ha)}\frac{y}{L}+\left(\frac{y}{L}\right)^2\right).
\end{equation}
The magnetic pressure field is equal to $r = 0$. This problem is characterized by the Hartmann number $Ha = B_0L/\sqrt{\rho\nu\eta}$ which is a ratio of the electromagnetic force to the viscous force. As the Hartmann number is increased, a sharp layer in $B_x$ begins to form near the two plates. In our numerical experiments, we set $L = \nu = \eta = 1$ and $G_0 = 5$, and we consider the flow domain $\Omega=(0,0.5)\times(-1,1)$. Dirichlet boundary conditions are applied along all of $\partial \Omega$ for both the velocity and magnetic fields using the analytical solution.

\begin{figure}[t!]
    \centering
    \includegraphics[width=0.45\textwidth]{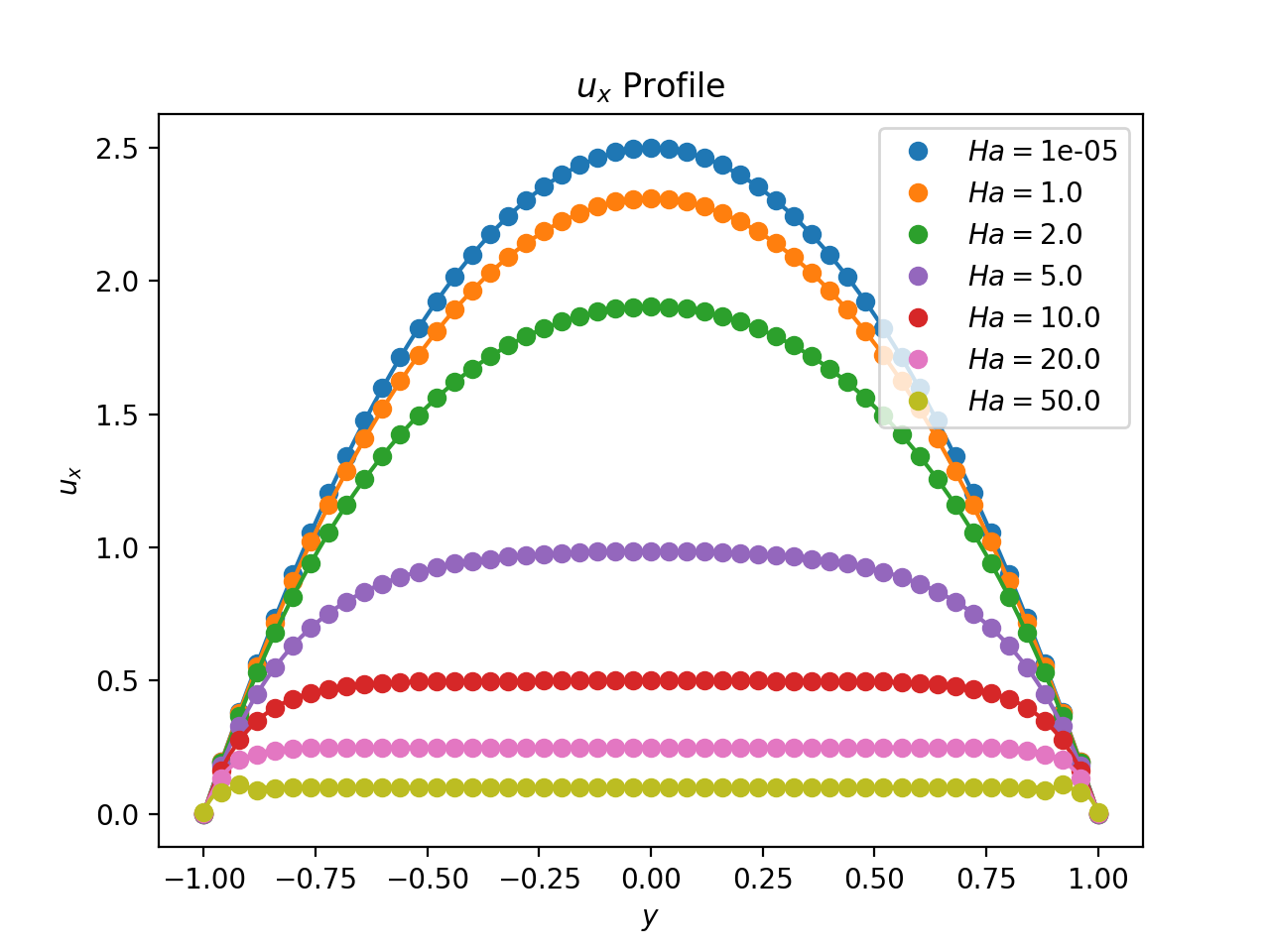}
    \includegraphics[width=0.45\textwidth]{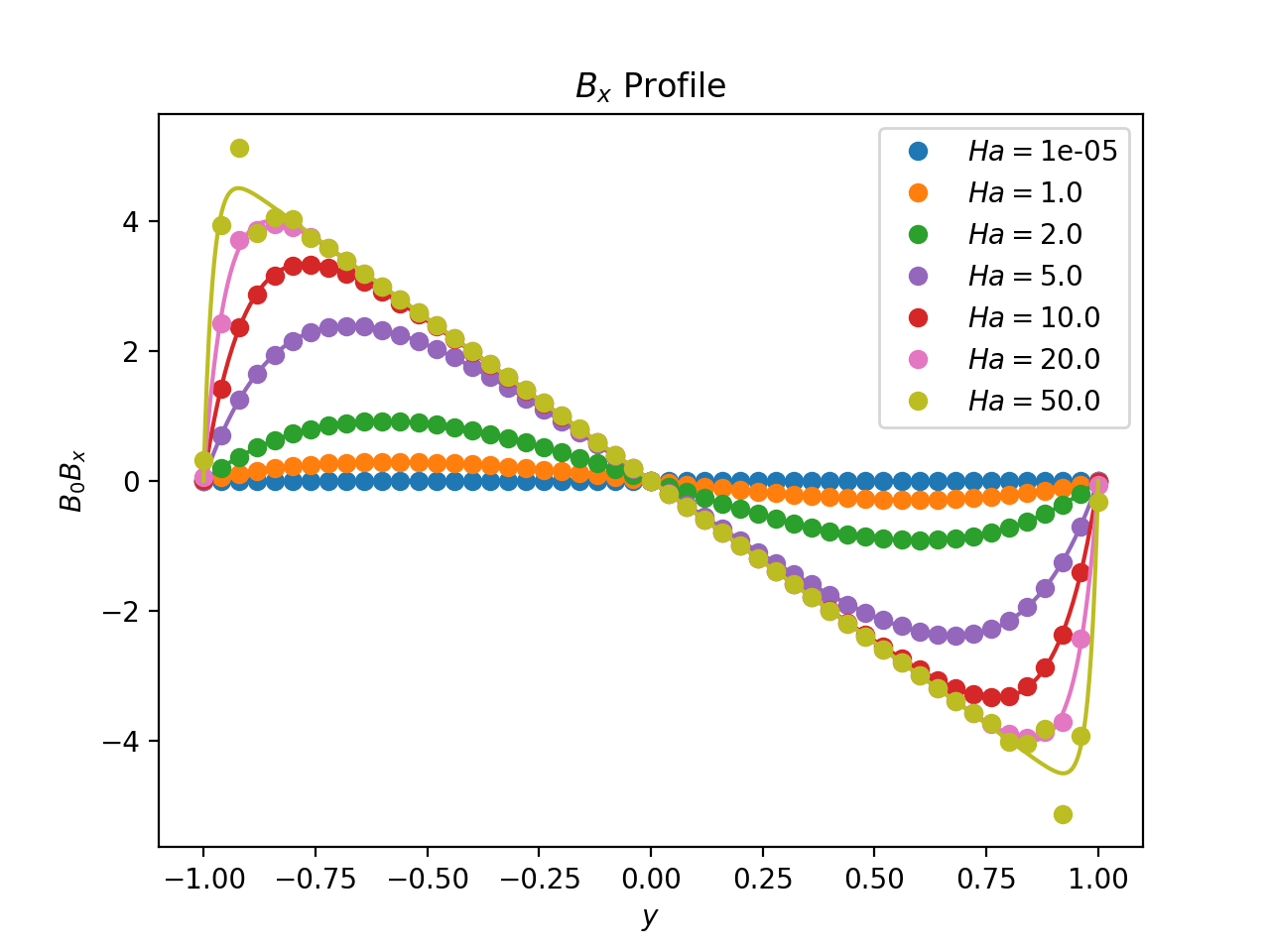} \\
    \hspace{5pt} (a) \hspace{195pt} (b) \\
    \includegraphics[width=0.45\textwidth]{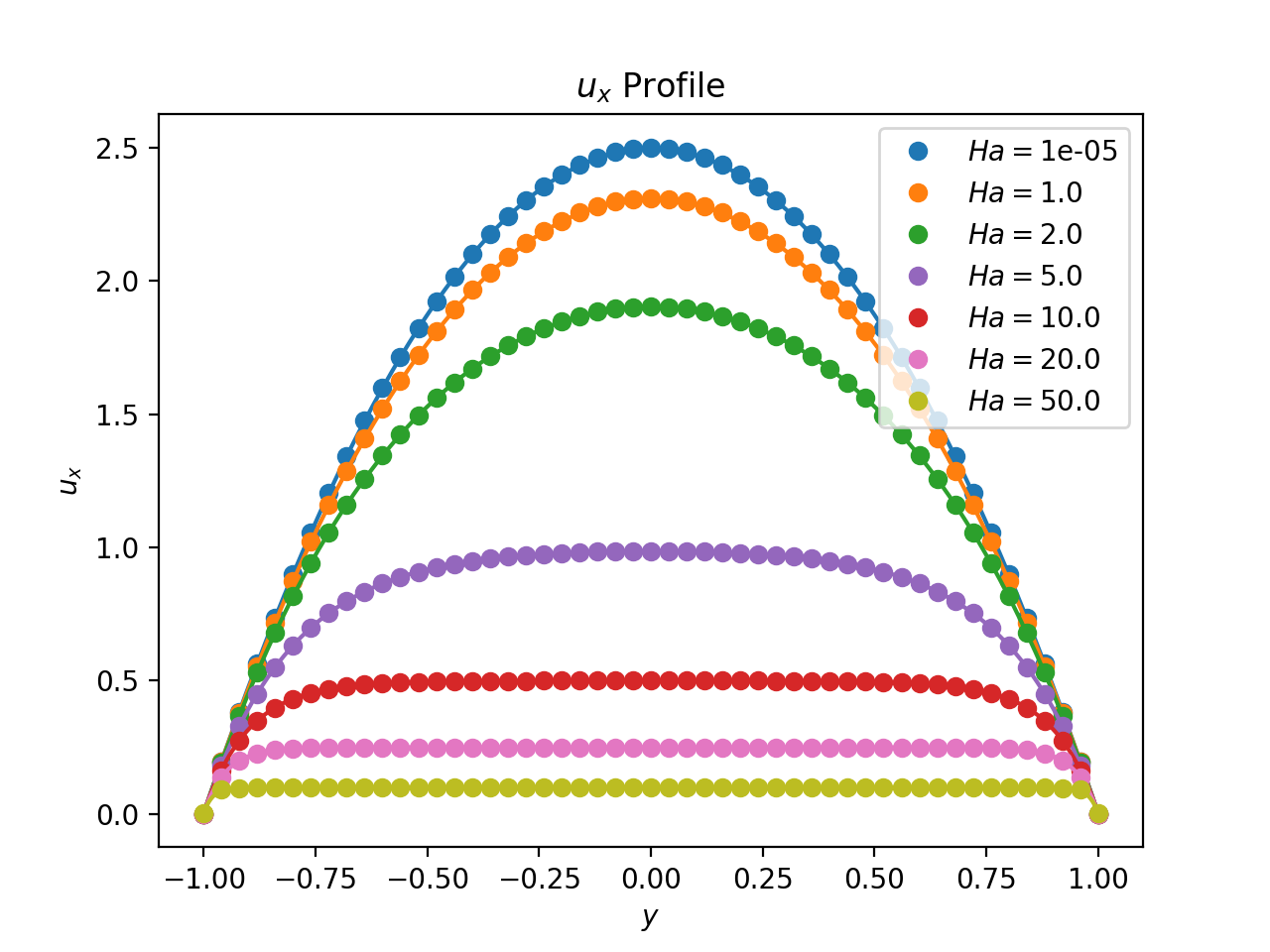}
    \includegraphics[width=0.45\textwidth]{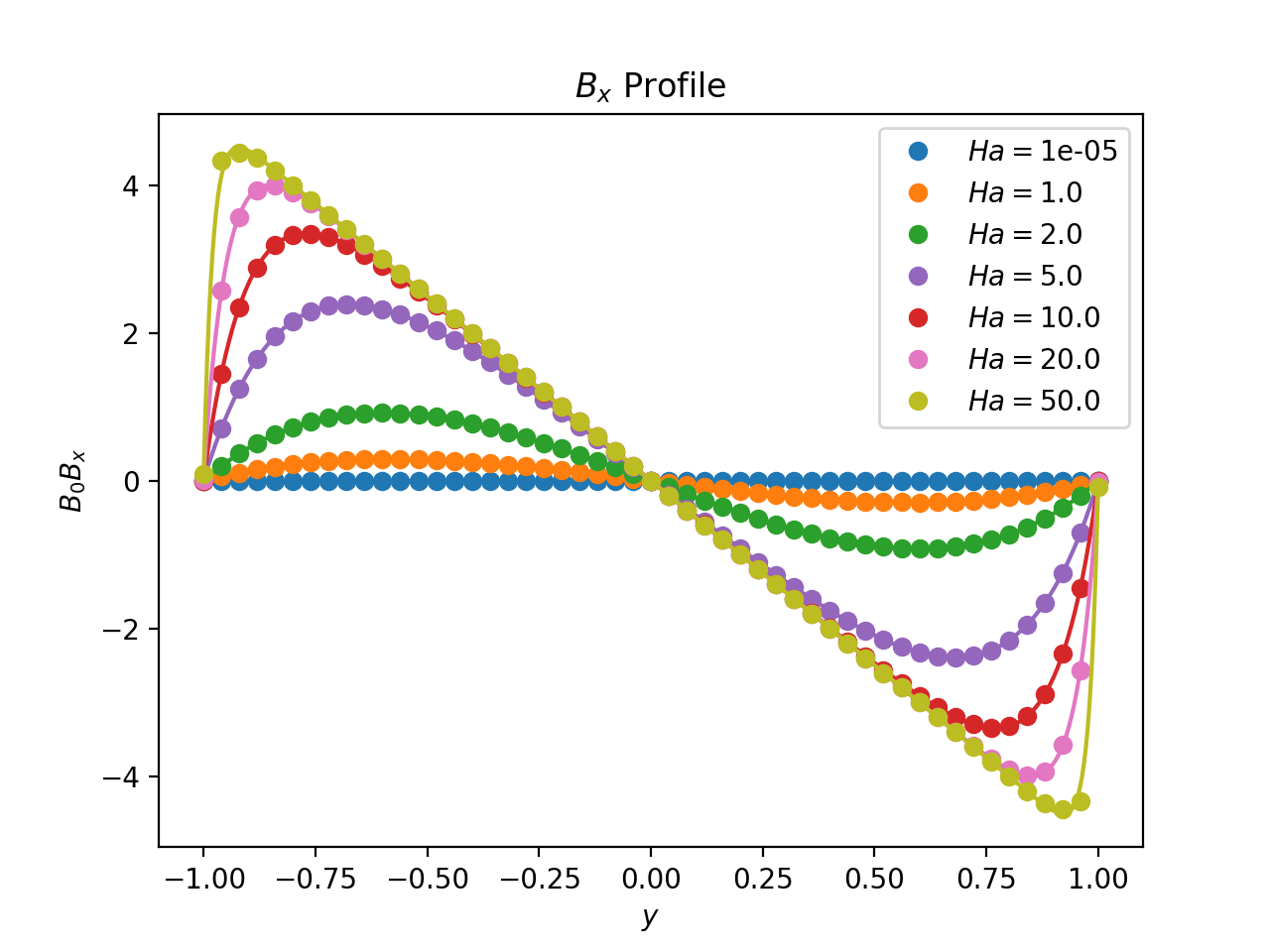} \\
    \hspace{5pt} (c) \hspace{195pt} (d) \\
    \includegraphics[width=0.45\textwidth]{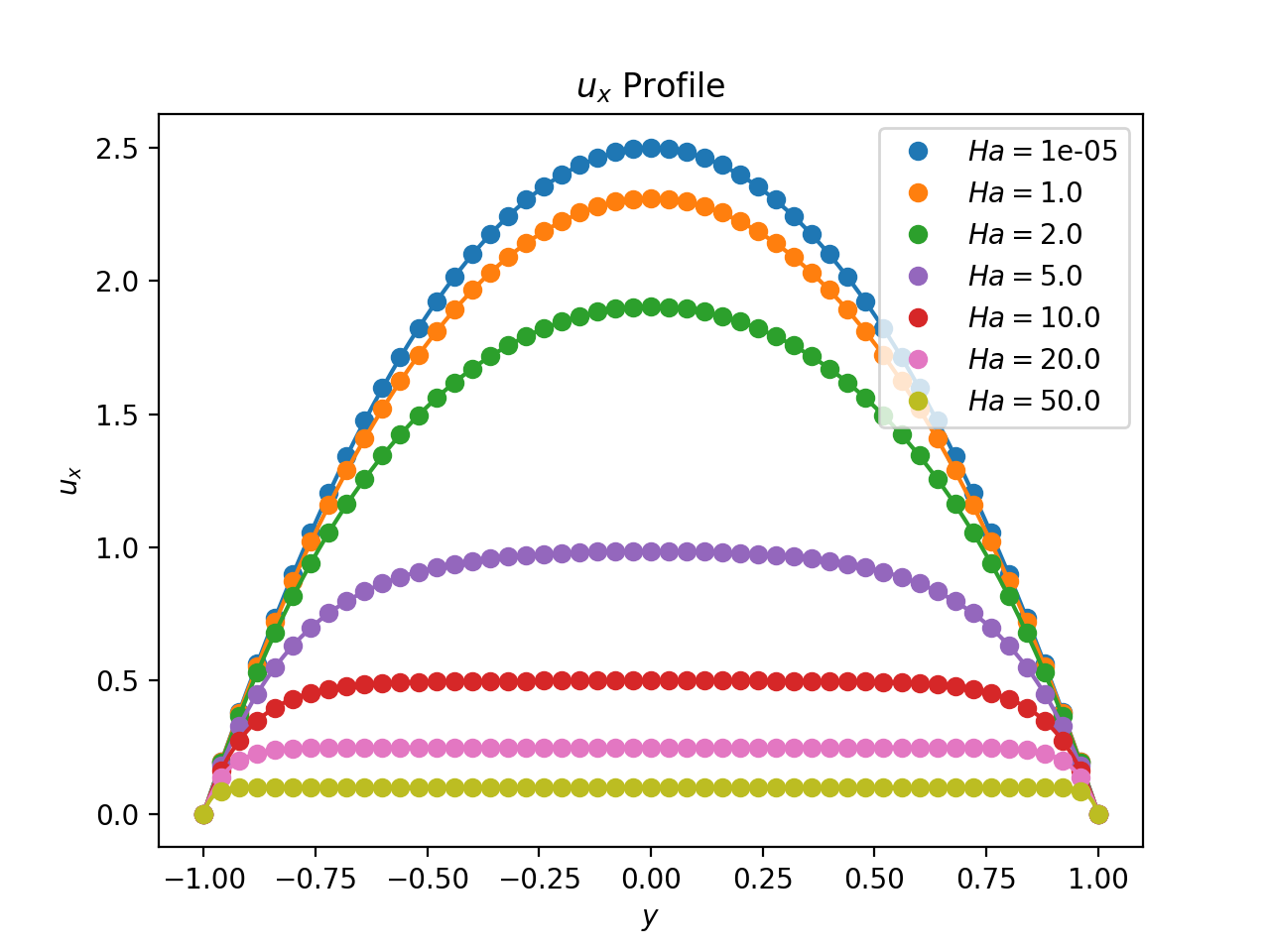}
    \includegraphics[width=0.45\textwidth]{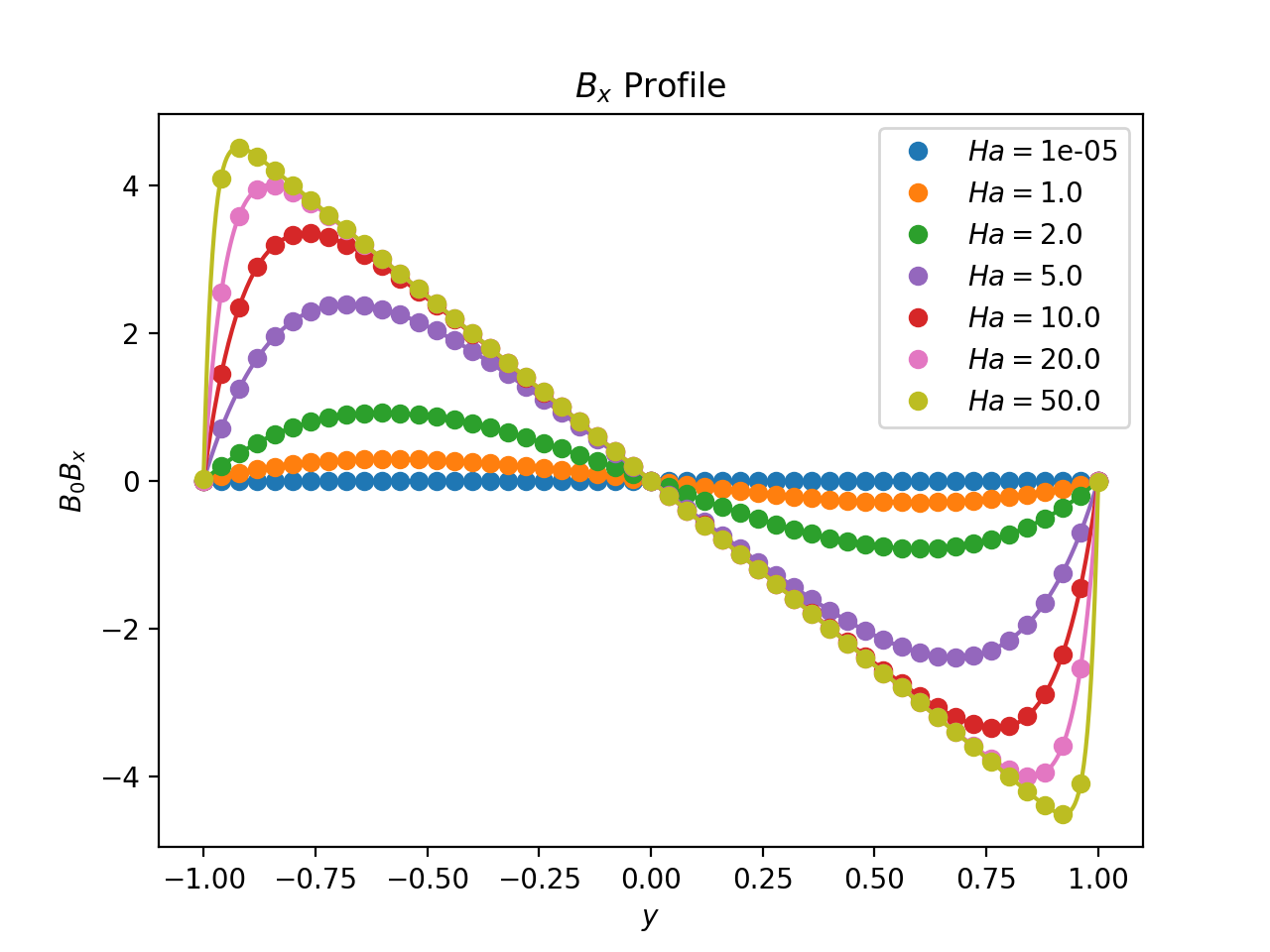} \\
    \hspace{5pt} (e) \hspace{195pt} (f)
    \caption{Computed $u_x$ and $B_x$ fields along the line $x = 0.25$ for the Hartmann channel flow problem for a variety of Hartmann numbers, a polynomial degree of $k = 2$, and meshes consisting of 128 elements (top row - (a) and (b)), 512 elements (middle row - (c) and (d)), and 2,048 elements (bottom row - (e) and (f)). Computed fields are displayed using dots, while analytical fields are displayed using lines.}
    \label{Hartmann_results}
\end{figure}

We have solved the Hartmann channel flow problem using a variety of Hartmann numbers, a polynomial degree of $k = 2$, and meshes with 128, 512, and 2,048 elements. The computed velocity and magnetic fields along the line $x = 0.25$ are displayed in Figures \ref{Hartmann_results}. Note that the computed velocity and magnetic fields are highly accurate for each mesh and Hartmann number, though slight overshoots and oscillations are seen in $B_x$ near the two plates for the coarsest mesh and highest Hartmann number considered due to underresolution of the sharp gradient in $B_x$.

\begin{figure}[t!]
    \centering
    \includegraphics[width=0.45\textwidth]{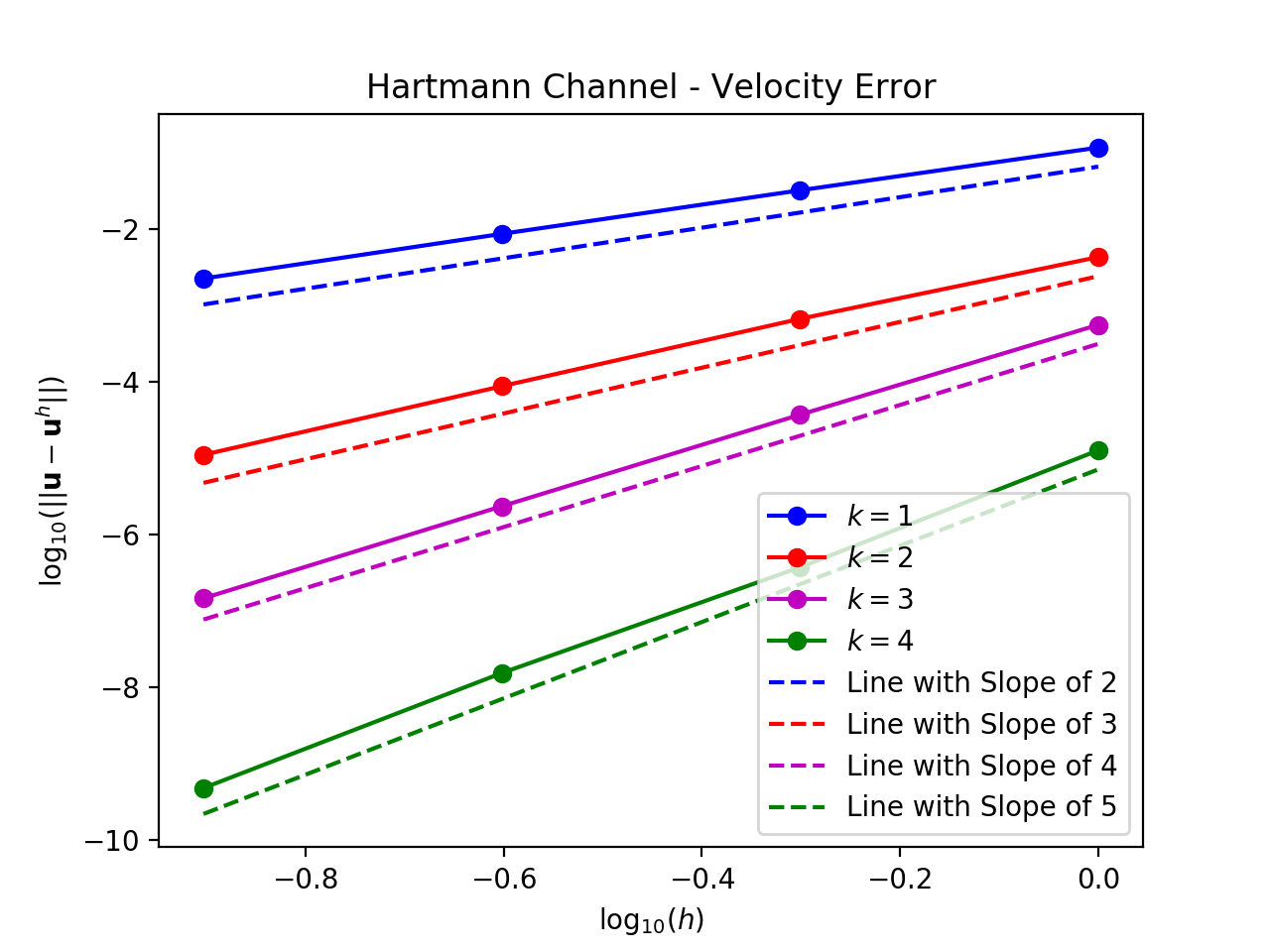}
    \includegraphics[width=0.45\textwidth]{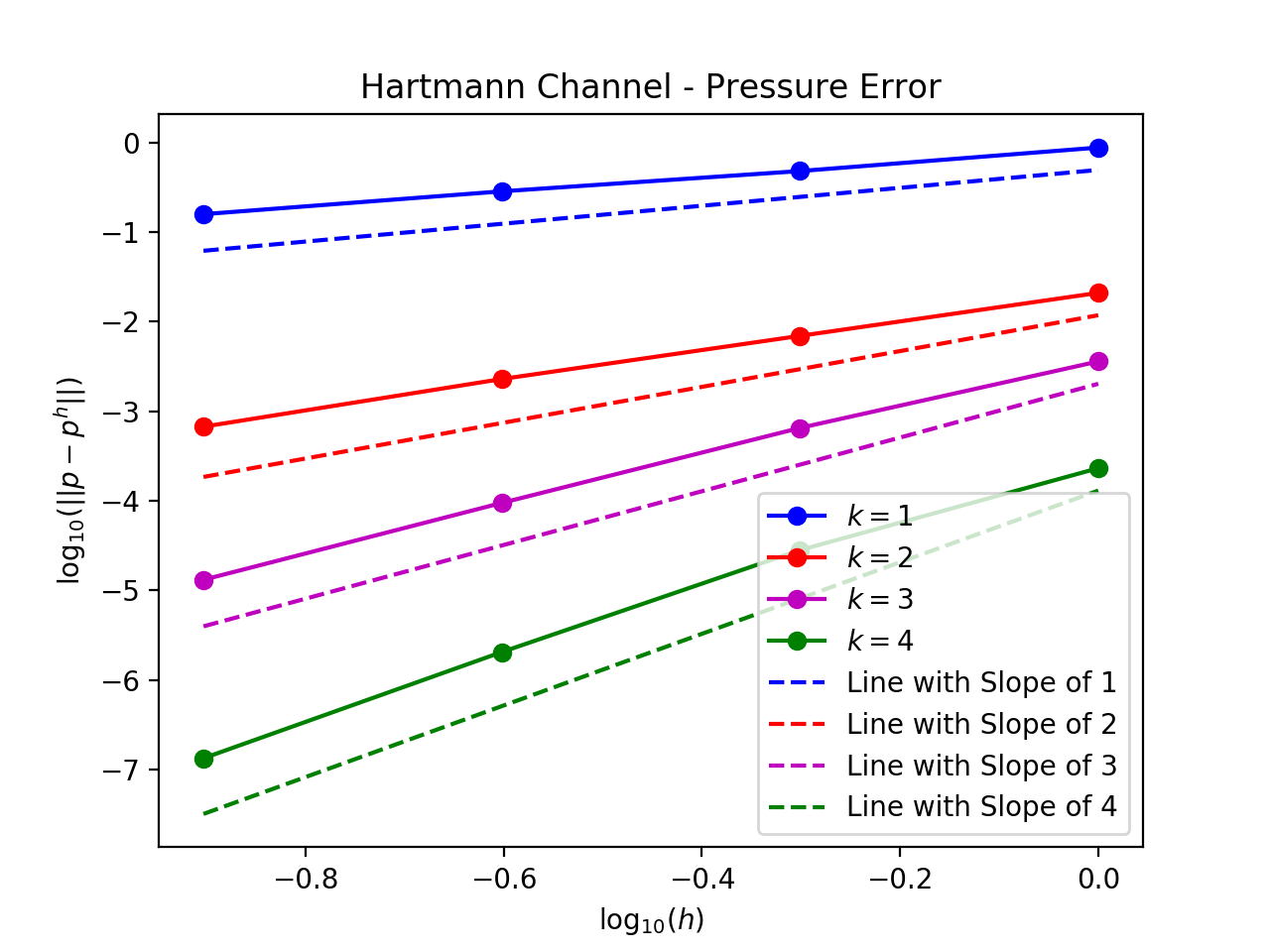}
    \includegraphics[width=0.45\textwidth]{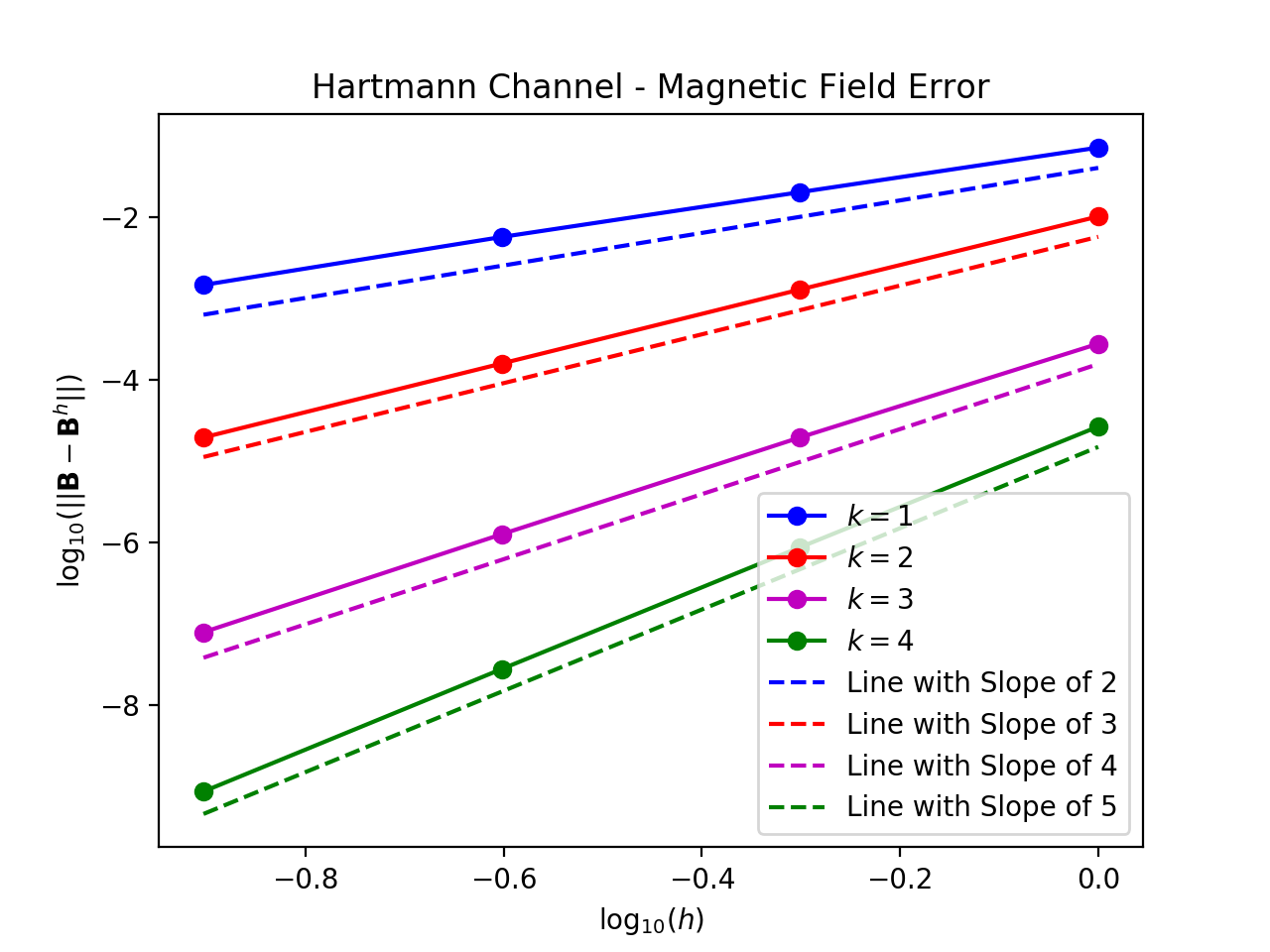}
    \includegraphics[width=0.45\textwidth]{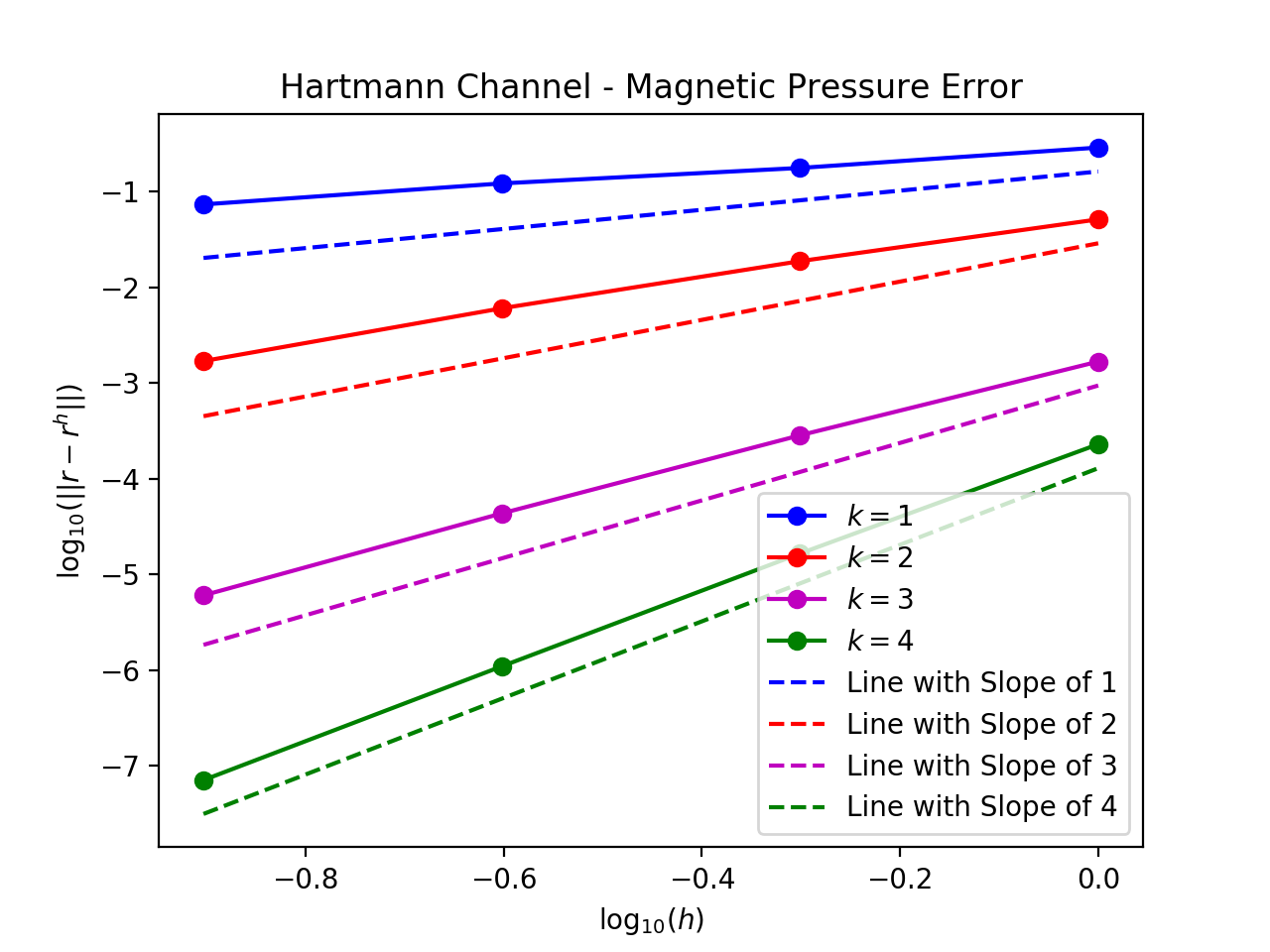}
    \caption{$L^2$-errors in the computed fields for the Hartmann channel flow problem for $Ha = 5$.}
    \label{hc_conv}
\end{figure}

We have also solved the Hartmann channel flow problem using a Hartmann number of $Ha = 5$, a series of meshes, and polynomial degrees $k = 1, 2, 3, 4$. Plots of the $L^2$-error in the computed velocity, pressure, magnetic, and magnetic pressure fields are shown in Fig. \ref{hc_conv}. Optimal rates are attained for each polynomial degree just as was the case for the manufactured solution problem.

\subsection{Alfv\'en Wave Propagation}

We next consider a generalization of Stokes' first problem to the setting of incompressible MHD. In this problem, an initially stationary fluid above a stationary flat plate located at $y = 0$ is suddenly brought into motion by impulsively moving the flat plate to the right at speed $U$ at time $t = 0$. The flow is further subjected to an applied magnetic field $B_y = B_0$ in the vertical direction. The fluid away from the plate is gradually brought into motion, and as this velocity profile is developed, a self-induced magnetic field $B_x$ in the horizontal direction is developed and an Alfv\'en wave is propagated through the fluid with velocity $A_0=B_0/\sqrt{\mu_0\rho}$. If the Reynolds and magnetic Reynolds numbers are equal (that is, if $\eta/\mu_0 = \nu$), the resulting velocity field $\textbf{u} = (u_x,u_y)$ and magnetic field $\textbf{B} = (B_x,B_y)$ take the form (see, e.g., \cite{Shadid1, Shadid2})
\begin{align}
    u_x =& \frac{U}{4}\left(e^{\frac{-A_0y}{d}}\left(1-\text{erf}\left(\frac{y-A_0t}{2\sqrt{dt}}\right)\right)-\text{erf}\left(\frac{y-A_0t}{2\sqrt{dt}}\right)\right) \nonumber \\
    &+\frac{U}{4}\left(e^{\frac{A_0y}{d}}\left(1-\text{erf}\left(\frac{y+A_0t}{2\sqrt{dt}}\right)\right)-\text{erf}\left(\frac{y+A_0t}{2\sqrt{dt}}\right)+2\right)
    \label{alfven_ux} \\
    u_y =& 0 \\
    B_x =& -\frac{1}{4}e^{\frac{-A_0y}{d}}\left(-1+e^{\frac{A_0y}{d}}\right)U\sqrt{\mu_0\rho}\left(\text{erfc}\left(\frac{y-A_0t}{2\sqrt{dt}}\right)+e^{\frac{A_0y}{d}}\text{erfc}\left(\frac{y+A_0t}{2\sqrt{dt}}\right)\right) \label{alfven_Bx} \\
    B_y =& B_0
\end{align}
where $d = \eta/\mu_0 = \nu$. Since this problem, which we refer to as the Alfv\'{e}n wave propagation problem, has a solution that depends on both space and time, it is a good candidate to assess the temporal accuracy of our method. In our numerical experiments, we set $U = 1$, $B_0 = 10$, and $\nu = \eta = 1$, and we consider the flow domain $\Omega=(0,1)\times(0,2.5)$. Dirichlet boundary conditions are again applied along all of $\partial \Omega$ for both the velocity and magnetic fields using the analytical solution. We use a uniform time step size $\Delta t$ in all of our numerical experiments, and our simulations start from a time of $t = 0.01$ rather than $t = 0$ as there is a singularity at $t = 0$ that limits the convergence rates of numerical methods.

\begin{figure}[t!]
    \centering
    \includegraphics[width=0.45\textwidth]{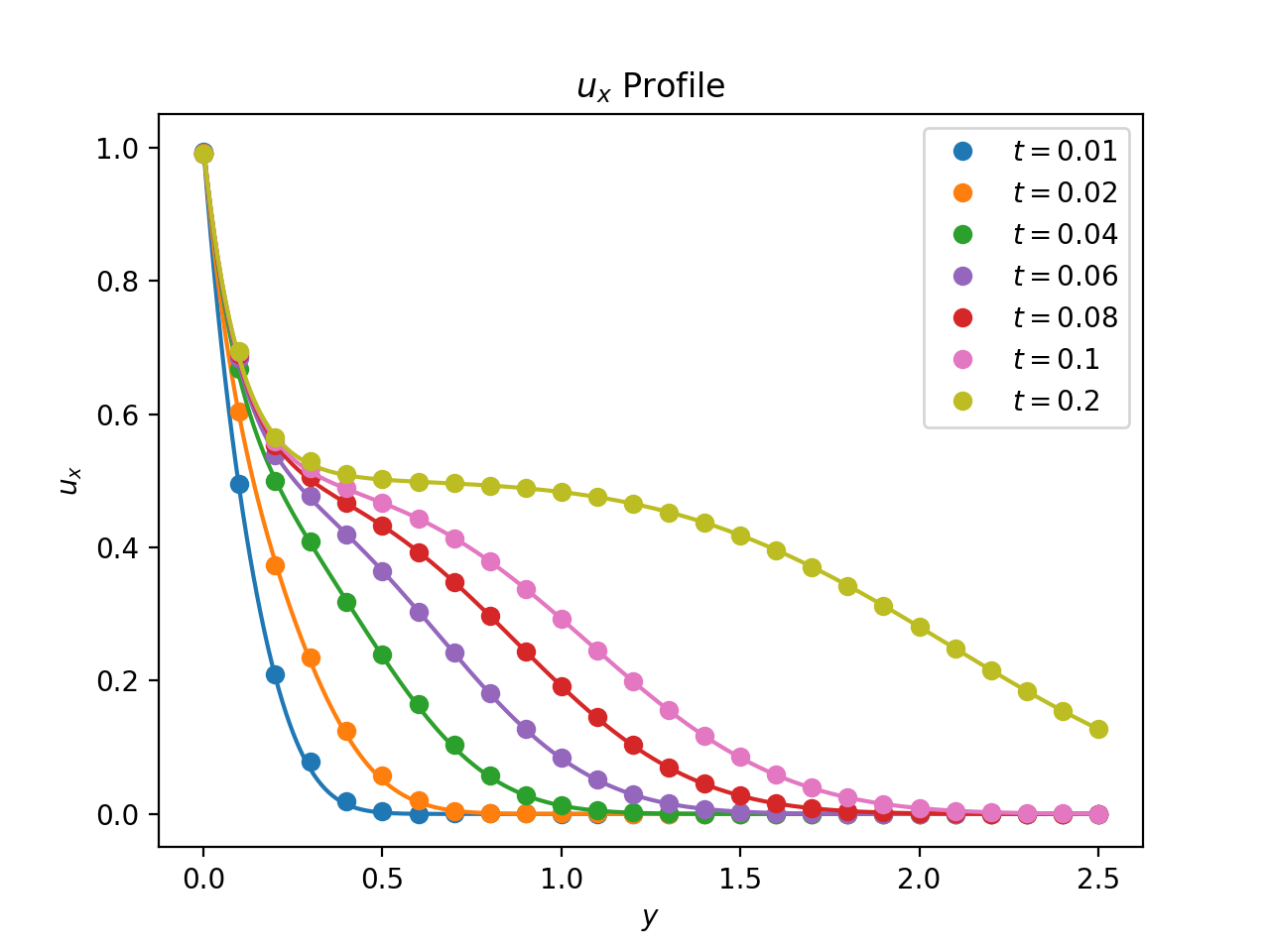}
    \includegraphics[width=0.45\textwidth]{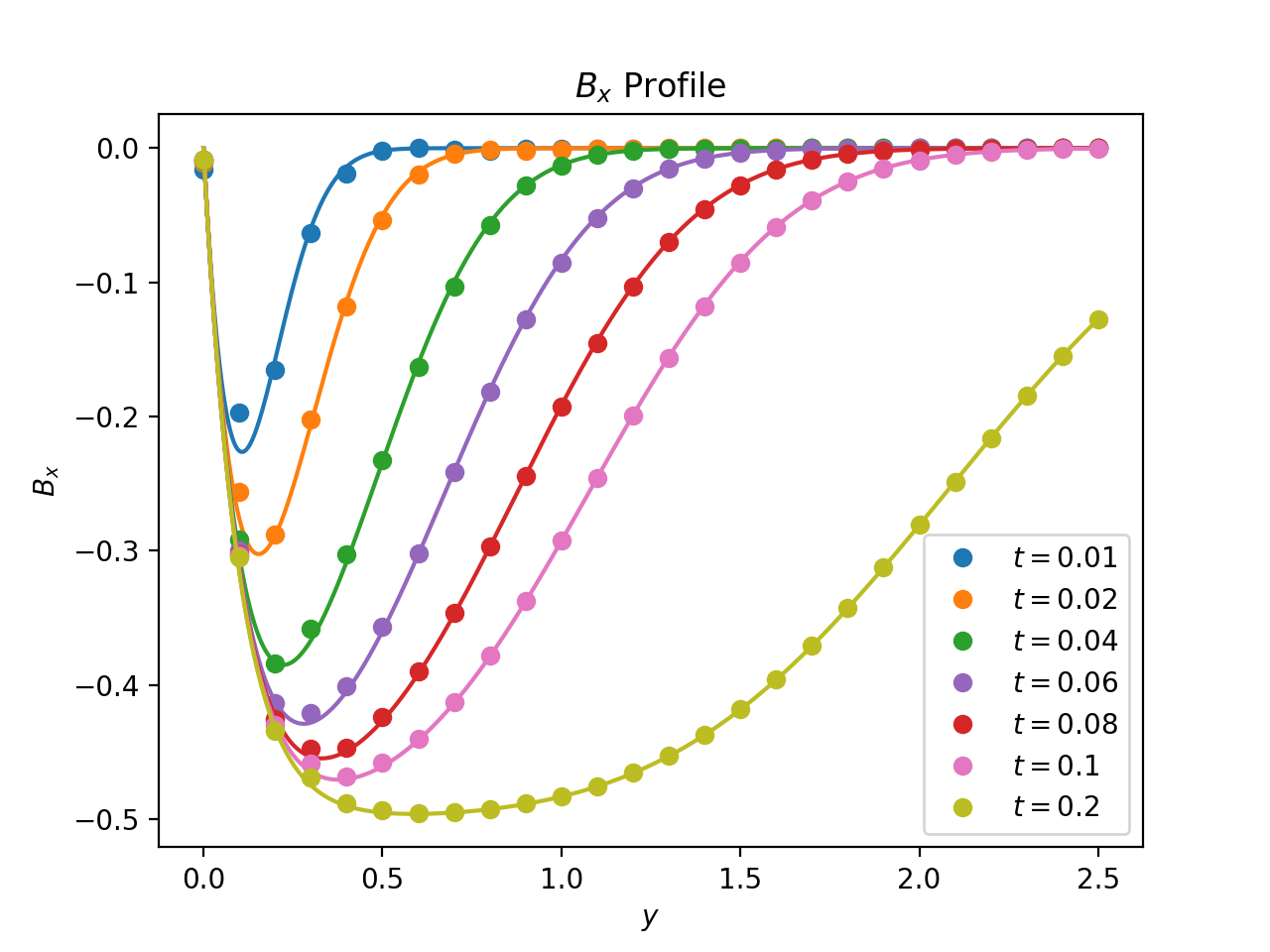} \\
    \hspace{5pt} (a) \hspace{195pt} (b) \\
    \includegraphics[width=0.45\textwidth]{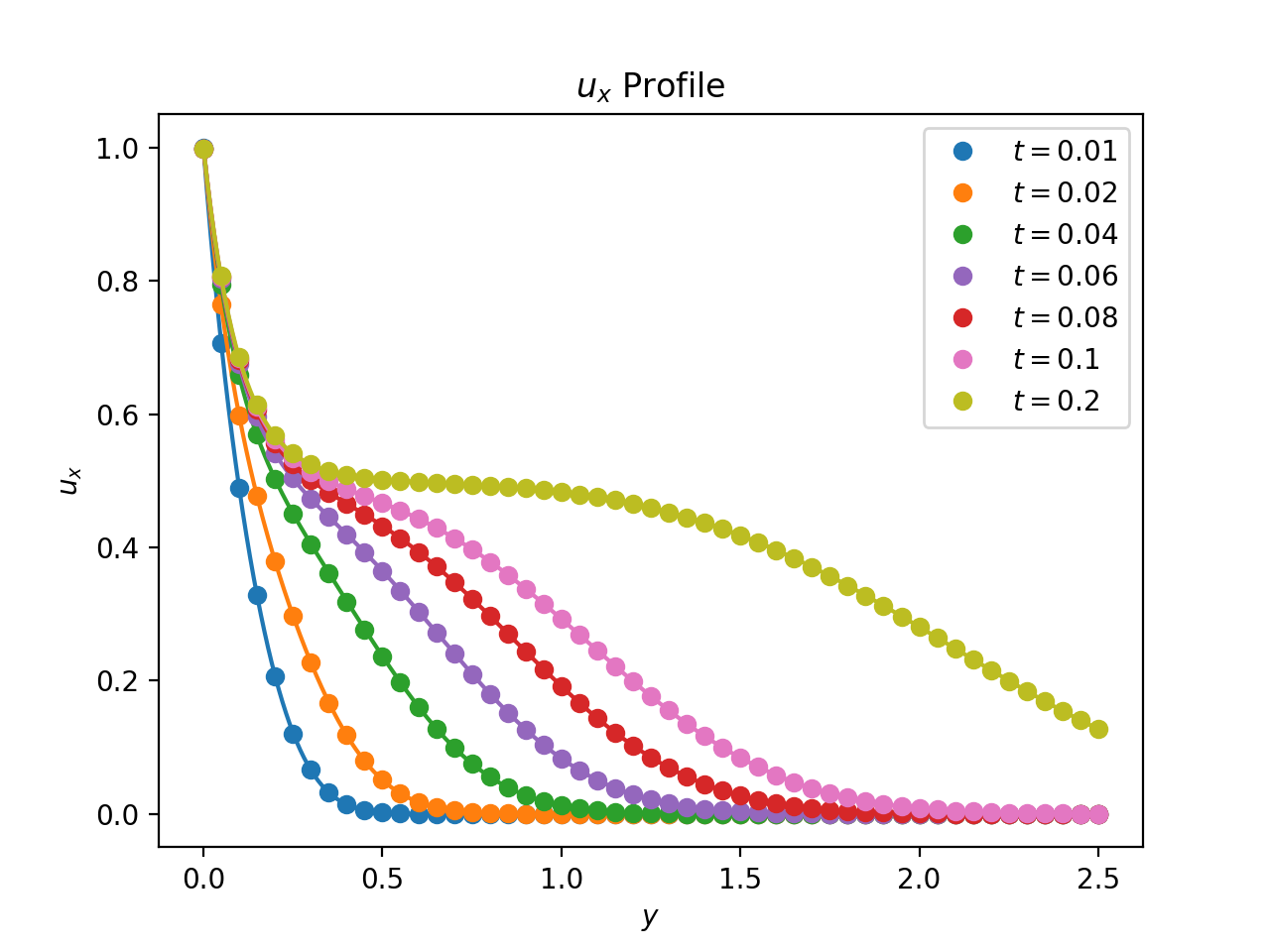}
    \includegraphics[width=0.45\textwidth]{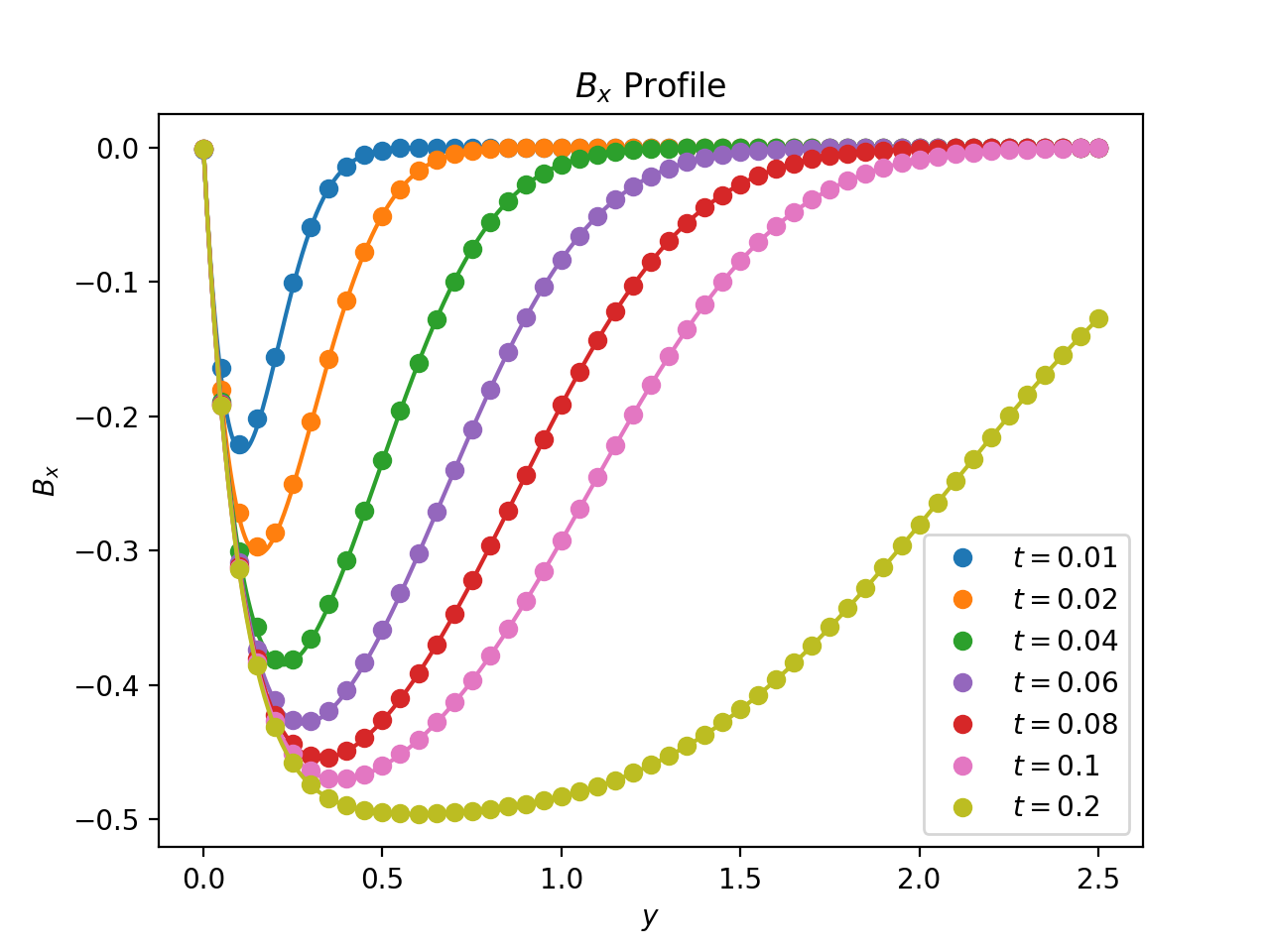} \\
    \hspace{5pt} (c) \hspace{195pt} (d) \\
    \includegraphics[width=0.45\textwidth]{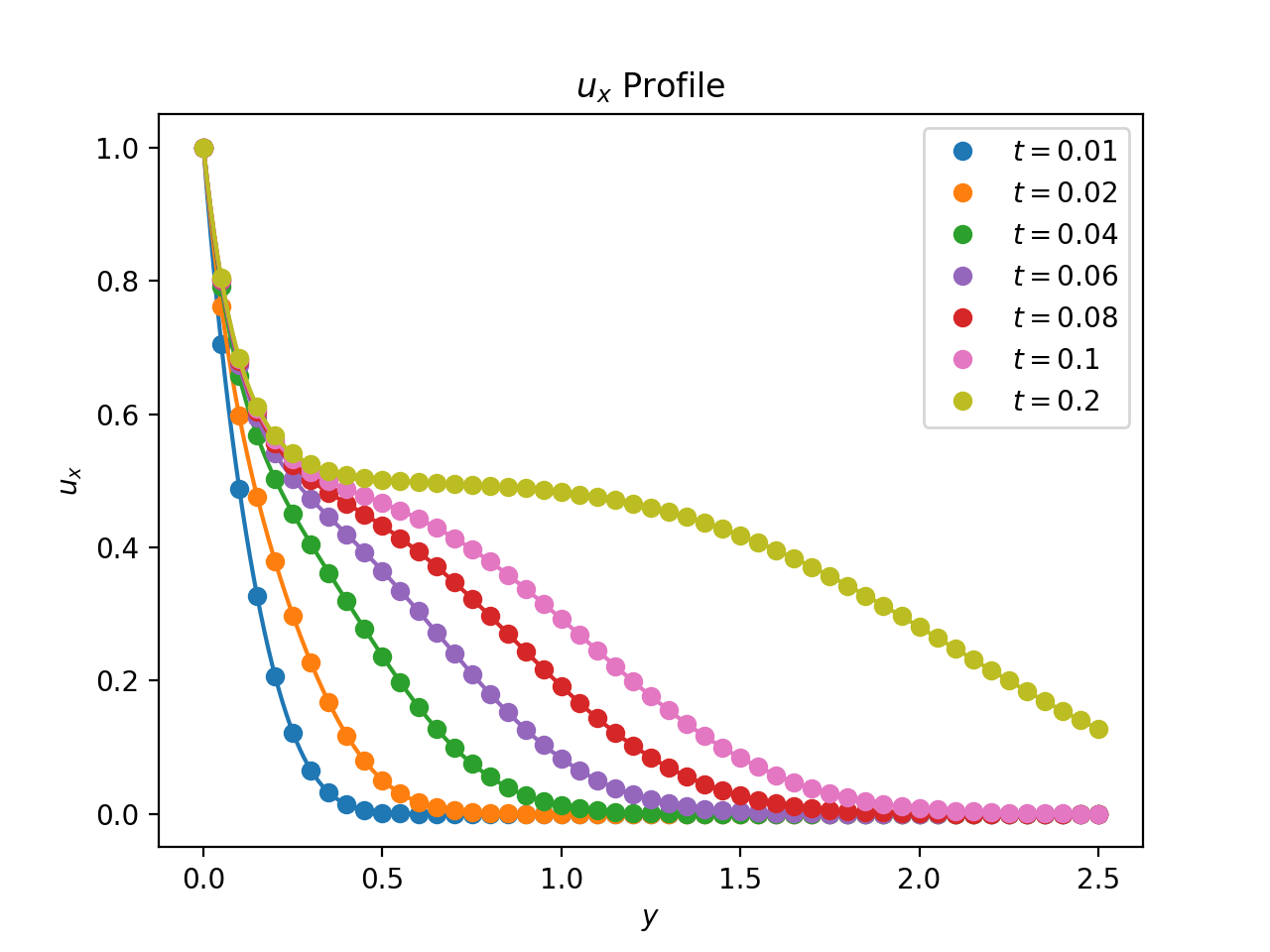}
    \includegraphics[width=0.45\textwidth]{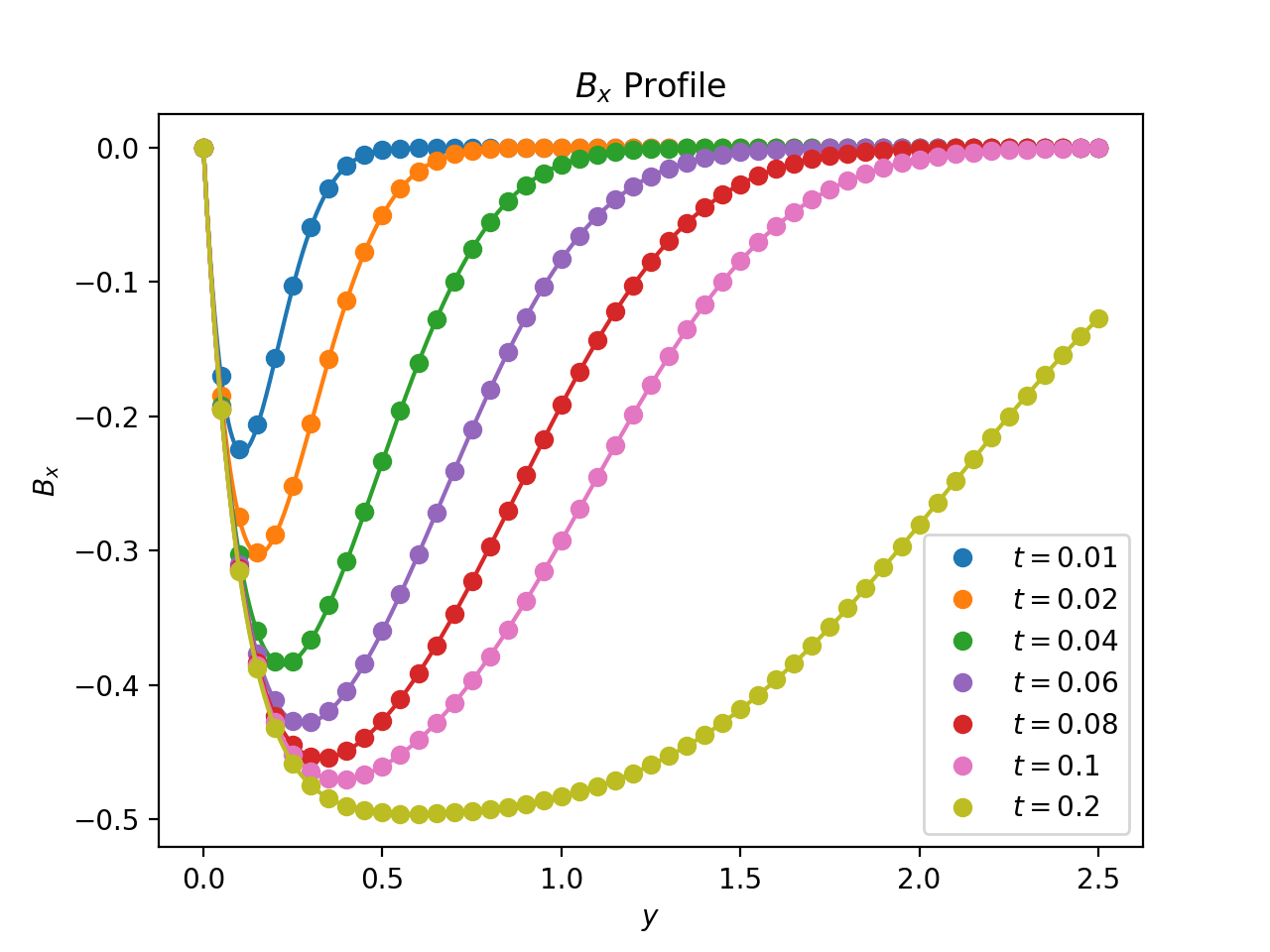} \\
    \hspace{5pt} (e) \hspace{195pt} (f) \\
    \caption{Computed $u_x$ and $B_x$ fields along the line $x = 0.5$ for the Alfv\'{e}n wave propagation problem for several time instances, a polynomial degree of $k = 4$, and a coarse mesh of 48 elements with $\Delta t = 1\cdot10^{-3}$ (top row - (a) and (b)), a medium mesh of 192 elements with $\Delta t = 5\cdot10^{-4}$ (middle row - (c) and (d)), and a fine mesh of 768 elements with $\Delta t = 2.5\cdot10^{-4}$ (bottom row - (e) and (f)). Computed fields are displayed using dots, while analytical fields are displayed using lines.}
    \label{alfven_results}
\end{figure}

We have solved the Alfv\'{e}n wave propagation problem using a polynomial degree of $k = 4$ and a coarse mesh of 48 elements with a time step size of $\Delta t = 1\cdot10^{-3}$, a medium mesh of 192 elements with a time step size of $\Delta t = 5\cdot10^{-4}$, and a fine mesh of 768 elements with a time step size of $\Delta t = 2.5\cdot10^{-4}$. The computed velocity and magnetic fields along the line $x = 0.5$ are displayed for several time instances in Fig. \ref{alfven_results}. The computed velocity and magnetic fields are highly accurate for each considered time instance and mesh, even the coarse mesh with only 48 elements.

\begin{figure}[t!]
    \centering
    \includegraphics[width=0.45\textwidth]{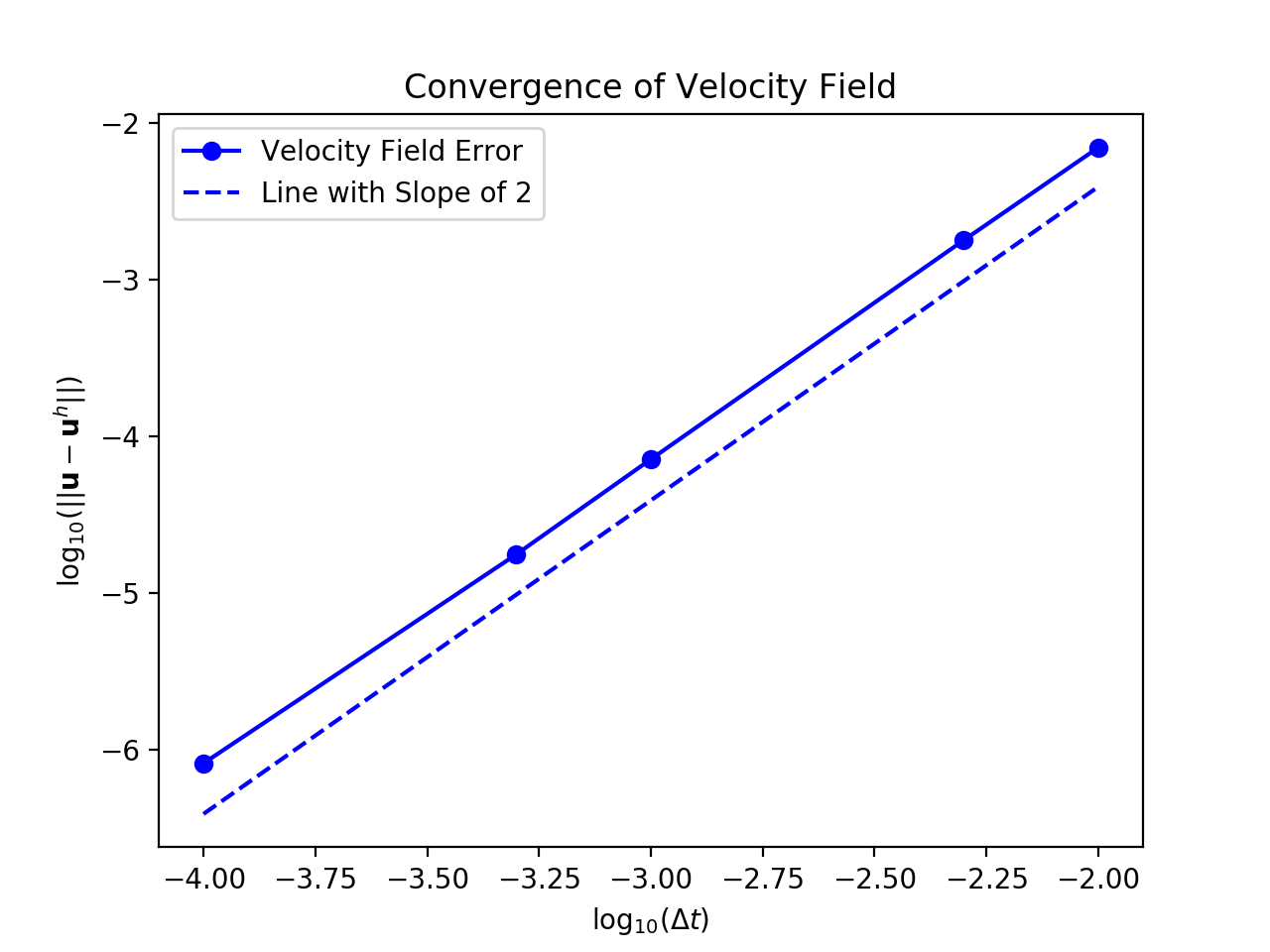}
    \includegraphics[width=0.45\textwidth]{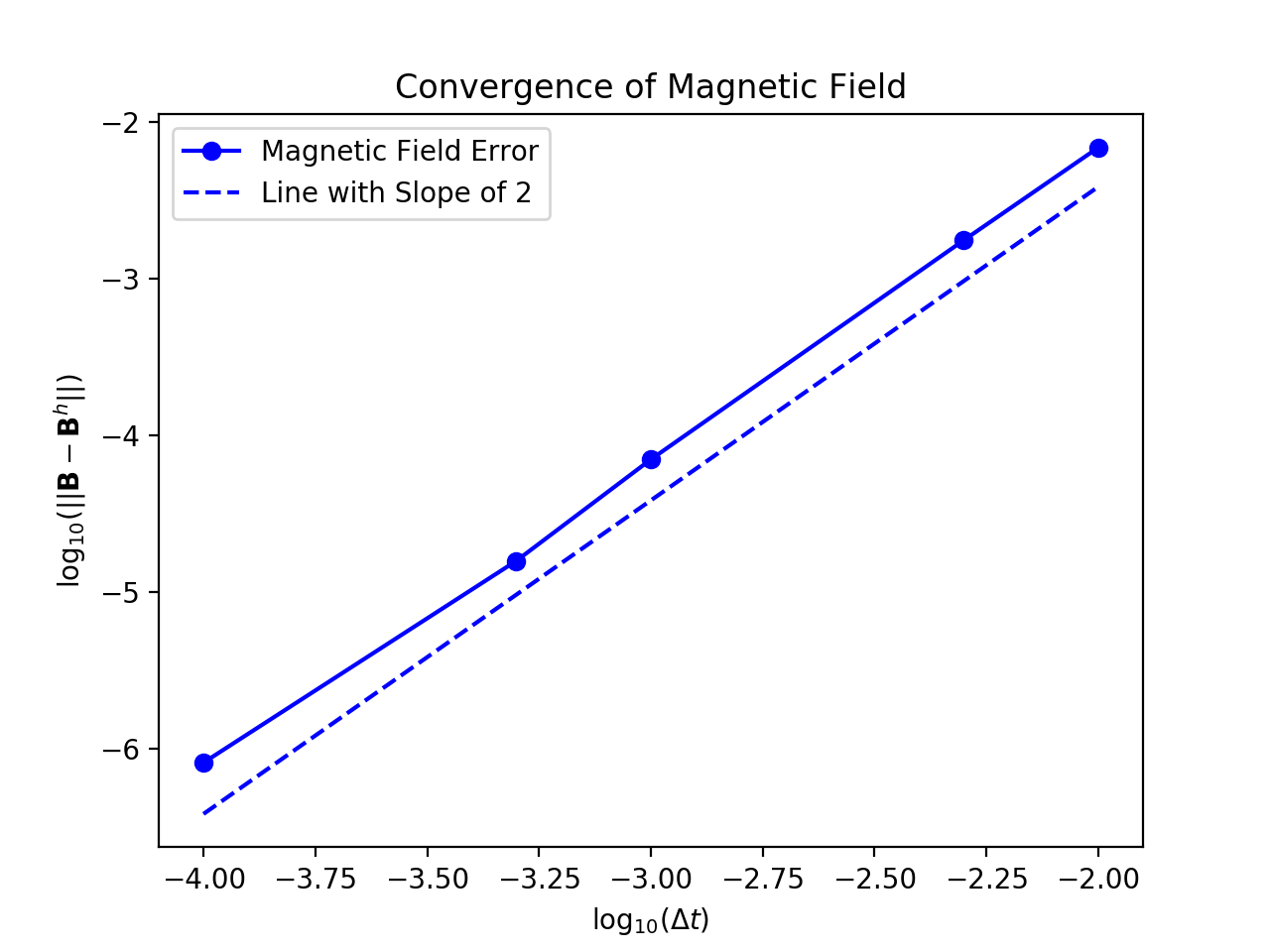}
    \caption{$L^2$-errors in the computed velocity and magnetic fields for the Alfv\'{e}n wave propagation problem for a polynomial degree of $k = 4$, a fine mesh of 768 elements, and a sequence of time step sizes.}
    \label{temp_conv}
\end{figure}

To better assess the temporal accuracy of our chosen time discretization scheme, we have solved the Alfv\'{e}n wave propagation problem using a polynomial degree of $k = 4$, the fine mesh of 768 elements, and a sequence of time step sizes. As a high polynomial degree and fine mesh is employed for this numerical experiment, we expect that temporal discretization error dominates spatial discretization error. Plots of the $L^2$-error in the computed velocity and magnetic pressure fields are shown in Fig. \ref{temp_conv}. Note that a temporal convergence rate of two is attained, indicating that our time discretization scheme is indeed second-order-in-time.

\begin{figure}[t!]
    \centering
    \includegraphics[width=0.45\textwidth]{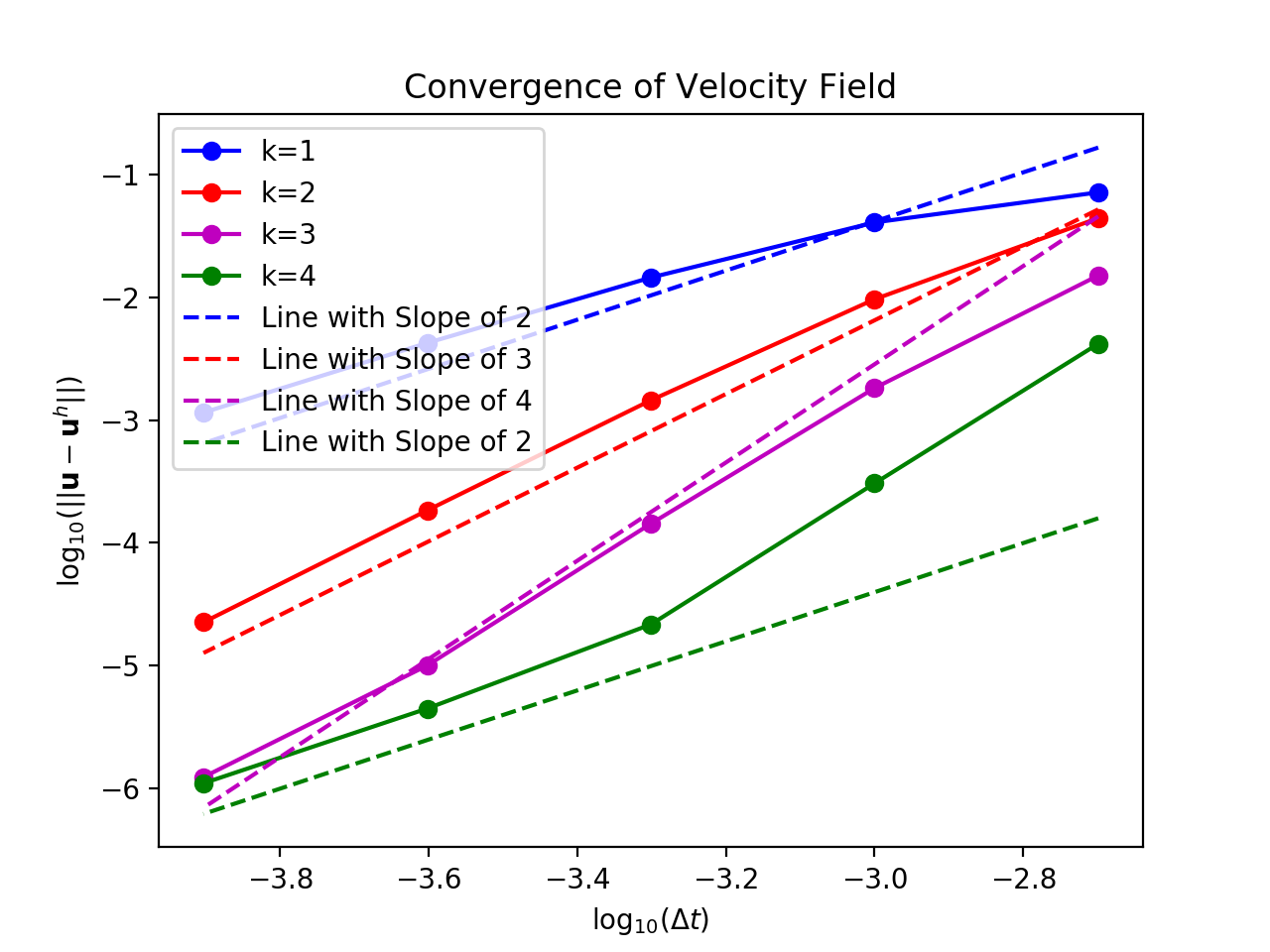}
    \includegraphics[width=0.45\textwidth]{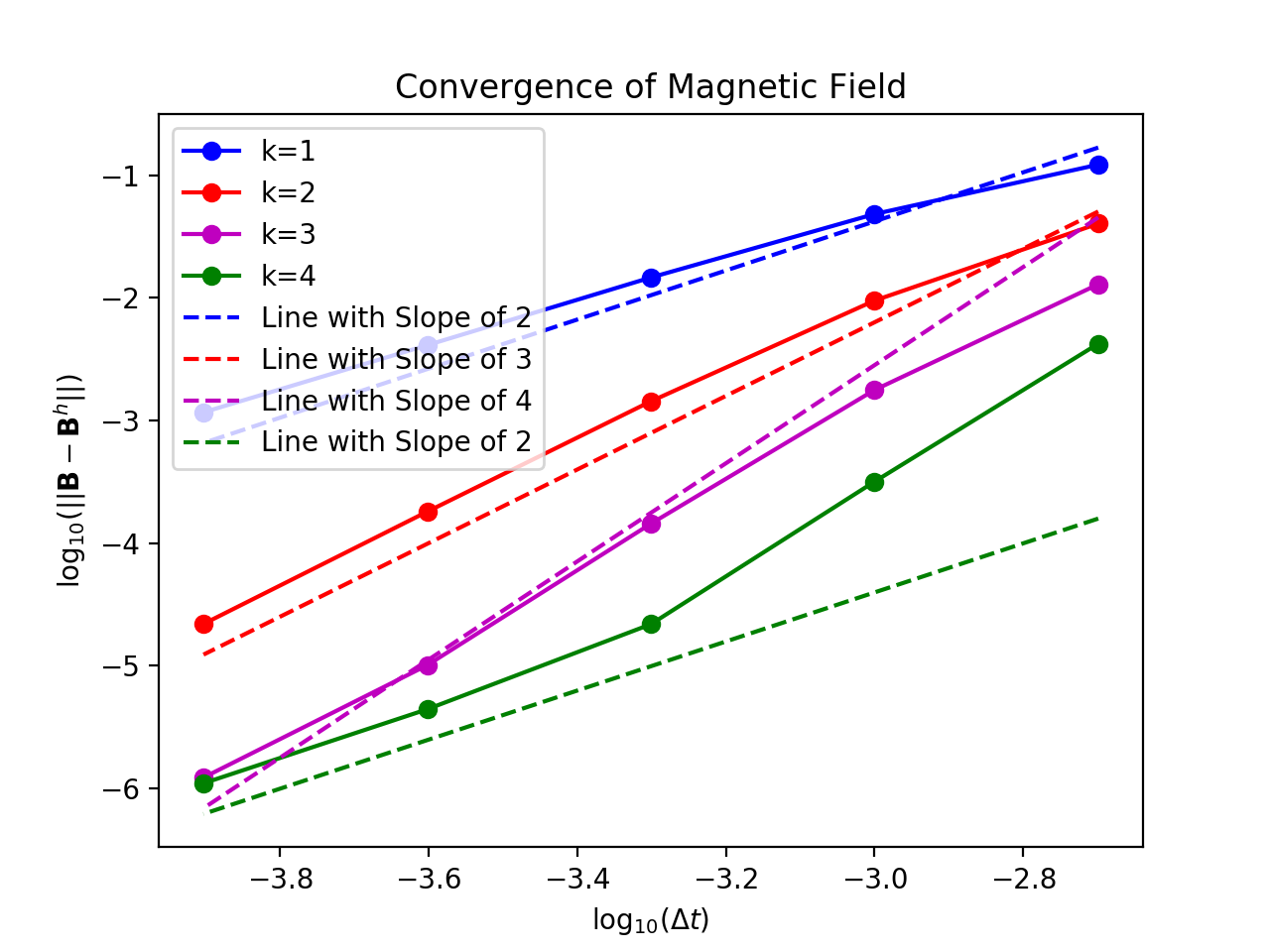}
    \caption{$L^2$-errors in the computed velocity and magnetic fields for the Alfv\'{e}n wave propagation problem for polynomial degrees of $k = 1, 2, 3, 4$ and a sequence of meshes and time step sizes with fixed CFL number $U \Delta t/h $.}
    \label{temp_conv2}
\end{figure}

To assess the combined spatial and temporal accuracy of our method, we have solved the Alfv\'{e}n wave propagation problem using polynomial degrees $k = 1, 2, 3, 4$ and a sequence of meshes and time step sizes with a fixed CFL number $U \Delta t/h$. Plots of the $L^2$-error in the computed velocity and magnetic pressure fields are shown in Fig. \ref{temp_conv2}. A convergence rate of two is observed for polynomial degree $k = 1$ as expected as the spatial and temporal discretization errors are balanced. Convergence rates of three and four are observed for polynomial degrees $k = 2$ and $k = 3$, indicating that spatial discretization error dominates temporal discretization error for the considered mesh sizes. An asymptotic convergence rate of two is observed for polynomial degree $k = 4$, indicating that temporal discretization error dominates spatial discretization error for sufficiently fine meshes as expected.

\subsection{Kelvin-Helmholtz Instability}

\begin{figure}[t!]
    \centering
    \includegraphics[width=0.45\textwidth]{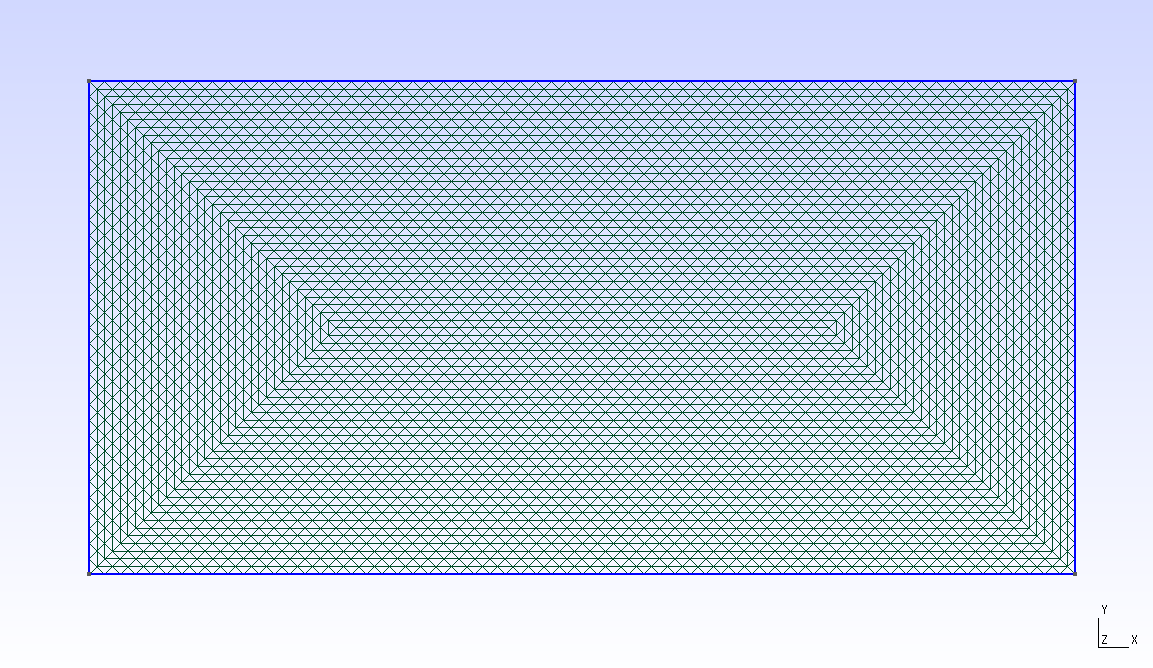}
    \includegraphics[width=0.45\textwidth]{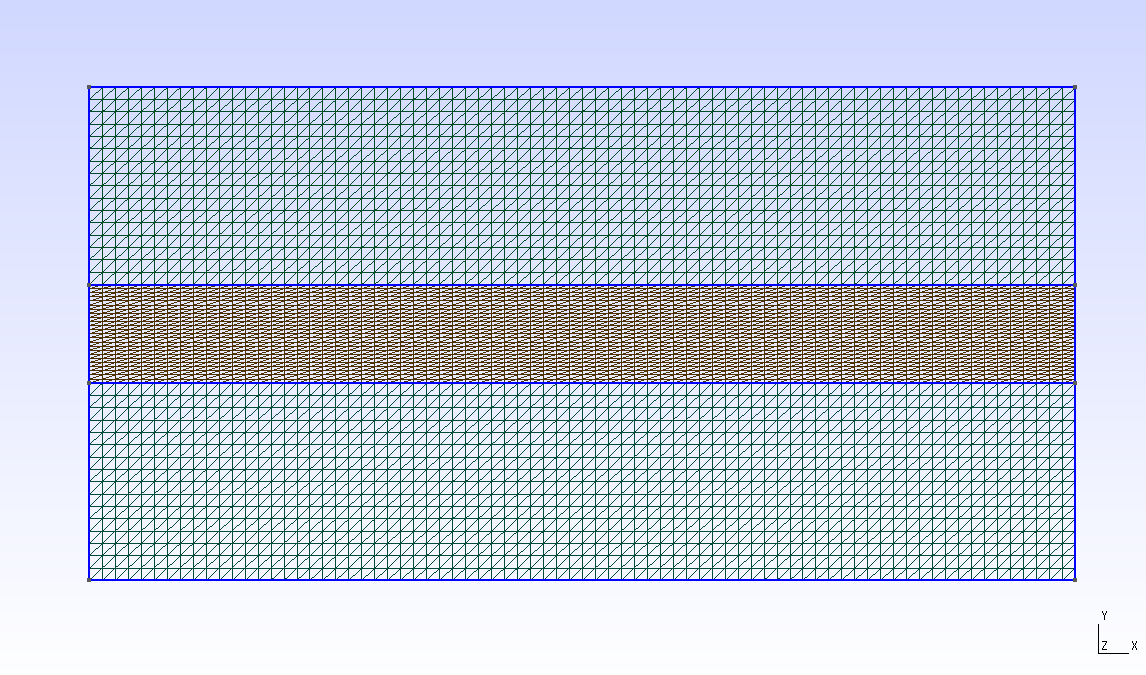}
    \caption{Meshes used for the Kelvin-Helmholtz instability problem for $Re = Re_m = 10^3$ (left) and $Re = Re_m = 10^4$ (right).}
    \label{KHmeshes}
\end{figure}

We finally consider a generalization of the classical Kelvin-Helmholtz instability problem to the setting of incompressible MHD. Like the classical Kelvin-Helmholtz instability problem, the initial velocity field consists of a shear flow that is sufficiently strong to be Kelvin-Helmholtz unstable. However, an initial magnetic field is also considered in the direction of flow. If the initial magnetic field is sufficiently strong, the magnetic stress stabilizes the shear layer. Otherwise, the shear layer breaks down into vortices. Like \cite{Shadid3,fu}, we use the Kelvin-Helmholtz instability problem to study the energy stability of our method. In our numerical experiments, we consider the domain $\Omega = (0,2)\times(0,1)$ and initial conditions $\textbf{u}_0 = (u_{x,0},u_{y,0})$ and $\textbf{B}_0 = (B_{x,0},B_{y,0})$ with
\begin{align}
u_{x,0} &=
\begin{cases}
U & y\geq 0.5 \\
-U & y < 0.5
\end{cases} \\
u_{y,0} &= 0 \\
B_{x,0} &= B_0\tanh(y/0.1) \\
B_{y,0} &= 0
\end{align}
at $t = 0$ where $U \in \mathbb{R}$ and $B_0 \in \mathbb{R}$. We consider periodic boundary conditions in the $x$-direction, zero normal boundary conditions $\textbf{u} \cdot \textbf{n} = 0$ and $\textbf{B} \cdot \textbf{n} = 0$ along the top and bottom of the domain, and zero tangential traction boundary conditions $\textbf{t}_v \times \textbf{n} = \textbf{0}$ and $\textbf{t}_m \times \textbf{n} = \textbf{0}$ along the top and bottom of the domain where
\begin{align}
\textbf{t}_v &= \mathbf{n}\cdot\left(-\frac{1}{\rho}p\mathds{1}+2\nu\nabla^s\mathbf{u}+\frac{1}{\rho \mu_0}\left(\mathbf{B}\otimes\mathbf{B}-\frac{1}{2}|\mathbf{B}|^2\mathds{1}\right)\right) - \min(\mathbf{u}\cdot\mathbf{n},0)\mathbf{u}\\
\textbf{t}_m &= \mathbf{n}\cdot\left(-r\mathds{1}+2\frac{\eta}{\mu_0} \nabla^a\mathbf{B}\right) - \min(\mathbf{u}\cdot\mathbf{n},0)\mathbf{B}.
\end{align}
Our method as presented in this paper is easily modified to be able to handle these boundary conditions, and Proposition \ref{prop:energy} extends to this case. The Alfv\'{e}n Mach number is defined as $M_A = U/U_A$ for this problem where $U_A = B_0/\sqrt{\rho\mu_0}$ is the Alfv\'{e}n wave speed. When $M_A < 1$, the magnetic field is strong enough to suppress the Kelvin-Helmholtz instability, and when $M_A > 1$, it is not strong enough to do so. In our numerical experiments, we set $U = 1$, and we consider $B_0 = \frac{1}{3}$ and $B_0 = 1.2$, corresponding to $M_A = 3$ and $M_A = \frac{5}{6}$ respectively. We also consider both $\nu = \eta = 10^{-3}$, corresponding to $Re = Re_m = 10^3$ for $\mathcal{U} = \mathcal{L} = 1$, and $\nu = \eta = 10^{-4}$, corresponding to $Re = Re_m = 10^4$. As the required spatial resolutions for these two Reynolds numbers differ, we employ two different meshes in our simulations, one for $Re = Re_m = 10^3$ and one for $Re = Re_m = 10^4$. These meshes are displayed in Fig. \ref{KHmeshes}. We employ a polynomial degree of $k = 2$ in each of our simulations.

\begin{figure}[t!]
    \centering
    \includegraphics[width=0.39\textwidth]{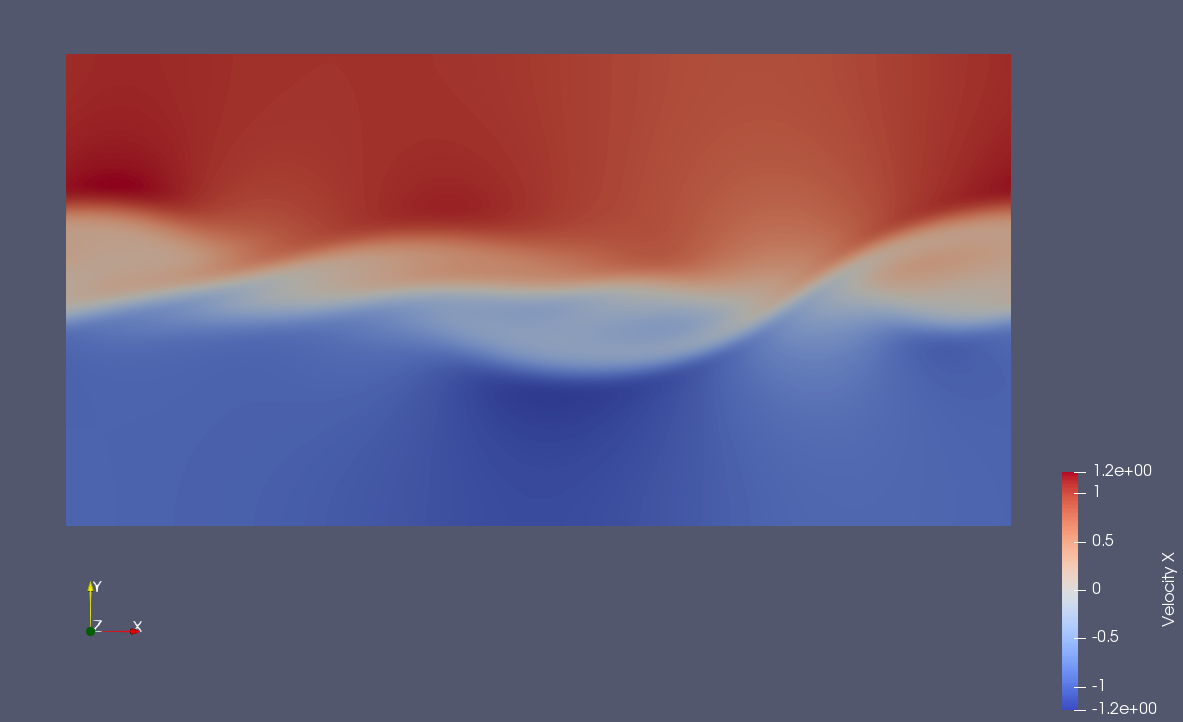}
    \includegraphics[width=0.39\textwidth]{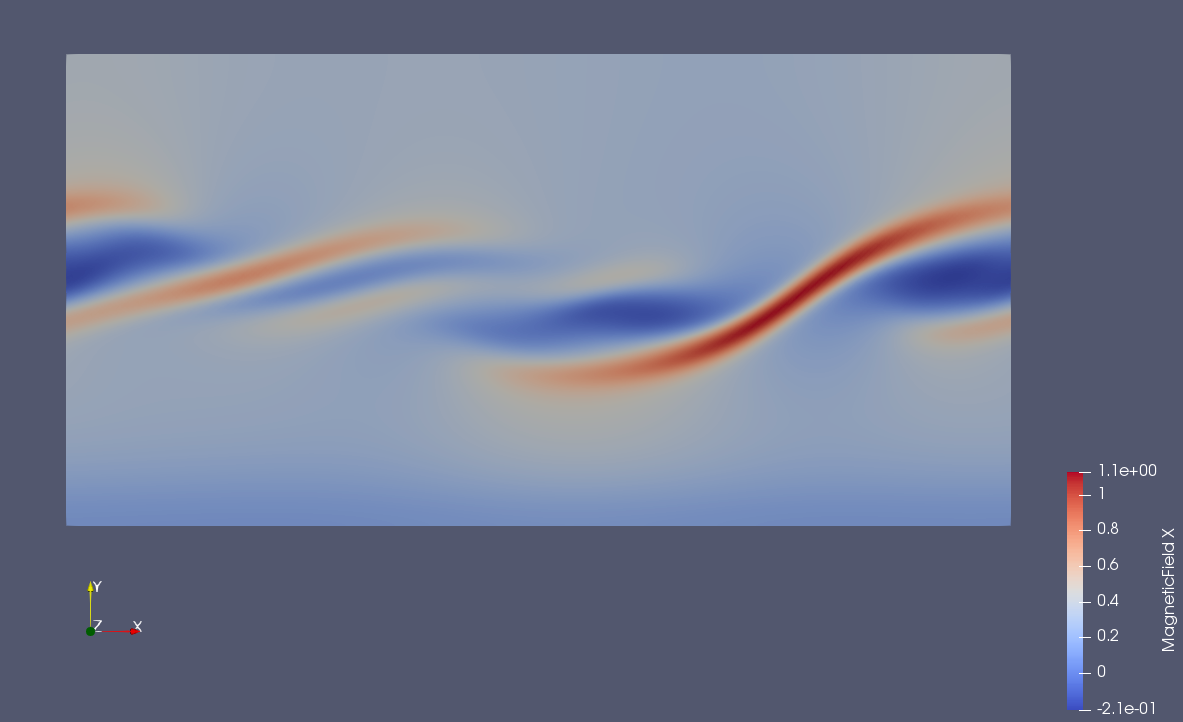}
    \caption{The $x$-component of the velocity (left) and magnetic (right) fields at $t = 2$ for the Kelvin-Helmholtz instability problem for $M_A = 3$ and $Re=Re_m=10^3$.}
    \label{KH1}
\end{figure}

\begin{figure}[t!]
    \centering
    \includegraphics[width=0.39\textwidth]{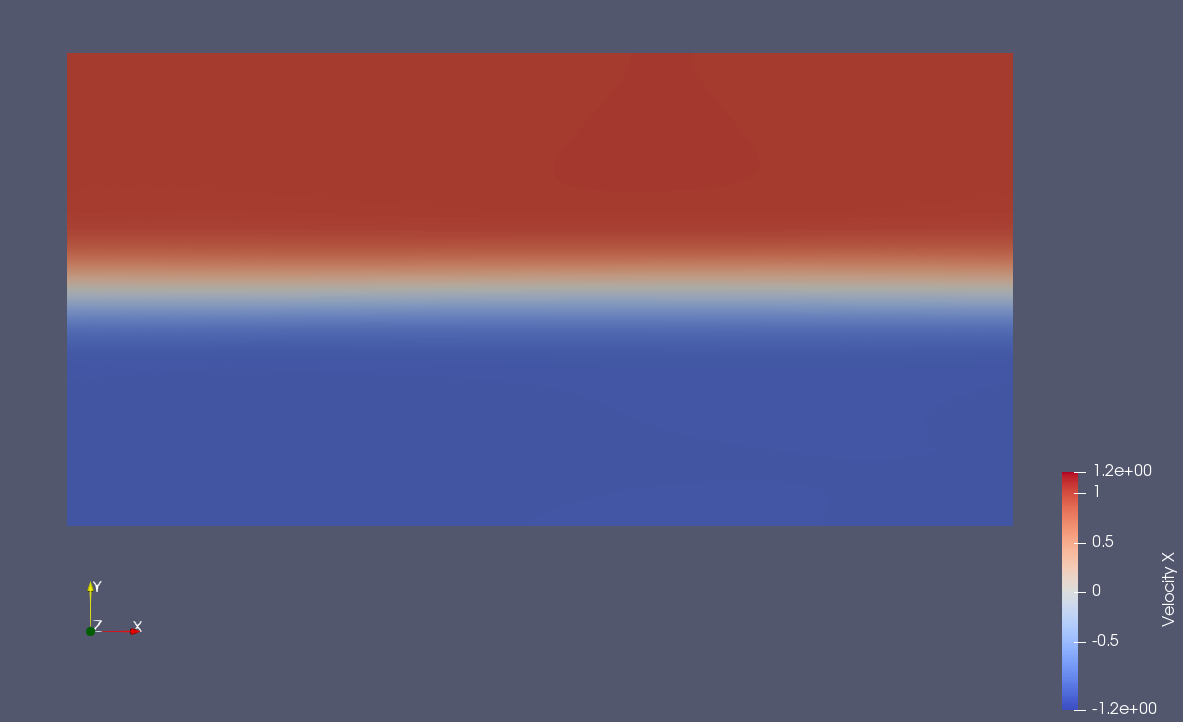}
    \includegraphics[width=0.39\textwidth]{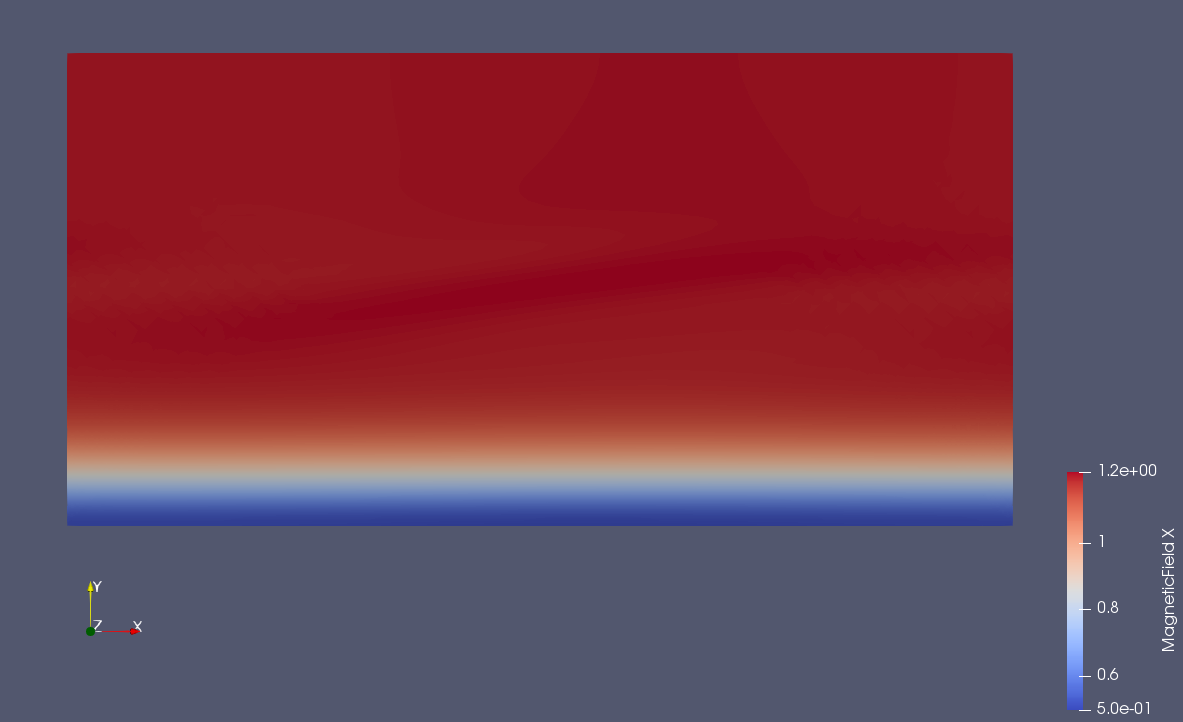}
    \caption{The $x$-component of the velocity (left) and magnetic (right) fields at $t = 2$ for the Kelvin-Helmholtz instability problem for $M_A = \frac{5}{6}$ and $Re=Re_m=10^3$.}
    \label{KH2}
\end{figure}

\begin{figure}[t!]
    \centering
    \includegraphics[width=0.39\textwidth]{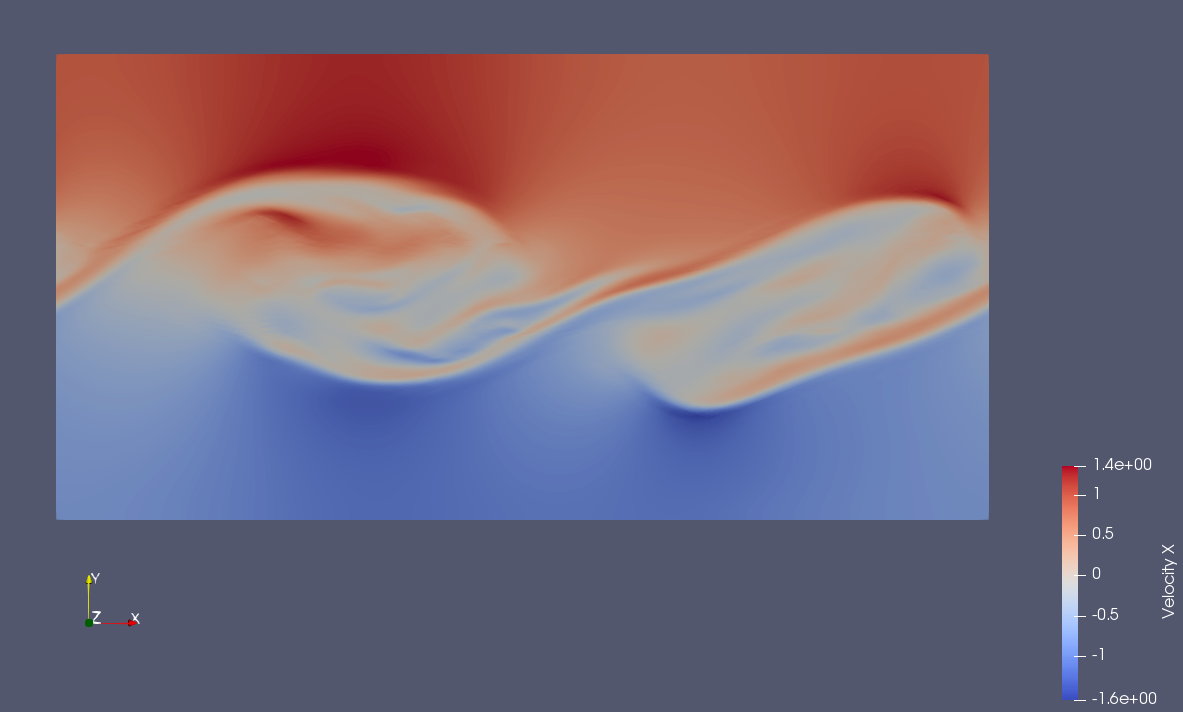}
    \includegraphics[width=0.39\textwidth]{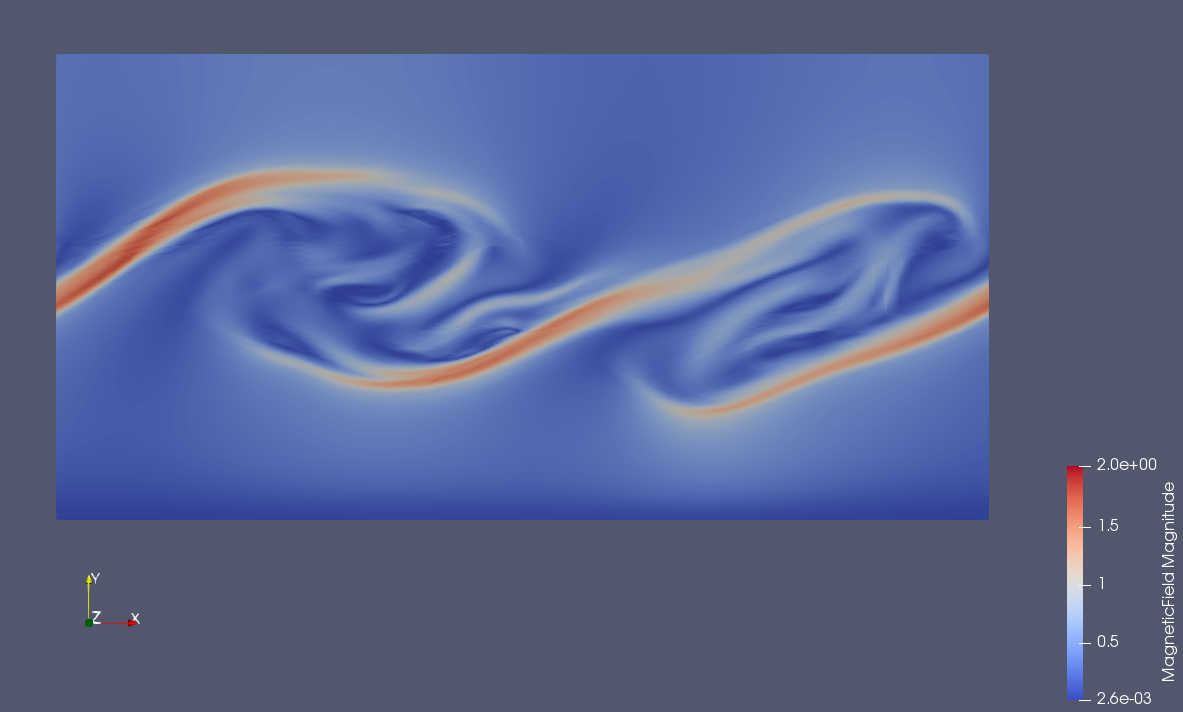}
    \caption{The $x$-component of the velocity (left) and magnetic (right) fields at $t = 2$ for the Kelvin-Helmholtz instability problem for $M_A = 3$ and $Re=Re_m=10^4$.}
    \label{KH3}
\end{figure}

\begin{figure}[t!]
    \centering
    \includegraphics[width=0.39\textwidth]{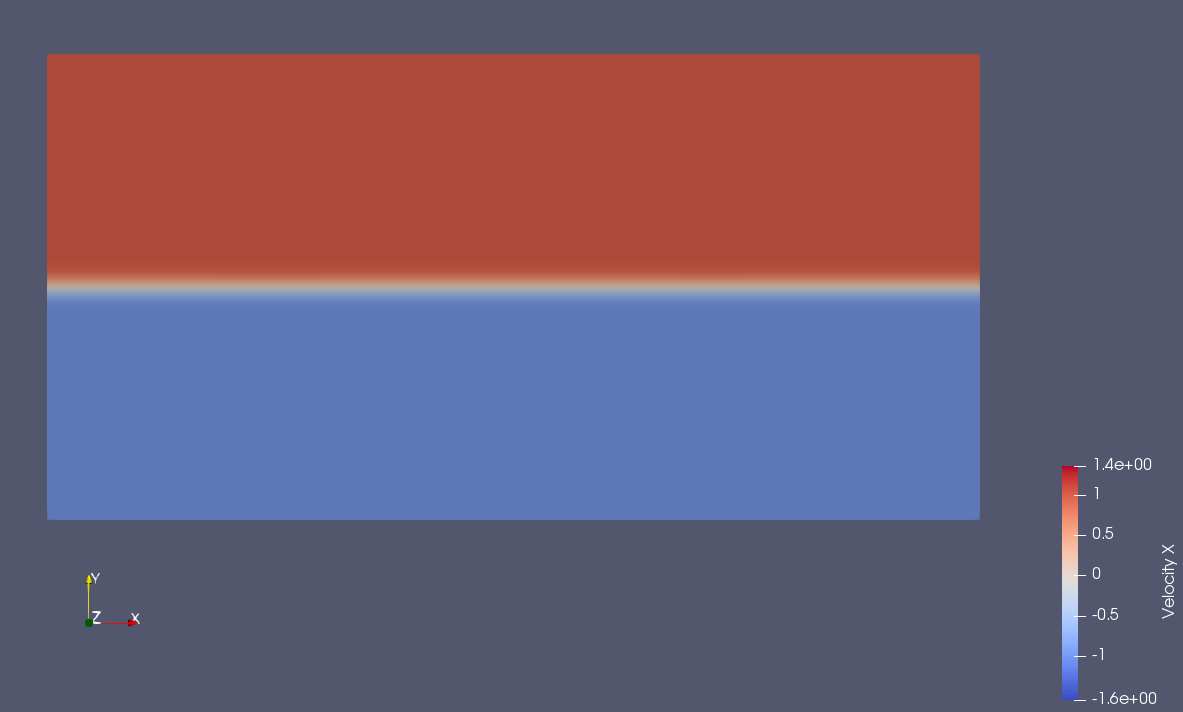}
    \includegraphics[width=0.39\textwidth]{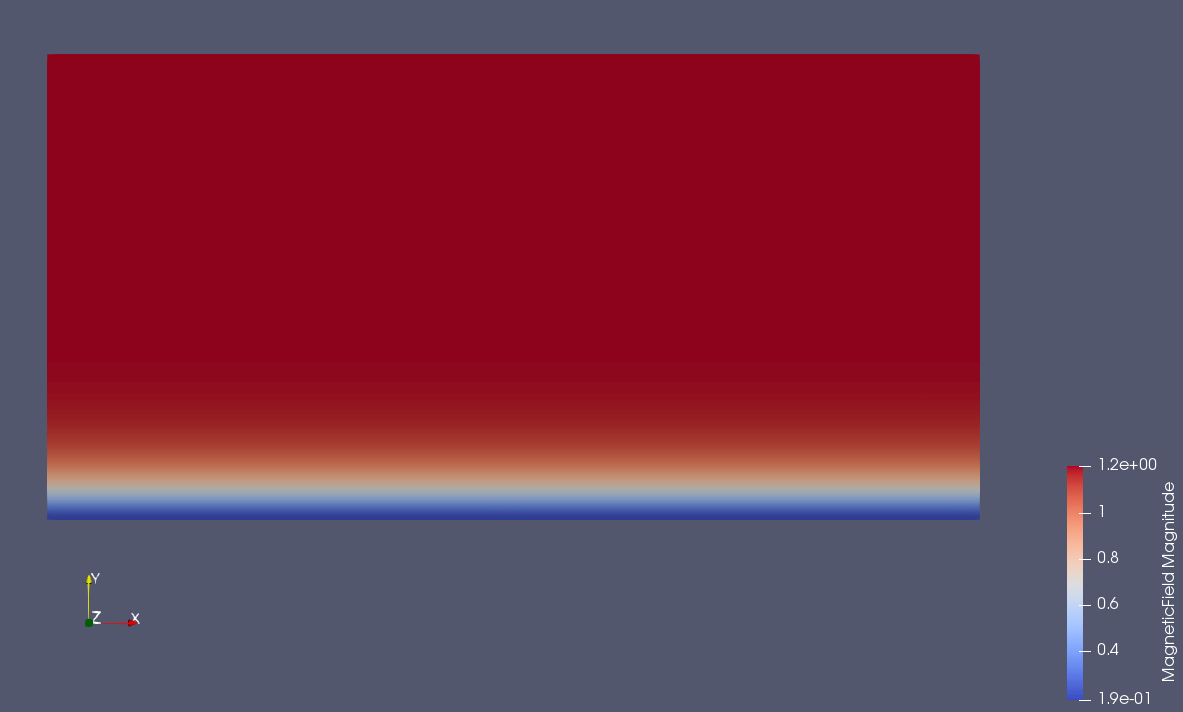}
    \caption{The $x$-component of the velocity (left) and magnetic (right) fields at $t = 2$ for the Kelvin-Helmholtz instability problem for $M_A = \frac{5}{6}$ and $Re=Re_m=10^4$.}
    \label{KH4}
\end{figure}

In Fig. \ref{KH1}-\ref{KH4}, the $x$-components of the velocity and magnetic fields as predicted by our simulations at time $t = 2$ are displayed. As expected, when $M_A < 1$, we see the shear layer remains stable both for $Re = Re_m = 10^3$ and $Re = Re_m = 10^4$. When $M_A > 1$, the shear layer curls and vortices are formed as expected. Finer structures have been generated for $Re = Re_m = 10^4$ as compared with $Re = Re_m = 10^3$.

\begin{figure}[t!]
    \centering
    \includegraphics[width=0.31\textwidth]{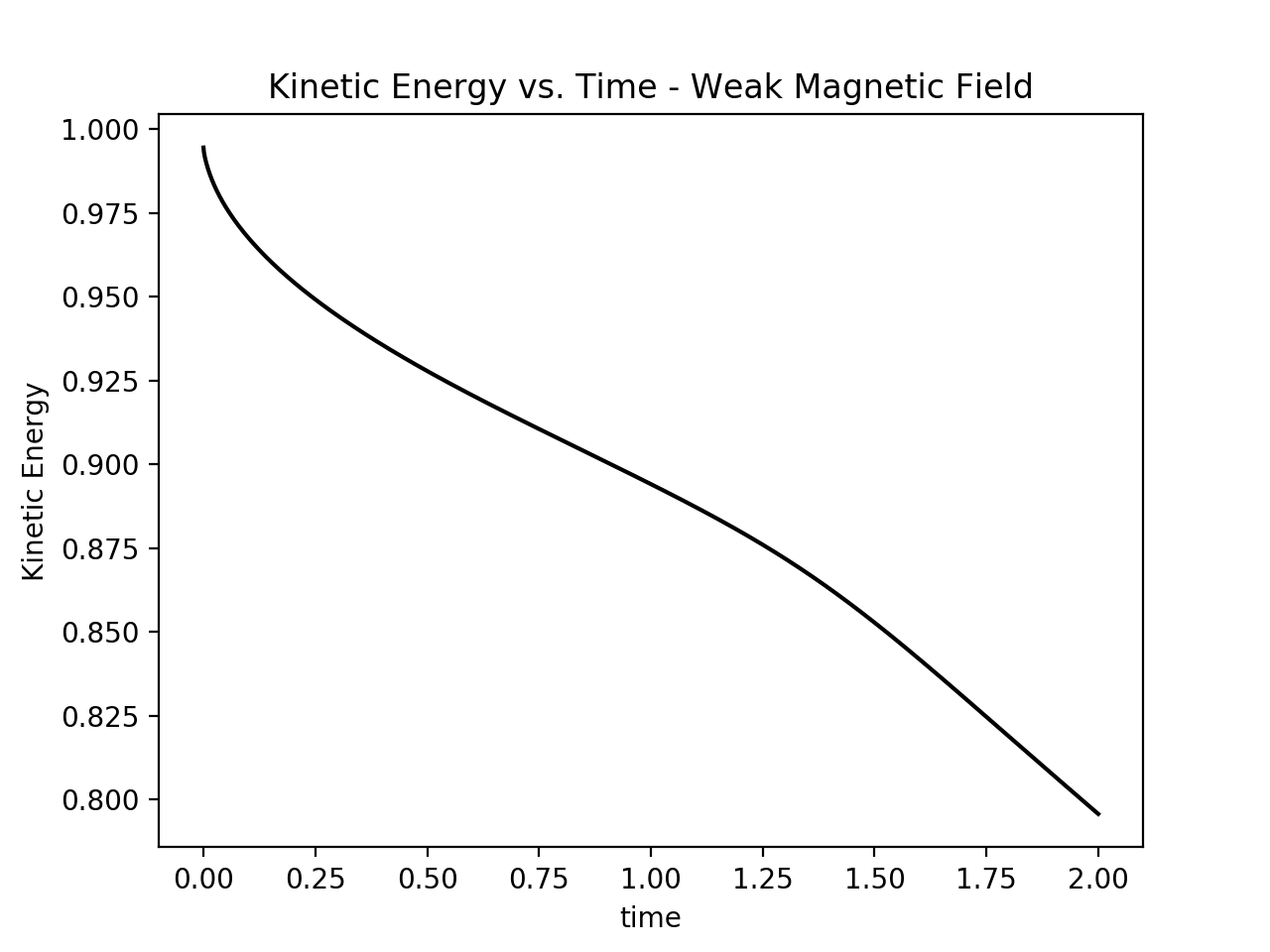}
    \includegraphics[width=0.31\textwidth]{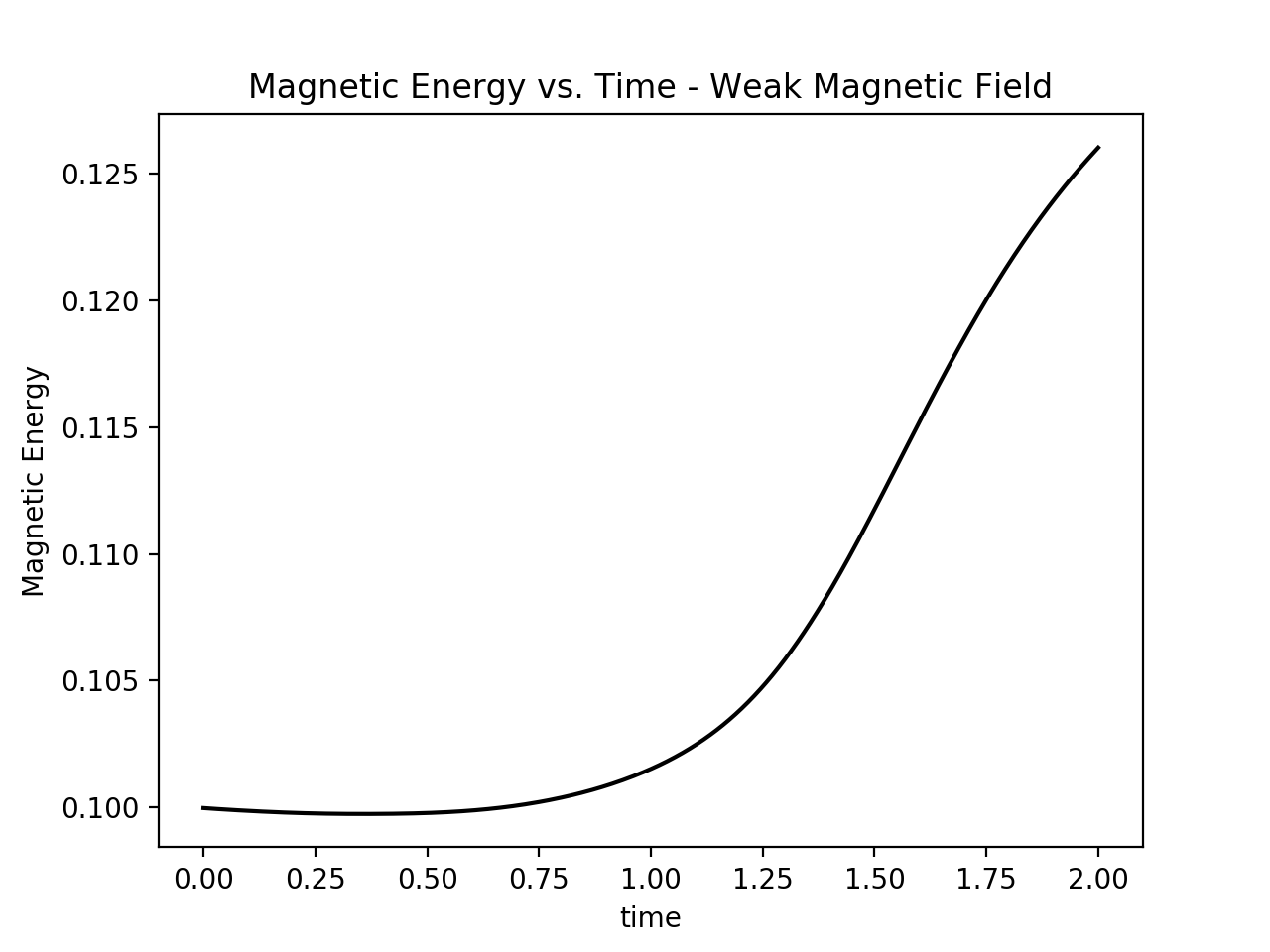}
    \includegraphics[width=0.31\textwidth]{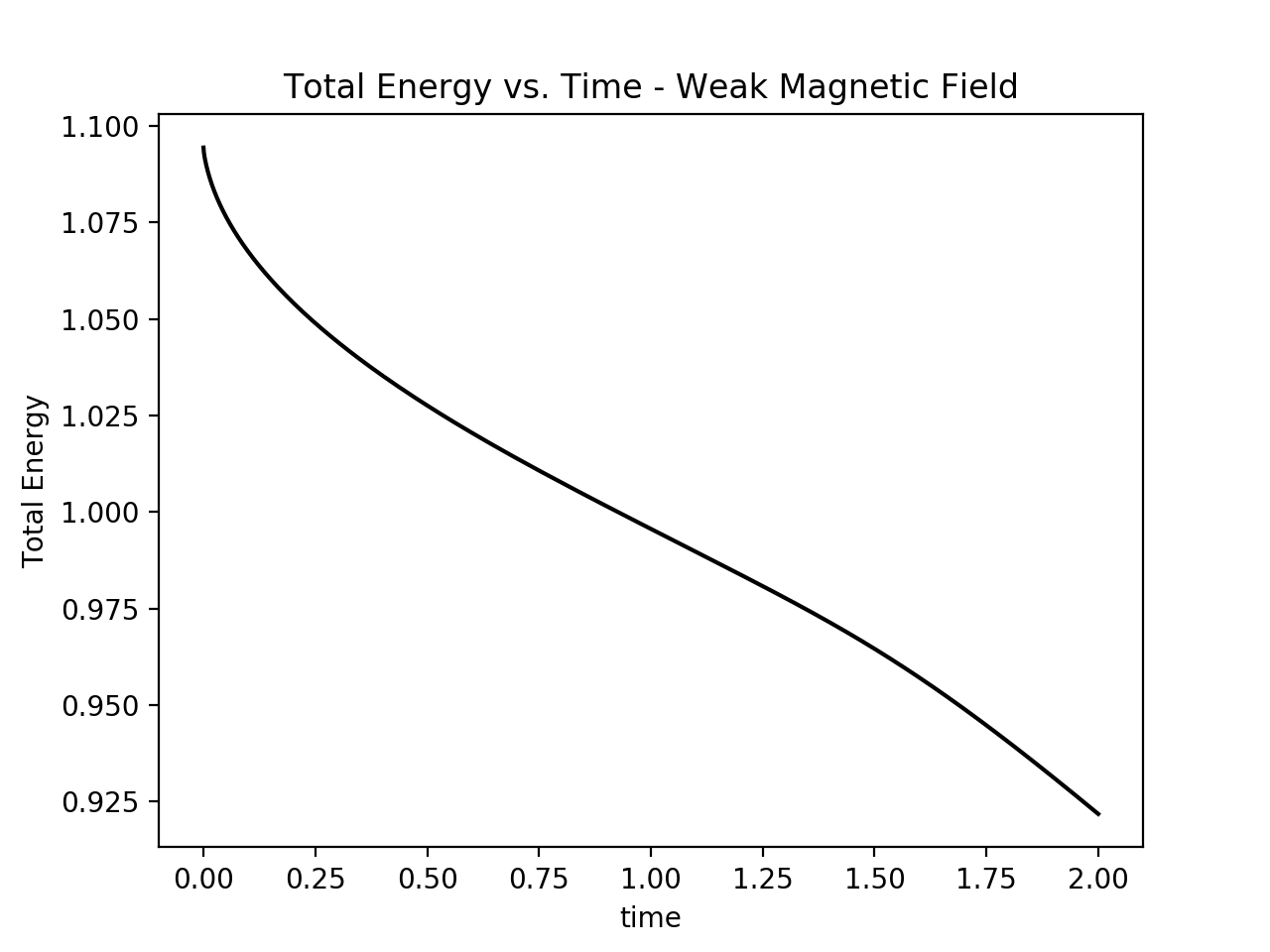}
    \caption{Time history of the kinetic (left), magnetic (middle), and total (right) energies for the Kelvin-Helmholtz instability problem for $M_A = 3$ and $Re=Re_m=10^3$.}
    \label{KHweakenergy}
\end{figure}

\begin{figure}[t!]
    \centering
    \includegraphics[width=0.31\textwidth]{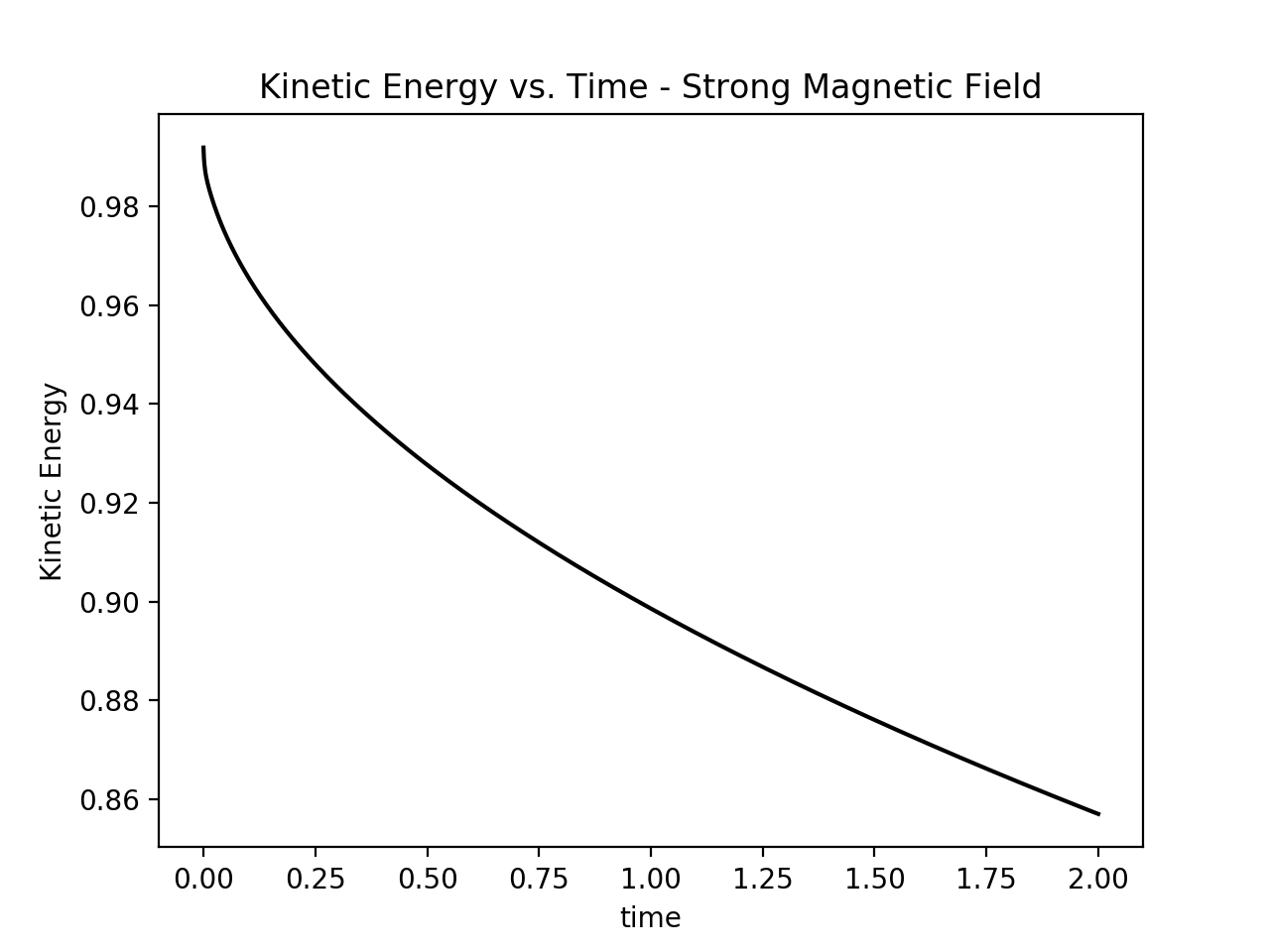}
    \includegraphics[width=0.31\textwidth]{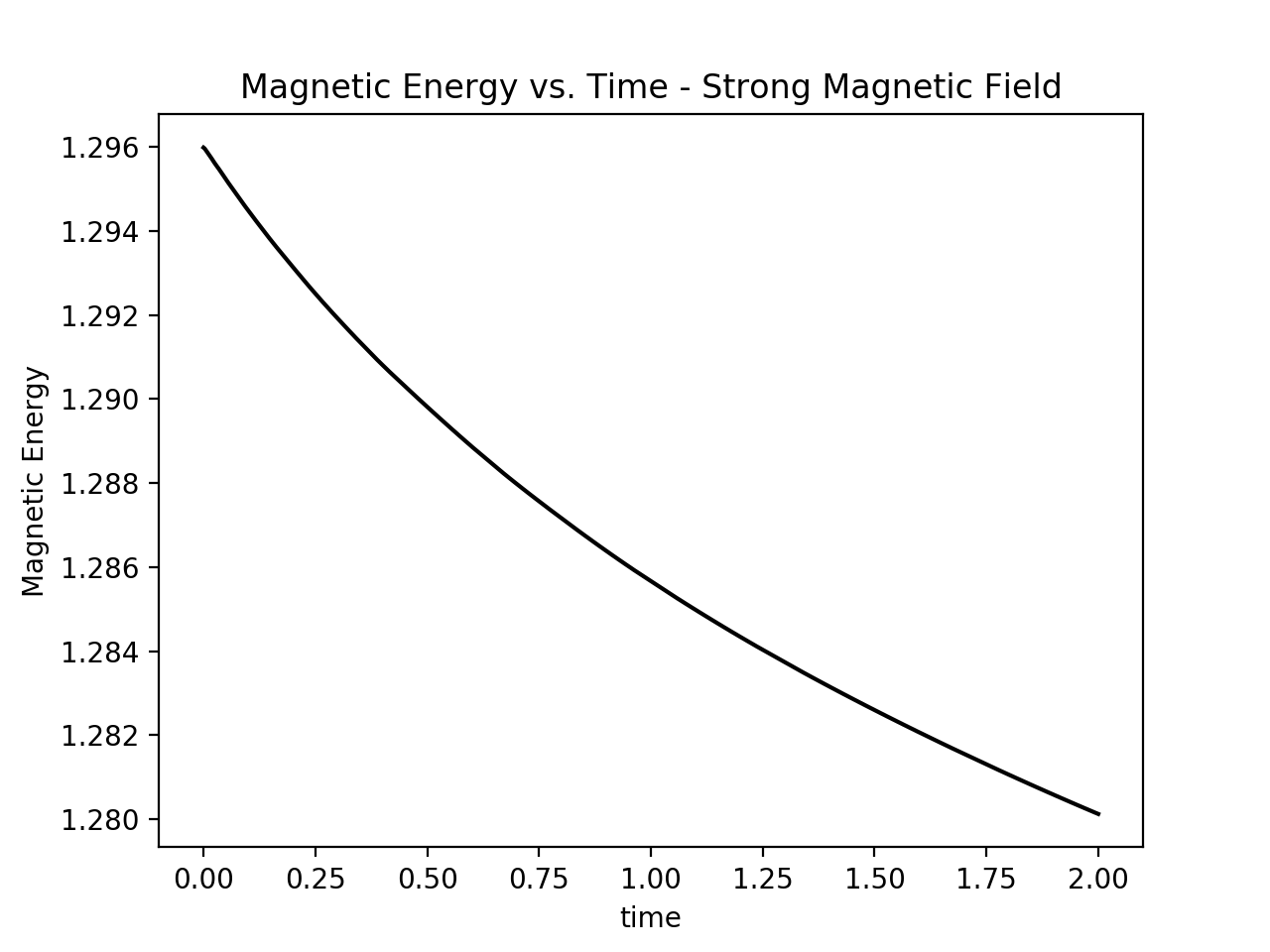}
    \includegraphics[width=0.31\textwidth]{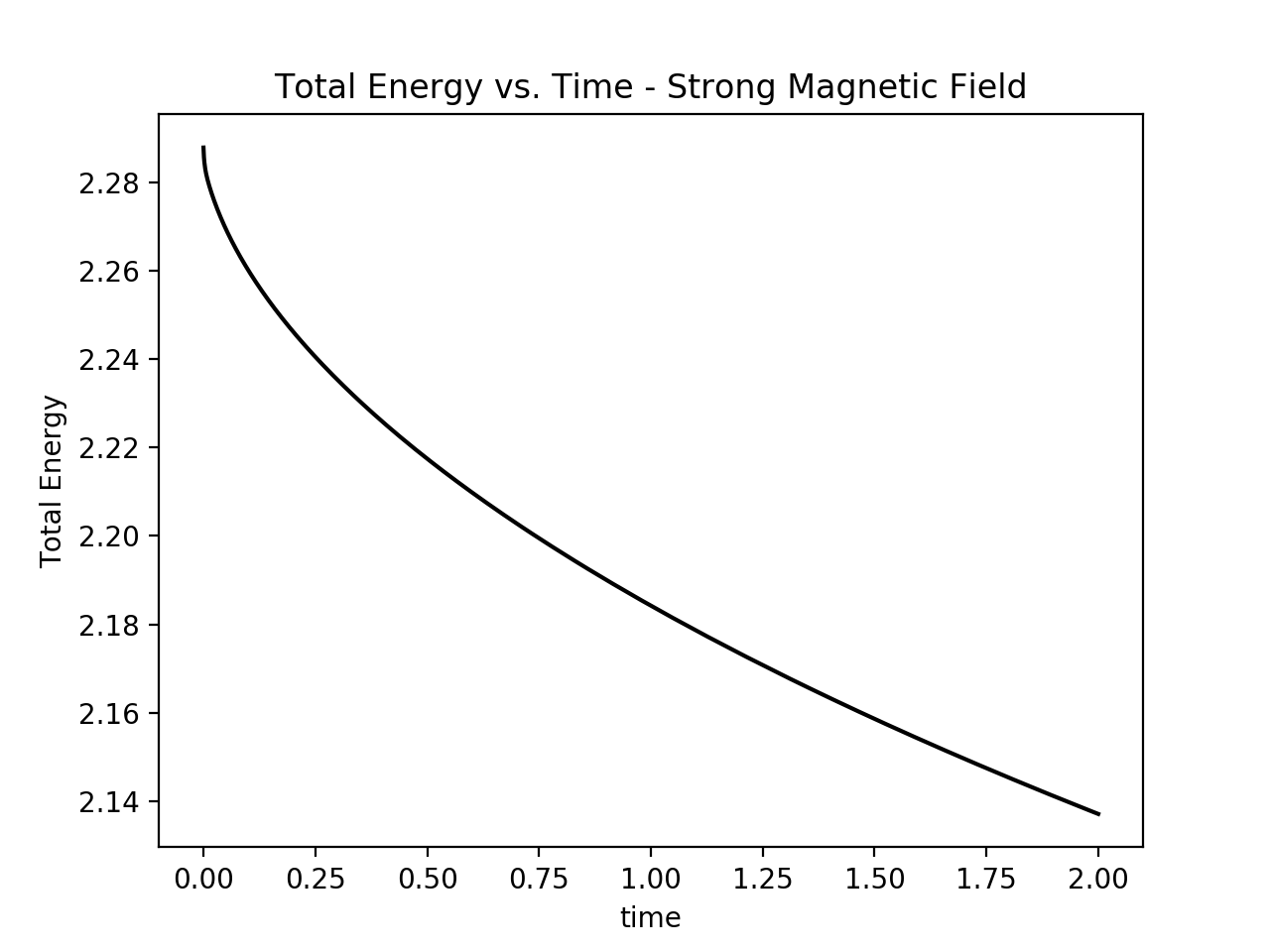}
    \caption{Time history of the kinetic (left), magnetic (middle), and total (right) energies for the Kelvin-Helmholtz instability problem for $M_A = \frac{5}{6}$ and $Re=Re_m=10^3$.}
    \label{KHstrongenergy}
\end{figure}

\begin{figure}[t!]
    \centering
    \includegraphics[width=0.31\textwidth]{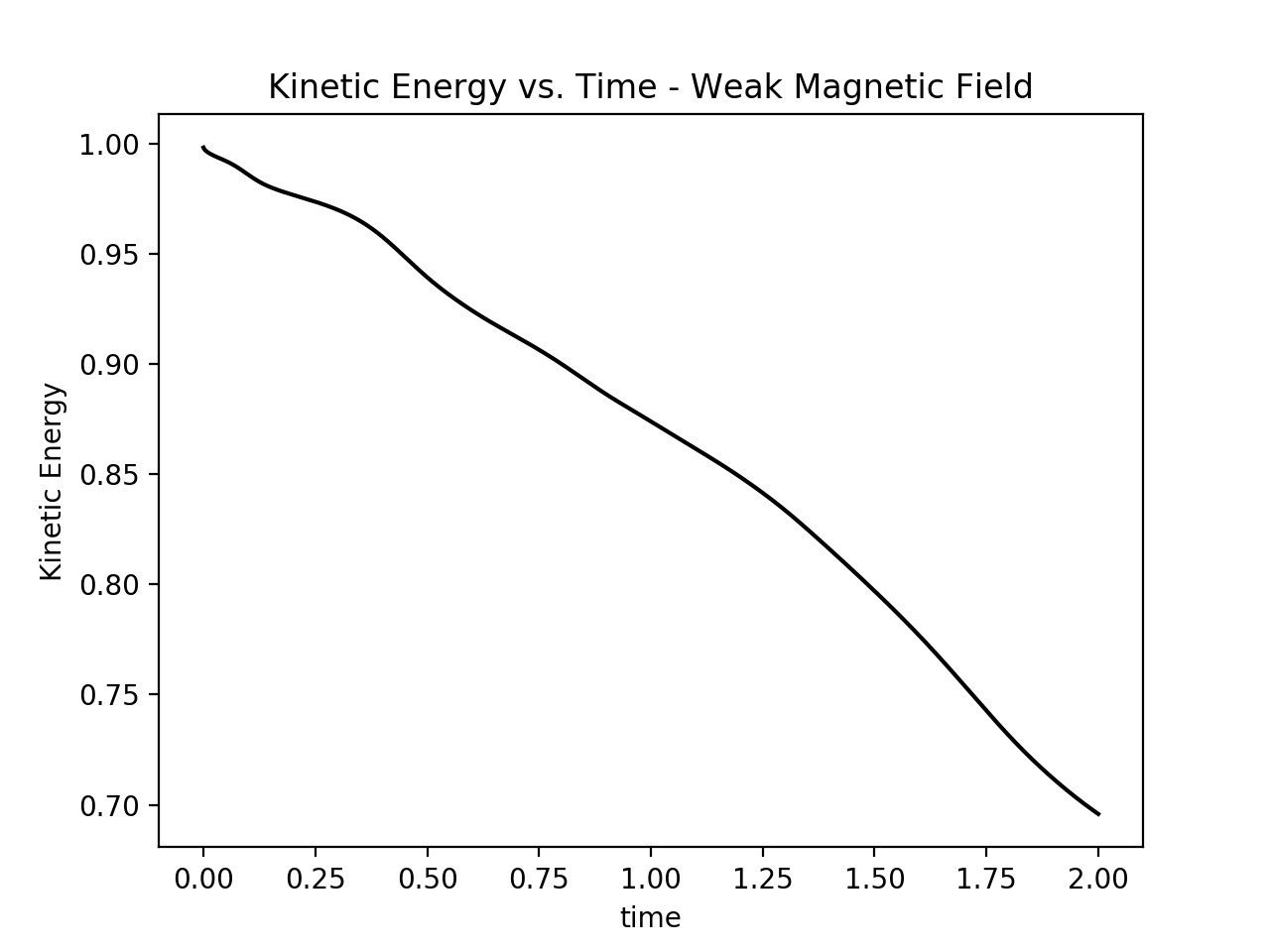}
    \includegraphics[width=0.31\textwidth]{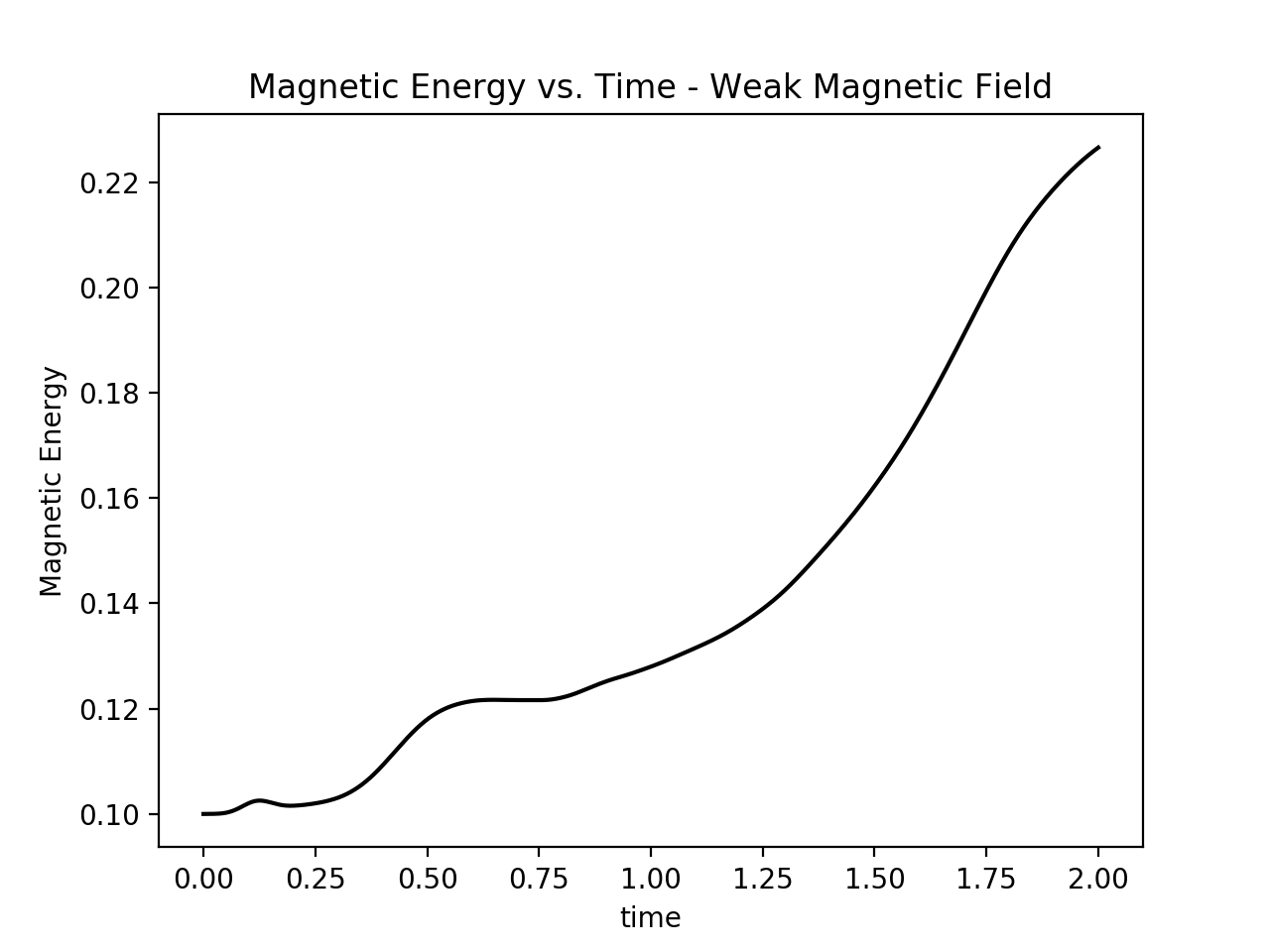}
    \includegraphics[width=0.31\textwidth]{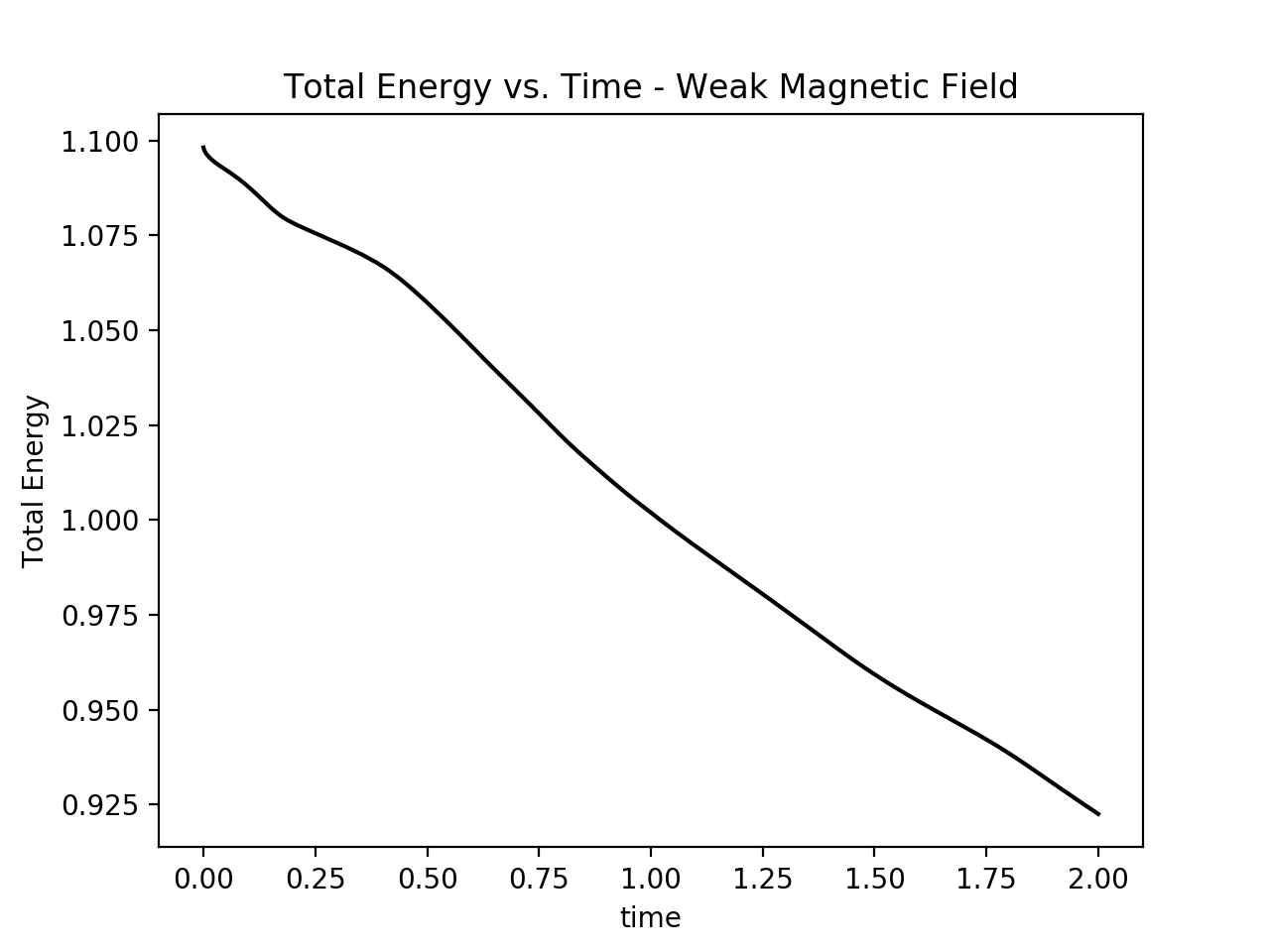}
    \caption{Time history of the kinetic (left), magnetic (middle), and total (right) energies for the Kelvin-Helmholtz instability problem for $M_A = 3$ and $Re=Re_m=10^4$.}
    \label{KHweakenergy4}
\end{figure}

\begin{figure}[t!]
    \centering
    \includegraphics[width=0.31\textwidth]{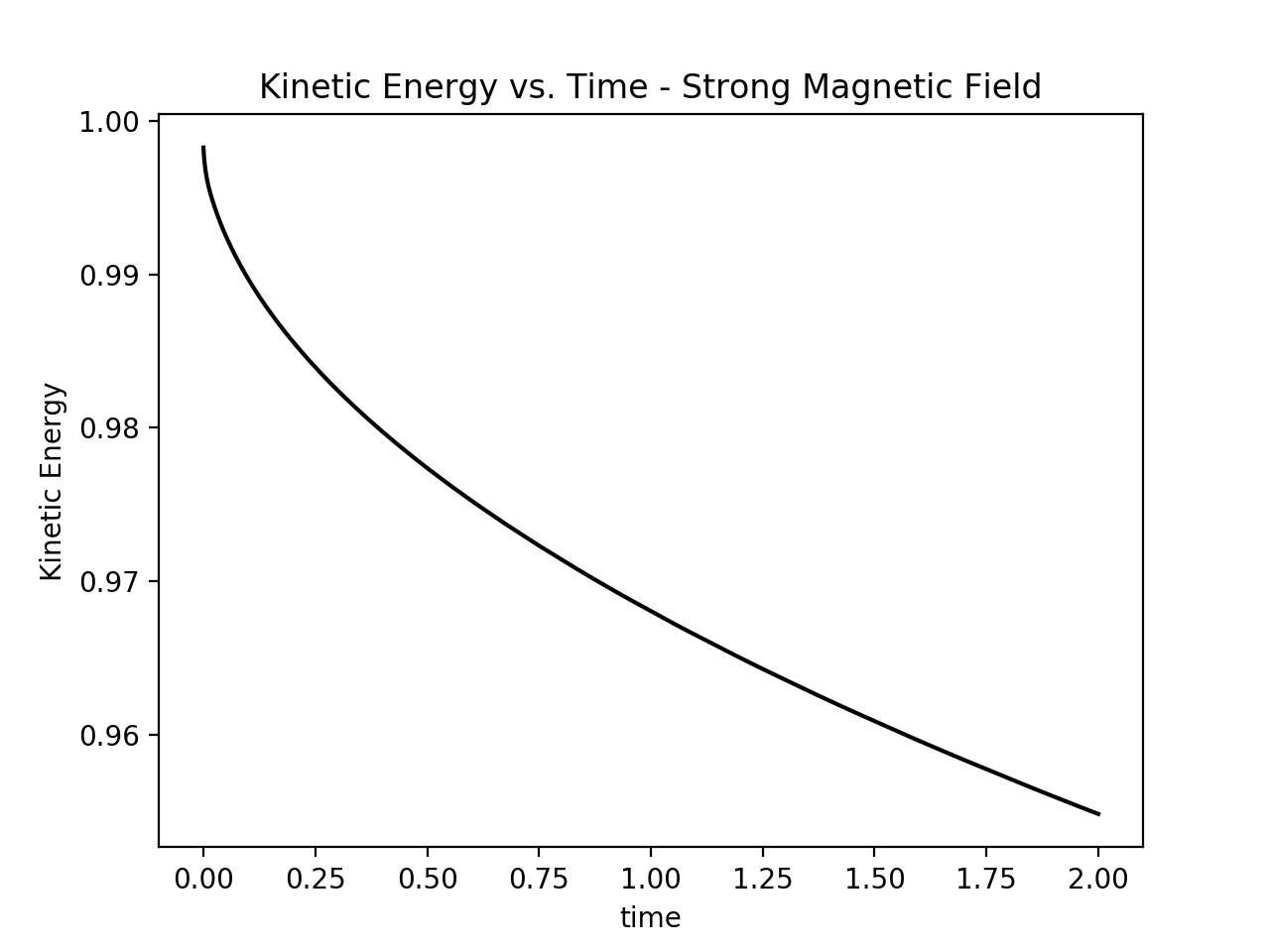}
    \includegraphics[width=0.31\textwidth]{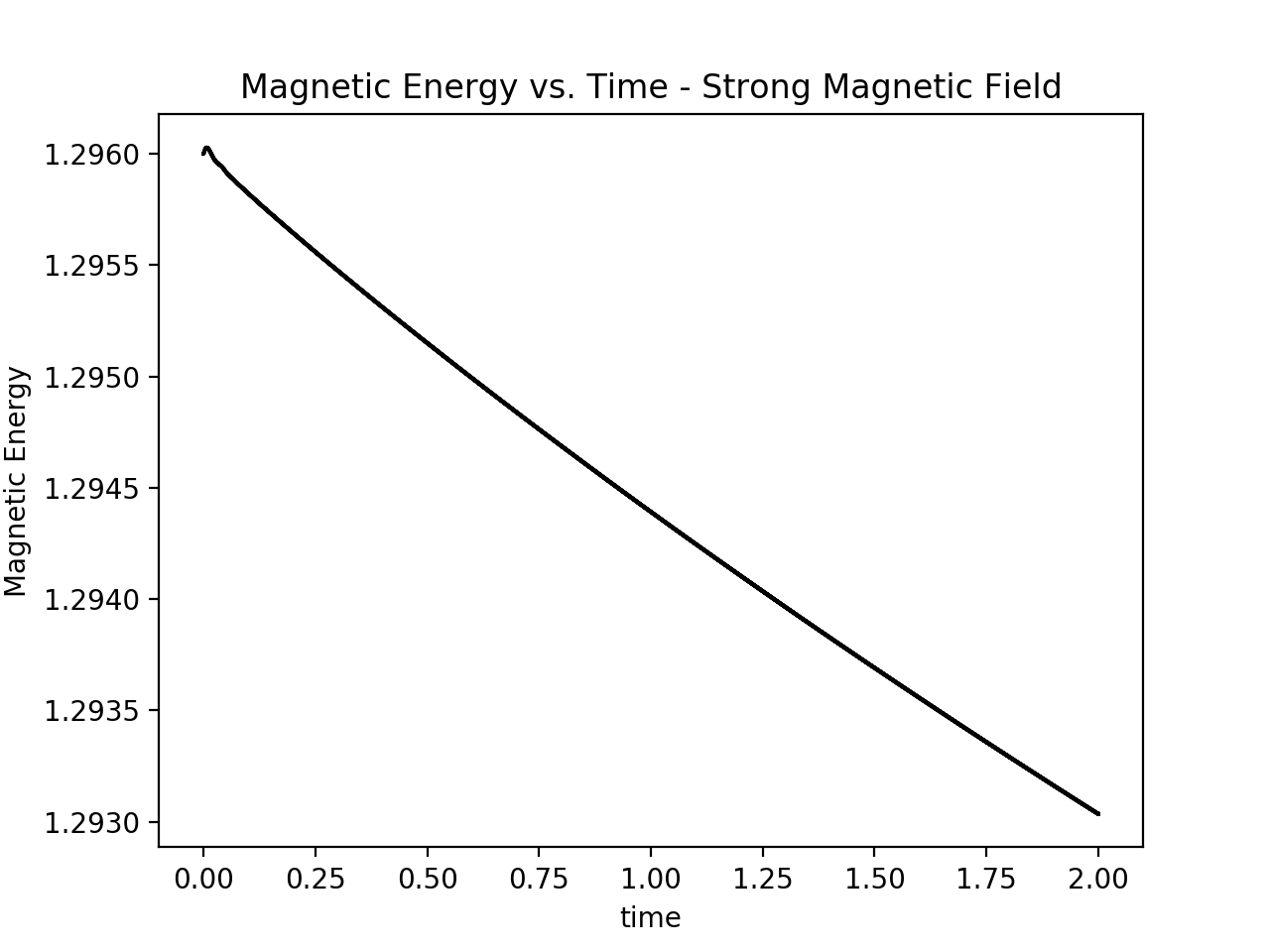}
    \includegraphics[width=0.31\textwidth]{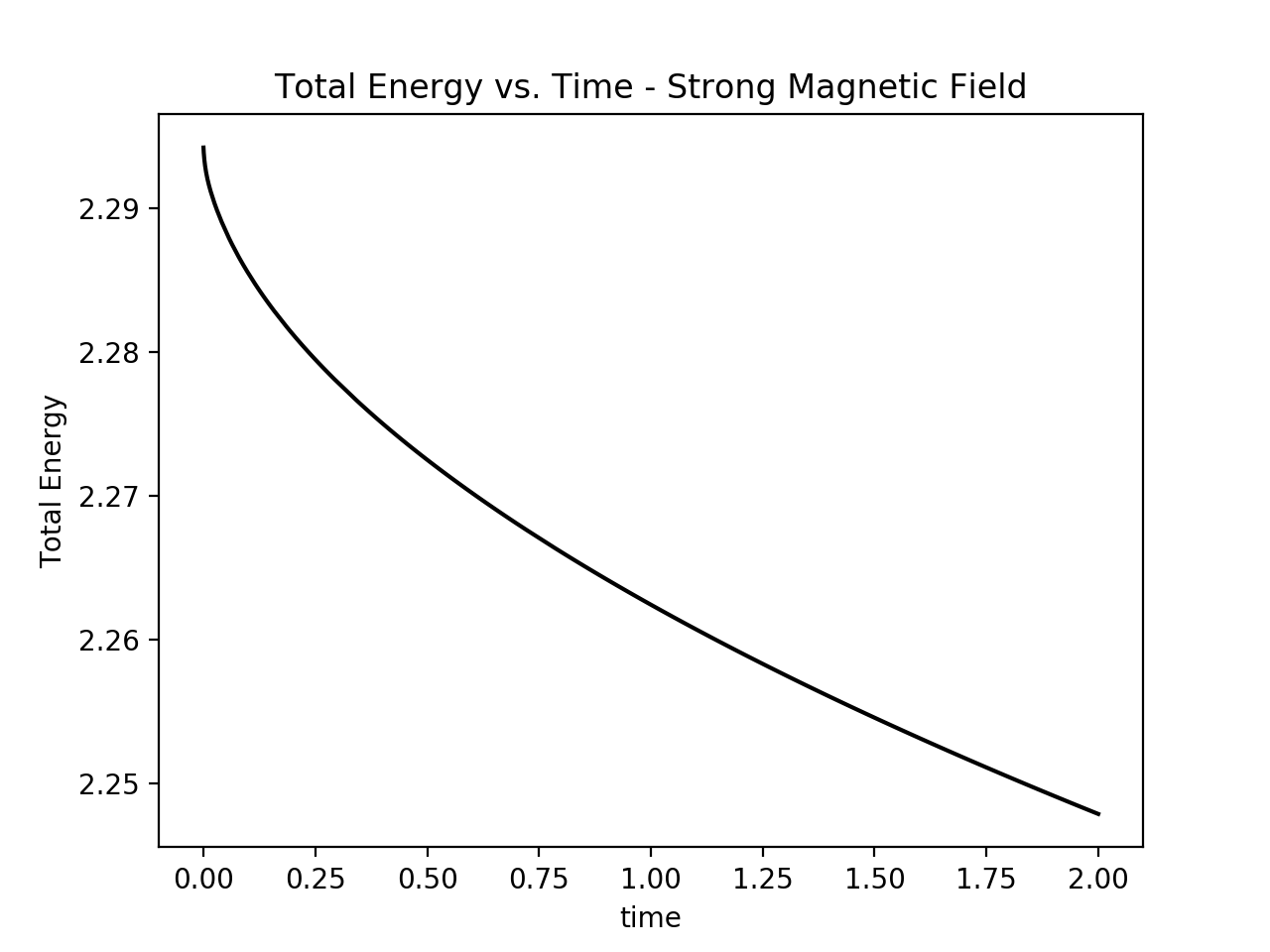}
    \caption{Time history of the kinetic (left), magnetic (middle), and total (right) energies for the Kelvin-Helmholtz instability problem for $M_A = \frac{5}{6}$ and $Re=Re_m=10^3$.}
    \label{KHstrongenergy4}
\end{figure}

In Fig. \ref{KHweakenergy}-\ref{KHstrongenergy4}, time histories of the kinetic, magnetic, and total energies from our simulations are displayed. Here, the kinetic energy is defined as $\int \frac{1}{2} \rho | \mathbf{u} |^2$, the magnetic energy is defined as $\int \frac{1}{2} \frac{1}{\mu_0} | \mathbf{B} |^2$, and the total energy is defined as the sum of kinetic energy and magnetic energy. Note that the total energy is non-increasing in time for each simulation, verifying our theoretical result that our method is energy stable. Both the kinetic energy and the magnetic energy are non-increasing in time as well for the two simulations at $M_A = \frac{5}{6}$ when the magnetic field is strong enough to stabilize the shear layer. However, there is accretion of magnetic energy for the two simulations at $M_A = 3$. In these simulations, energy is transferred from the kinetic energy to the magnetic energy using the Lorenz force $-\nabla\cdot\left(\frac{1}{\rho \mu_0}\left(\mathbf{B}\otimes\mathbf{B}-\frac{1}{2}|\mathbf{B}|^2\mathds{1}\right)\right)$ appearing in the conservation of momentum equation and the coupling term $-\nabla\cdot\left(\mathbf{B}\otimes\mathbf{u}\right)$ appearing in the magnetic induction equation.




\section{Conclusions}

In this paper, we introduced a hybridized discontinuous Galerkin method for the incompressible magnetohydrodynamics (MHD) equations. We first presented our semi-discrete method for the incompressible MHD equations and proved that this method is consistent, returns pointwise divergence-free velocity and magnetic fields, conserves momentum globally, and is energy stable. We then discretized in time using a second-order-in-time unconditionally stable generalized-$\alpha$ method and presented a block-iterative predictor-multicorrector scheme for solving the nonlinear algebraic system of equations at each time step, and we demonstrated how static condensation can be used to reduce the size of the matrix systems to be solved at each step in our predictor-multicorrector scheme. We confirmed our method yields optimal convergence rates for steady state problems using a manufactured solution as well as the Hartmann channel problem, and we confirmed our method is second-order-in-time using an Alfv\'en wave propagation problem. We also used the manufactured solution problem to verify our method is pressure robust. Finally, we demonstrated our method is energy stable using the Kelvin-Helmholtz instability problem.

In future work, there are several directions that could be explored. While the method and theoretical results appearing in this paper apply to both two- and three-dimensional incompressible MHD, it would be valuable to numerically confirm the spatial accuracy, temporal accuracy, pressure robustness, and energy stability of our method for three-dimensional incompressible MHD. It would also be useful to design robust and efficient solvers for the matrix systems attained at each step in our predictor-multicorrector scheme, perhaps building upon the work of \cite{rhebergen2018preconditioning}, as well as reduce the size of the velocity and magnetic field trace spaces using an embedded-hybridized DG method, as is done in \cite{rhebergen2020embedded}. Finally, the method presented in this paper could be extended to compressible MHD. The velocity field is not divergence-free in this setting, but preservation of the magnetic field divergence constraint is critical to the development of robust schemes for compressible MHD \cite{div-free}. Consequently, methods such as the elliptic projection method \cite{div-free}, eight and nine-wave formulations \cite{eight-wave, nine-wave}, the constrained transport method \cite{cs1, cs2}, locally divergence-free methods \cite{local-cdg, local-dg}, magnetic vector potential formulations \cite{A-form}, and Lagrange multiplier formulations \cite{mhd-hdg1, div-clean} have been used to enforce this condition. However, most of these methods still only enforce the magnetic field divergence constraint in an approximate sense as opposed to the method introduced in this paper.

\section{Acknowledgments}
This work was completed while the first two authors were undergraduate and graduate students respectively at the University of Colorado Boulder. The first author was supported by the Undergraduate Research Opportunities Program at the University of Colorado Boulder while the second and third authors were supported by the National Science Foundation under Grant CBET-1710670.

\bibliography{references}
\bibliographystyle{unsrt}

\end{document}